\documentclass[12pt]{amsart}

\usepackage{mathrsfs}
\usepackage{amsfonts}
\usepackage{amssymb}
\usepackage{amsfonts, amscd, amsmath, mathrsfs, amssymb, amsthm, amsxtra, bbding, epsfig, graphicx, latexsym, url, mathbbol, bbold}
\usepackage[papersize={8.3in,11.8in},textwidth=16.8cm,textheight=23cm,centering]{geometry}
\usepackage{enumerate}
 \usepackage{caption} 

\usepackage{graphicx}
\usepackage{tikz}
\usepackage{float}

\usepackage{xcolor}
\definecolor{cite}{rgb}{0.00,0.60,1.00}
\definecolor{url}{rgb}{1.00,0.10,0.80}
\definecolor{link}{rgb}{0.00,0.00,1.00}
\usepackage[colorlinks,linkcolor=link,urlcolor=url,citecolor=cite,pagebackref,breaklinks]{hyperref}

\hypersetup{
pdfstartpage=1,
pdfstartview=FitH}

\DeclareFontFamily{U}{mathx}{\hyphenchar\font45}
\DeclareFontShape{U}{mathx}{m}{n}{
      <5> <6> <7> <8> <9> <10>
      <10.95> <12> <14.4> <17.28> <20.74> <24.88>
      mathx10
      }{}
\DeclareSymbolFont{mathx}{U}{mathx}{m}{n}
\DeclareMathAccent{\widecheck}{\mathalpha}{mathx}{"71}

\numberwithin{equation}{section}

\allowdisplaybreaks

\newtheorem{theorem}{Theorem}[section]
\newtheorem{lemma}{Lemma}[section]

\newtheorem{definition}{Definition}[section]
\newtheorem{proposition}{Proposition}[section]

\newtheorem{conjecture}{Conjecture}[section]

\makeatletter
\newcounter{roem}
\renewcommand{\theroem}{\Roman{roem}}

\newcommand{\c@org@eq}{}
\let\c@org@eq\c@equation
\newcommand{\org@theeq}{}
\let\org@theeq\theequation

\newcommand{\setroem}{
\let\c@equation\c@roem
 \let\theequation\theroem}

\newcommand{\setarab}{
\let\c@equation\c@org@eq
\let\theequation\org@theeq}
\makeatother

\newtheorem*{claim*}{Claim}

\theoremstyle{theorem}
\newtheorem{remark}{\bf Remark}

\newcommand{\ud}{\mathrm{d}}
\newcommand{\ue}{\mathrm{e}}
\newcommand{\ft}{\mathrm{FT}}

\newcommand{\Li}{\mathrm{Li}}
\newcommand{\tr}{\mathrm{tr}}

\newcommand{\rank}{\mathrm{rank}}
\newcommand{\Swan}{\mathrm{Swan}}
\newcommand{\kl}{\mathrm{Kl}}
\newcommand{\Frob}{\mathrm{Frob}}

\newcommand{\GL}{\mathrm{GL}}
\newcommand{\SL}{\mathrm{SL}}

\DeclareMathOperator{\Mod}{mod}

\renewcommand{\bmod}[1]{\,(\Mod{ #1})}

\newcommand{\bn}{\mathbf{n}}

\newcommand{\balpha}{\boldsymbol{\alpha}}
\newcommand{\bbeta}{\boldsymbol{\beta}}
\newcommand{\bgamma}{\boldsymbol{\gamma}}
\newcommand{\bdelta}{\boldsymbol{\delta}}
\newcommand{\blambda}{\boldsymbol{\lambda}}

\newcommand{\bA}{\mathbf{A}}
\newcommand{\bC}{\mathbf{C}}
\newcommand{\bF}{\mathbf{F}}
\newcommand{\bN}{\mathbf{N}}
\newcommand{\bP}{\mathbf{P}}
\newcommand{\bQ}{\mathbf{Q}}
\newcommand{\bR}{\mathbf{R}}
\newcommand{\bZ}{\mathbf{Z}}

\newcommand{\sF}{\mathscr{F}}
\newcommand{\sX}{\mathscr{X}}

\newcommand{\cA}{\mathcal{A}}
\newcommand{\cB}{\mathcal{B}}
\newcommand{\cC}{\mathcal{C}}
\newcommand{\cD}{\mathcal{D}}
\newcommand{\cE}{\mathcal{E}}
\newcommand{\cF}{\mathcal{F}}
\newcommand{\cG}{\mathcal{G}}

\newcommand{\cK}{\mathcal{K}}
\newcommand{\cL}{\mathcal{L}}
\newcommand{\cM}{\mathcal{M}}
\newcommand{\cN}{\mathcal{N}}
\newcommand{\cP}{\mathcal{P}}
\newcommand{\cQ}{\mathcal{Q}}
\newcommand{\cR}{\mathcal{R}}

\newcommand{\cT}{\mathcal{T}}

\newcommand{\cZ}{\mathcal{Z}}

\newcommand{\fc}{\mathfrak{c}}

\newcommand{\fI}{\mathfrak{I}}

\newcommand{\fP}{\mathfrak{P}}
\newcommand{\fS}{\mathfrak{S}}

\def\le{\leqslant}
\def\leq{\leqslant}
\def\ge{\geqslant}
\def\geq{\geqslant}

\usepackage{graphicx}
\usepackage{tikz}

\begin{document}

\title{On the Brun--Titchmarsh theorem. I}

\author{Ping Xi}

\address{School of Mathematics and Statistics, Xi'an Jiaotong University, Xi'an 710049, P. R. China}
\address{Institute of Pure Mathematics, Xi'an Jiaotong University, Xi'an 710049, P. R. China}

\email{pingxi.cn@gmail.com, ping.xi@xjtu.edu.cn}

\author{Junren Zheng}

\address{School of Mathematics and Statistics, Xi'an Jiaotong University, Xi'an 710049, P. R. China}

\address{D{\'e}partment  de Math{\'e}matiques et Statistique,
Universit{\'e} de Montr{\'e}al, CP 6128 succ Centre-Ville, Montr{\'e}al, QC  H3C 3J7, Canada}
\email{junrenzheng03@gmail.com, junren@stu.xjtu.edu.cn}



\begin{abstract} 
The classical Brun--Titchmarsh theorem gives an upper bound, which is of correct order of magnitude in the full range, for the number of primes $p\leqslant x$ satisfying $p\equiv a\bmod q$. We strengthen this inequality for different ranges of $\log q/\log x$, improving upon previous works by Motohashi, Goldfeld, Iwaniec, Friedlander and Iwaniec, and Maynard for general or special moduli. In particular, we are able to beat Iwaniec's barrier $q<x^{9/20-}$, and improve all existing inequalities in the range $x^{9/20}\ll q<x^{1/2-}$ by utilizing bilinear or trilinear structures in the remainder terms of linear sieve.
The proof is based on various estimates for character and exponential sums, which we derive by appealing to arithmetic exponent pairs and bilinear forms with algebraic trace functions from $\ell$-adic cohomology, trilinear forms with Kloosterman fractions, and sums of Kloosterman sums from spectral theory of automorphic forms, as well as large value theorem for Dirichlet polynomials.
\end{abstract}

\dedicatory{{\it \small Professor Chengdong PAN in memoriam}}
\maketitle

\setcounter{tocdepth}{1}


\section{Backgrounds and main results}
\label{sec:backgrounds&results}

\subsection{Main results}
This paper aims to understand the upper bound for the counting function
\begin{align*}
\pi(x; q, a):=|\{p\leqslant x:p\equiv a\bmod q\}|,
\end{align*}
where $q\in\bZ^+$, $(a,q)=1$ and $x\geqslant3$. It is widely expected that 
\begin{align}\label{eq:pi(x;q,a)-approximation}
\pi(x ; q, a)=(1+o(1)) \frac{1}{\varphi(q)}\frac{x}{\log x}
\end{align}
as $q,x\rightarrow+\infty$ in a very broad range, and the classical Siegel--Walfisz theorem asserts that this is valid for all $q\leqslant (\log x)^A$ with large but fixed $A$. No asymptotic formula is known for larger $q$ compared to $x$. Note that the Generalized Riemann Hypothesis may guarantee \eqref{eq:pi(x;q,a)-approximation} for all $q\le x^{1/2-\varepsilon}$.

With the aid of Brun's sieve, Titchmarsh \cite{Ti30} proved that
\begin{align}\label{eq:BT-Titchmarsh}
\max_{(a,q)=1}\pi(x;q,a)\ll \frac{x}{\varphi(q)\log{(x/q)}}
\end{align}
for all $1\leqslant q<x,$ and this bears the name {\it Brun--Titchmarsh theorem}. 
In what follows, we assume $q\sim x^\varpi$ with $\varpi\in[0,1[$, and seek a positive constant $C(\varpi)$, as small as possible, such that
\begin{align}\label{eq:BT}
\max_{(a,q)=1}\pi(x;q,a)\leqslant (C(\varpi)+\varepsilon)\frac{1}{\varphi(q)}\frac{x}{\log x}
\end{align}
holds for any $\varepsilon>0$ and $x>x_0(\varepsilon).$ Note that \eqref{eq:BT-Titchmarsh} guarantees the existence of $C(\varpi)$ for all $\varpi\in[0,1[$.

In this paper, we establish several choices of $C(\varpi)$ for general moduli and smooth moduli $q\sim x^\varpi$.
We now quickly state our first result. 

\begin{theorem}\label{thm:BT-generalmoduli<1/2}
Let $\varpi\in[9/20,1/2[.$ Then we may take
\begin{align}\label{eq:C(varpi):varpi<1/2}
C(\varpi)&=\frac{66}{33-16\varpi}-C^*(\varpi)
\end{align}
in $\eqref{eq:BT}$ with some constant $C^*(\varpi)\geqslant0.$ In general, we may choose
\begin{equation}\label{eq:C*(varpi):varpi<1/2}
\begin{split}
C^*(\varpi)&=-\frac{66}{33-16\varpi}\int_2^4\frac{\log(t-1)}{t}\ud t+\frac{8}{4-(3+\vartheta)\varpi}\int_{\frac{8-(7+2\vartheta)\varpi}{4\varpi}}^{\frac{165(4-(3+\vartheta)\varpi)}{4(33-16\varpi)}-\frac{1}{4}}\frac{\log(t-1)}{t}\ud t\\
&\ \ \ \ +\frac{8}{4-(1+\vartheta)\varpi}\int_{\max\{\frac{4-(1+\vartheta)\varpi}{2(2-3\varpi)}-\frac{5}{4},~2\}}^{\frac{8-(7+2\vartheta)\varpi}{4\varpi}}\frac{\log(t-1)}{t}\ud t,
\end{split}
\end{equation}
and for prime $q$ we may choose
\begin{equation}\label{eq:C*(varpi):varpi<1/2,primemoduli}
\begin{split}
C^*(\varpi)&=-\frac{66}{33-16\varpi}\int_2^4\frac{\log(t-1)}{t}\ud t+\frac{16}{8-7\varpi}\int_{\frac{8-7\varpi}{4\varpi}}^{\frac{165(8-7\varpi)}{8(33-16\varpi)}}\frac{\log(t-1)}{t}\ud t\\
&\ \ \ \ +\frac{16}{8-3\varpi}\int_{\max\{\frac{9\varpi}{4(2-3\varpi)},~2\}}^{\frac{8-7\varpi}{4\varpi}}\frac{\log(t-1)}{t}\ud t.
\end{split}
\end{equation}
Here $\vartheta$ denotes the exponent towards the Ramanujan--Petersson conjecture  for $\GL_2(\bQ).$ 
In particular, one may take $\vartheta=7/64$ thanks to Kim and Sarnak {\rm\cite{KS03}.}
\end{theorem}

It is worthwhile to make a remark on the merit of Theorem \ref{thm:BT-generalmoduli<1/2}. Note that the choices \eqref{eq:C*(varpi):varpi<1/2} and \eqref{eq:C*(varpi):varpi<1/2,primemoduli} for $C^*(\varpi)$ are always positive for all $\varpi\in[9/20,1/2[$, and Theorem \ref{thm:BT-generalmoduli<1/2} improves all previously known bounds even if taking $C^*(\varpi)=0$ in \eqref{eq:C(varpi):varpi<1/2}. The gain coming from positive choices of $C^*(\varpi)$ arises from a suitable Buchstab iteration in the application of sieves, so that one may create {\it trilinear} structures in the remainder terms with the aid of the well-factorization of the linear sieve weight constructed by Iwaniec \cite{Iw80}.
For general moduli, we estimate a twisted moment of character sums (see Lemma \ref{lm:weightedL-secondmoment-specialcoefficient}) by appealing to Kloostermania developed by Deshouillers and Iwaniec \cite{DI82}; and in the prime moduli case, the above choice for $C(\varpi)$ is independent of the Ramanujan--Petersson conjecture, and it will follow from another different approach, by appealing to a recent work of Blomer, Humphries, Khan and Milinovich \cite{BHKM20} on twisted fourth moment of Dirichlet $L$-functions to prime moduli
(see Lemma \ref{lm:weightedL-fourthmoment} below).

In fact, 
Motohashi \cite{Mo73} took $C(\varpi)={16}/(8-3\varpi)$ for $\varpi<1/3,$ and the admissible range was later extended by Goldfeld \cite{Go75} to $\varpi<24/71.$ Iwaniec \cite{Iw82} obtained the same constant for all $\varpi<9/20,$ for which he utilized Burgess's estimate for short character sums:
\begin{align}\label{eq:Burgesscharactersum}
\sum_{n\leqslant N}\chi(n)\ll N^{1-\frac{1}{r}}q^{\frac{r+1}{4r^2}+\varepsilon}
\end{align}
holds for any non-trivial Dirichlet character modulo $q$ and $r\in\{1,2,3\}$ (see \cite{Bu63,Bu86}).
Taking $r=2$, \eqref{eq:Burgesscharactersum} implies the so-called Burgess bound $L(\frac{1}{2},\chi)\ll q^{3/16+\varepsilon}$ for Dirichlet $L$-functions.
The restriction $\varpi<9/20$ is quite essential in Iwaniec's approach. 
On the other hand, the choice $C(\varpi)=4/(2-\varpi)$ for $\varpi\in[9/20,1/2[$ has been obtained implicitly by Iwaniec \cite{Iw82}. This has been illustrated by Baker \cite{Ba96}.
One may compare our bounds in Theorem \ref{thm:BT-generalmoduli<1/2} with Iwaniec--Baker's constant. In particular, we may take $C(\frac{9}{20})\approx 2.4814$ for general $q$ and $C(\frac{9}{20})\approx 2.441$
for prime $q$, whereas Iwaniec--Baker had to take $C(\frac{9}{20})\approx 2.5806.$

For $\varpi\geqslant1/2$, we have the following choices for $C(\varpi)$.

\begin{theorem}\label{thm:BT-generalmoduli>1/2}
Let $q$ be a prime. Then we may take $C(\varpi)$ in $\eqref{eq:BT}$ as follows:
\begin{align*}
C(\varpi)=
\begin{cases}
8/(5-5\varpi), \quad &\varpi\in[1/2,~12/23[,\\
32/(32-43\varpi), \quad &\varpi\in[1/2,~32/61[,\\
24/(16-17\varpi), &\varpi\in[32/61,~8/15[,\\
48/(40-49\varpi), &\varpi\in[8/15,~7/13[,\\
16/(11-12\varpi), & \varpi\in[7/13,~6/11[,\\
32/(28-35\varpi), &\varpi\in[ 6/11,~4/7[. 
\end{cases}
\end{align*}
\end{theorem}

To see the strength of Theorems \ref{thm:BT-generalmoduli<1/2} and \ref{thm:BT-generalmoduli>1/2}, we recall some previous admissible values of $C(\varpi)$ obtained by a series of big names.
\begin{table}[htbp]
\renewcommand{\arraystretch}{1.1}
\renewcommand{\tabcolsep}{1.7mm}
\begin{tabular}{ | c | c | c | }\hline
$\varpi$ & Admissible  $C(\varpi)$ & Authors  \\ \hline
$]0,1[$ & existence & Titchmarsh \cite{Ti30} \\ \hline 
$]0,1[$ & $2/(1-\varpi)$ & van Lint \& Richert \cite{LR65} \\
&  & Montgomery \& Vaughan \cite{MV73}\\
&  & Selberg \cite{Se91}\\ \hline 
$]0,1/3[$ & $16/(8-3\varpi)$ & Motohashi \cite{Mo73} \\ \hline 
$[1/3,2/5]$ & $4/(2-\varpi)$ & Motohashi \cite{Mo74} \\ \hline 
$]2/5,1/2]$ & $2/(2-3\varpi)$ & Motohashi \cite{Mo73} \\ \hline 
$]0,24/71[$ & $16/(8-3\varpi)$ & Goldfeld \cite{Go75} \\ \hline 
$]0,9/20[$ & $16/(8-3\varpi)$ & Iwaniec \cite[Theorem 3]{Iw82} \\ \hline 
$[9/20,2/3]$ & $8/(6-7\varpi)$ & Iwaniec \cite[Theorem 6]{Iw82} \\ \hline 
$[9/20,1/2[$ & $4/(2-\varpi)$ & Baker \cite[Theorem 1]{Ba96} \\ \hline 
$[6/11,1[$ & $(2-(\frac{1-\varpi}{4})^6)/(1-\varpi)$ & Friedlander \& Iwaniec \cite[Theorem 1]{FI97} \\ \hline 
$[1-\delta,1[$ & $(2-c_0(1-\varpi)^2)/(1-\varpi)$ & Bourgain \& Garaev \cite[Theorem 4]{BG14} \\ \hline 
$]0,1/8[$ & $2$ & Maynard \cite{Ma13} \\ \hline 
\end{tabular}
\caption{Admissible values of $C(\varpi)$} \label{table:BTconstants}
\end{table}

The second aim of this paper  is to reduce the admissible value of $C(\varpi)$ further by assuming that $q$ has good factorizations.
For a small constant $\eta>0$, we say $q$ is {\it $\eta$-smooth} if all prime factors of $q$ do not exceed $q^\eta.$
The smoothness assumption allows us to utilize the factorization of $q$ to control the resultant algebraic exponential sums very effectively. In particular, one may go far beyond the so-called P\'olya--Vinogradov barrier for short exponential and character sums, as a consequence of the general form of the $q$-analogue of van der Corput method and arithmetic exponent pairs developed by Wu, Xi and Sawin \cite{WX21}. We refer the reader to Section \ref{subsec:exponentpairs} for details.

The following theorem shows that we may go further when $\varpi$ is close to $2/3$, and in particular for explicit ranges with $\varpi>2/3$, provided that $q$ is squarefree and sufficiently smooth.

\begin{theorem}\label{thm:BT-smoothmoduli2/3specialexponentpair}
Suppose $q$ is $\eta$-smooth and squarefree. Then we may take
\begin{align*}
C(\varpi)=\frac{24}{15-16\varpi}
\end{align*}
in $\eqref{eq:BT}$ for all $\varpi\in[3/10,3/4],$
provided that $\eta$ is small enough.
\end{theorem}

Theorem \ref{thm:BT-smoothmoduli2/3specialexponentpair} is in fact a special case of the following general formulation with $(\kappa,\lambda)=(1/6,2/3)$.

\begin{theorem}\label{thm:BT-smoothmoduli2/3generalexponentpair}
Suppose $q$ is $\eta$-smooth and squarefree. For each exponent pair $(\kappa,\lambda,*)$ as defined by Lemma $\ref{lm:arithmeticexponentpairs}$ and $\eqref{eq:arithmeticexponentpairs-set},$  we may take
\begin{align*}
C(\varpi)=\frac{4}{(3+\kappa-\lambda)-(3+2\kappa-\lambda)\varpi}
\end{align*}
in $\eqref{eq:BT}$ for all $\varpi$ with
\begin{align*}
\frac{1+\kappa-\lambda}{2+2\kappa-\lambda} 
\leq \varpi \leq \frac{1+\kappa-\lambda}{1+2\kappa-\lambda},
\end{align*}
provided that $\eta$ is small enough.
\end{theorem}

Theorems \ref{thm:BT-smoothmoduli2/3specialexponentpair} and \ref{thm:BT-smoothmoduli2/3generalexponentpair} are stronger than the results of van Lint and Richert \cite{LR65} and that of Iwaniec \cite{Iw82} when $\varpi$ is around $2/3$, assuming that $q$ is $\eta$-smooth and squarefree. 
Put
\begin{align*}
f(\kappa,\lambda)
=\frac{1+\kappa-\lambda}{1+2\kappa-\lambda},\ \ \
g_\varpi(\kappa,\lambda)
=(3+\kappa-\lambda)-(3+2\kappa-\lambda)\varpi.
\end{align*}
To have $f(\kappa,\lambda)$ or $g_\varpi(\kappa,\lambda)$ as large as possible, one needs to choose $(\kappa,\lambda)$ by applying $A$- and $B$-processes to $(0,1)$ suitably; we refer to Section \ref{subsec:exponentpairs} for the terminology and conventions.
The optimization process can be conducted following the algorithm in \cite[Section 5.3]{GK91}.
In particular, for $\varpi=2/3,$ we find
\begin{align*}
g_\varpi(\kappa,\lambda)
&=1-\frac{\kappa+\lambda}{3},
\end{align*}
and it suffices to minimize $\kappa+\lambda.$ This task matches exactly the problem of determining Rankin's constant for Riemann zeta function.
According to Rankin (see also \cite[Section 5.4]{GK91}), the minimal value for $\kappa+\lambda$
should be $0.829021\cdots$ by choosing
\begin{align*}
(\kappa,\lambda)=ABA^3BA^2BABABA^2BABABA^2BA^2\cdots B(0,1),
\end{align*}
 in which case we get $g_\varpi(\kappa,\lambda)=0.723659\cdots$.
This yields $C(2/3)=4/g_\varpi(\kappa,\lambda)\approx5.52746$ in Theorem \ref{thm:BT-smoothmoduli2/3generalexponentpair}. Note that van Lint and Richert \cite{LR65} obtained the constant $2/(1-\varpi)=6$ for $\varpi=2/3.$

Recall the choice of $C(\varpi)$ in Theorem \ref{thm:BT-smoothmoduli2/3generalexponentpair} and Iwaniec's choice $8/(6-7\varpi)$ in \cite{Iw82}.
We expect that
\begin{align*}
\frac{4}{(3+\kappa-\lambda)-(3+2\kappa-\lambda)\varpi}<\frac{8}{6-7\varpi},
\end{align*}
which is equivalent to
\begin{align*}
\varpi>\frac{2(\lambda-\kappa)}{1+2\lambda-4\kappa}.
\end{align*}
Taking
\begin{align*}
BA^{k-1}B\cdot (0,1)&=\Big(\frac{1}{2}-\frac{k}{2^{k+1}-2},~\frac{1}{2}+\frac{1}{2^{k+1}-2}\Big)
\end{align*}
with sufficiently large $k,$ we find Theorem \ref{thm:BT-smoothmoduli2/3generalexponentpair} improves 
Iwaniec's choice $8/(6-7\varpi)$ as long as $\varpi>1/2+\varepsilon$, and $q$ is $\eta$-smooth and squarefree.

The following theorem gives alternative choices for $C(\varpi)$ when $q$ is smooth and $\varpi<1/2.$
\begin{theorem}\label{thm:BT-smoothmoduli<1/2}
Suppose $q$ is $\eta$-smooth and squarefree. Then we may take 
\begin{align*}
C(\varpi)=
\begin{cases}
2, \quad &\varpi\in[1/8,~5/12[,\\
5/(5-6\varpi), &\varpi\in[5/12,~9/20[
\end{cases}
\end{align*}
in $\eqref{eq:BT},$ provided that $\eta$ is small enough.
\end{theorem}

Theorem \ref{thm:BT-smoothmoduli<1/2} can be compared with the following conditional result due to Motohashi  \cite{Mo74}.
He proved that the Lindel\"{o}f Hypothesis for Dirichlet $L$-functions allows one to take in \eqref{eq:BT}
\begin{align*}
C(\varpi)=
\begin{cases}
2, \quad &\varpi\in[0,~1/3],\\
2/(2-3\varpi), &\varpi\in[1/3,~1/2].
\end{cases}
\end{align*}

\subsection{Ingredients of proof and comparisons with previous works}

The statements of Theorems \ref{thm:BT-generalmoduli<1/2} and \ref{thm:BT-generalmoduli>1/2} are a bit complicated at first glance. The following table presents our choices of $C(\varpi)$, for some specific $\varpi>1/2,$ and shows how much we can improve upon Iwaniec \cite{Iw82}. It is worthwhile to record that any small improvements on the constants rely heavily on various estimates for exponential and character sums, individually or on average. More precisely, we appeal to spectral theory of automorphic forms for sums of Kloosterman sums, and to $\ell$-adic cohomology for short character sums and algebraic exponential sums.
\begin{table}[htbp]
\renewcommand{\arraystretch}{1.1}
\renewcommand{\tabcolsep}{1.99mm}
\begin{tabular}{ | c | c | c | c |}\hline
$\varpi$ & Our $C(\varpi)$ & Iwaniec's $C(\varpi)$  & Improvement \\ \hline
$16/31$ & $3.3067$ & $3.3514$  & $1.4\%$ \\ \hline
$12/23$ & $3.3455$ & $3.4074$  & $1.8\%$ \\ \hline
$32/61$ & $3.3889$ & $3.4366$  & $1.3\%$ \\ \hline
$8/15$ & $3.4615$ & $3.5294$  & $1.9\%$ \\ \hline
$7/13$ & $3.5254$ & $3.5862$  & $1.6\%$ \\ \hline
$6/11$ & $3.5918$ & $3.6667$  & $2.0\%$ \\ \hline
\end{tabular}
\caption{Improved $C(\varpi)$} \label{table:improvingIwaniec}
\end{table}

The following two figures present our new choices of $C(\varpi)$
for $\varpi<1/2$ and $\varpi>1/2,$ respectively, and one may see when we can succeed or not (we choose the very special exponent pair $(\kappa,\lambda)=(1/6,2/3)$ for smooth moduli and $\varpi>1/2$).

\begin{figure}[ht]
  \centering
  \begin{minipage}[b]{0.46\textwidth}
    \includegraphics[width=\textwidth]{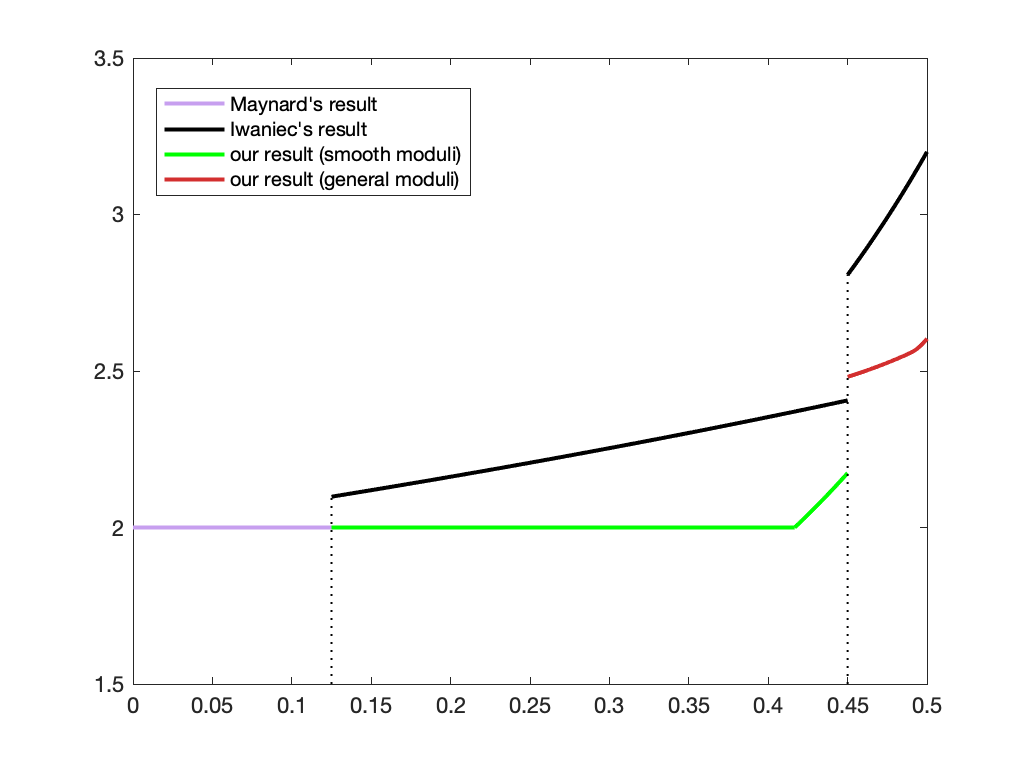}
    \caption{$C(\varpi)$ for $\varpi<1/2$}
    \label{fig:image1}
  \end{minipage}
  \hfill 
  \begin{minipage}[b]{0.46\textwidth}
    \includegraphics[width=\textwidth]{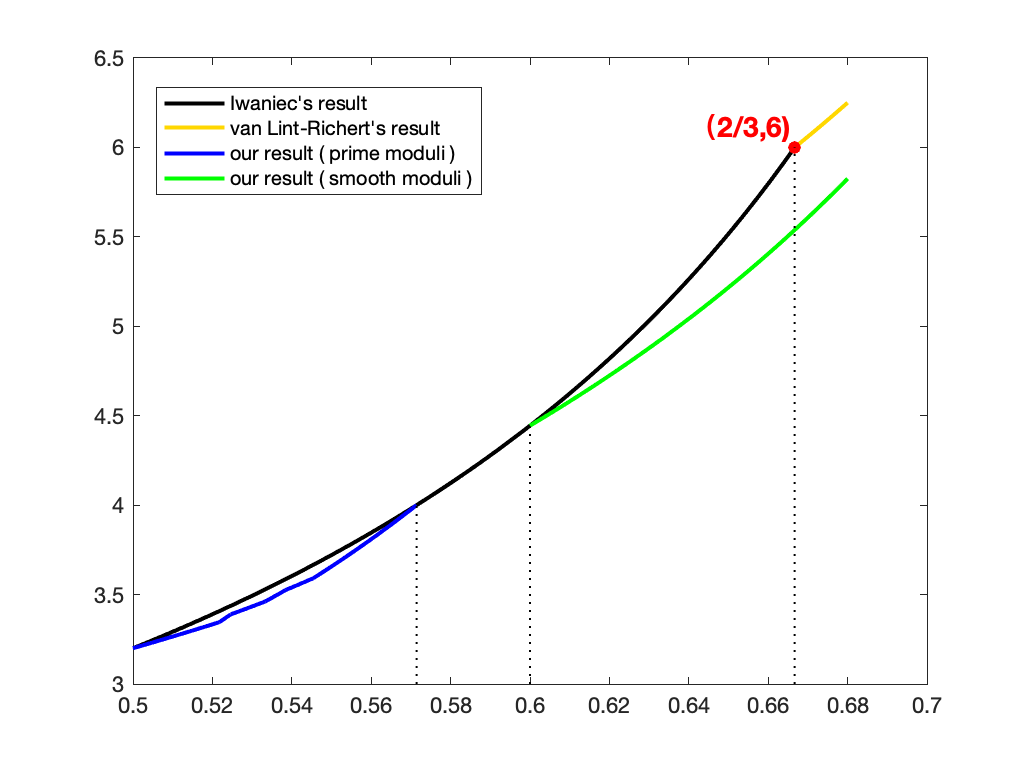}
    \caption{$C(\varpi)$ for $\varpi>1/2$}
    \label{fig:image2}
  \end{minipage}
\end{figure}

Before discussing the ingredients of our proof, it would be better to give a short survey of the methods employed in previous works.

\begin{itemize}
\item Both van Lint and Richert \cite{LR65} and Motohashi \cite{Mo74,Mo73} utilized Selberg's sieve, which is more convenient than Brun's sieve at that time. Recall that one takes $C(\varpi)=2/(1-\varpi)$ in \cite{LR65}, and the constant $2$ resulted from applications of the linear sieve with its best possible main term. 
The above denominator $(1-\varpi)$ comes from the level $(x/q)^{1-\varepsilon}$ after a {\it trivial} treatment of error terms. It was Motohashi \cite{Mo74,Mo73} who first dealt with the error term in a non-trivial manner, and he successfully employed the bilinear structure in Selberg's sieve. He employed multiplicative characters, and a basic tool in his work was the mean value of Dirichlet $L$-functions. Such an approach did not allow him to improve upon the result of van Lint and Richert \cite{LR65} for $\varpi>1/2$. 

\item We also mentioned the works of Montgomery and Vaughan \cite{MV73} and Selberg \cite{Se91} in Table \ref{table:BTconstants}. Selberg proved that
\begin{align}\label{eq:Selberg-MV}
\max_{(a,q)=1}\pi(x;q,a)<\frac{2x}{\varphi(q)(\log (q/x)+c)}
\end{align}
holds with $c=2.8$, for all $q<x/c_0$ with some $c_0>1.$ Selberg employed his own $\Lambda^2$ sieve, and published this inequality in \cite[Page 233]{Se91}. In an another direction, Montgomery and Vaughan \cite{MV73} developed the large sieve, from which they deduced that \eqref{eq:Selberg-MV} is valid for all $1\leqslant q<x$ and $c=0$.

\item Iwaniec \cite{Iw82} developed a new linear sieve of combinatorial type with well-factorable weights. This far-reaching work was originally motivated by Chen \cite{Ch75} and Motohashi \cite{Mo74}, and the merit is that one may optimize the balances between convolutions in the error term, and the sign changes among them can be captured more effectively. In particular, this well-factorization in the error term allows Iwaniec \cite{Iw82} to employ (Weil's bound for) Kloosterman sums to beat the barrier $\varpi\approx1/2$. He was able to beat van Lint and Richert \cite{LR65} even if $\varpi$ is close to $2/3-.$ On the other hand, it is worthwhile to mention that Motohashi \cite{Mo99} showed another bilinear structure in the Selberg $\Lambda^2$ sieve, which has a similar flavor to that in \cite{Iw82}. In particular, he proved that $C(\varpi)=16/(8-3\varpi)$ is admissible for all $\varpi<9/20,$ which recovered the work of Iwaniec.

\item For $\varpi\geqslant2/3,$ the above approach of Iwaniec \cite{Iw82} is no longer valid since the incomplete Kloosterman sums are extremely short, and
Friedlander and Iwaniec \cite{FI97} found another surprising approach to beat the record in \cite{LR65}. They utilized an estimate of Karatsuba \cite{Ka95} on the estimates for the bilinear form
\begin{align}\label{eq:bilinearform-Karatsuba}
\sum_m\sum_n\alpha_m\beta_n\ue\Big(\frac{h\overline{mn}}{q}\Big),
\end{align}
where $m$ and $n$, as prime variables, are restricted to very special but short segments.
In particular, they succeeded in improving van Lint and Richert \cite{LR65} for all $\varpi\in[2/3,1[.$
To apply Karatsuba's estimate for \eqref{eq:bilinearform-Karatsuba}, Friedlander and Iwaniec \cite{FI97} introduced Buchstab's iteration so that one may create prime variables in the sifting functions. By employing ideas from additive combinatorics (sum-product phenomenon), Bourgain and Garaev \cite{BG14} can bound \eqref{eq:bilinearform-Karatsuba} without restricting $m,n$ to primes. Therefore, Buchstab's iteration is not necessary, and they can give a further improvement over \cite{FI97}.

\item All above works are built on sieve methods (Brun sieve, Selberg sieve, or Rosser--Iwaniec sieve). The approach  of Maynard \cite{Ma13} is purely analytic (zero density estimates, Landau--Siegel zeros, etc), and he  proved that $\pi(x;q,a)<2\Li(x)/\varphi(q)$ for all $c<q\leqslant x^{1/8},$ where $\Li(x)=\int_2^x\ud t/\log t$ and $c$ is a computable constant. Moreover, he is also able to produce lower bounds for $\pi(x;q,a)$ with effectively computable constants.

\end{itemize}

Our work improves the above results in several aspects for different values of $\varpi.$ In principle, we need to bound the difference of convolutions
\begin{align*}
\mathop{\sum_{m\sim M}\sum_{n\sim N}\sum_{l\sim L}}_{mnl\equiv a\bmod q}\alpha_m\beta_n-\frac{1}{\varphi(q)}\mathop{\sum_{m\sim M}\sum_{n\sim N}\sum_{l\sim L}}_{(mnl,q)=1}\alpha_m\beta_n
\end{align*}
with $MN$ as large as possible, where $MNL=x.$ The coefficient $\balpha*\bbeta$ comes from the linear sieve of Iwaniec \cite{Iw82}, and each of $\balpha,\bbeta$ is bounded by divisor functions. There are at least two approaches to detect the congruence condition $mnl\equiv a\bmod q,$ via multiplicative (Dirichlet) characters or additive characters.

{\it Multiplicative Approach}: After a direct application of orthogonality of multiplicative characters, the problem might reduce to proving
\begin{align*}
\sum_{\substack{\chi\bmod q\\ \chi\neq \chi_0}}\Big|\sum_{m\sim M}\sum_{n\sim N}\sum_{l\sim L}\alpha_m\beta_n\chi(mnl)\Big|\ll x(\log x)^{-2024}.
\end{align*}
One may easily see this approach does not work unless $\varpi<1/2,$ by expecting that there exist (at most) square-root cancellations on the left hand side.
After suitable applications of H\"older's inequality, one may arrive at bounding the twisted moment
\begin{align}\label{eq:twistedmoment-ingredients}
\sum_{\substack{\chi\bmod q\\ \chi\neq \chi_0}}\Big|\sum_{r\sim R}c_r\chi(r)\Big|^2\Big|\sum_{l\sim L}\chi(l)\Big|^{2k}.
\end{align}
The underlying idea is to prove square-root cancellations among the averages over $r$ and $l$.
This expectation cannot be realized when $R$ or $k$ is suitably large, hence it becomes necessary to control large values of Dirichlet polynomials.

In his original proof, Iwaniec \cite{Iw82} studied the fourth moment of character sums twisted by Dirichlet polynomials (i.e., $k=2$ in \eqref{eq:twistedmoment-ingredients}).
In contrast with Iwaniec’s original argument, one of the novelties in proving Theorem 1.1 is that we use the Buchstab identity in the implementation of linear sieves. As a consequence, we are led to study two types of second moments of character sums twisted by Dirichlet polynomials (i.e., $k=1$ in \eqref{eq:twistedmoment-ingredients}). For the first type, the coefficient $(c_r)$ is arbitrary, and we appeal to the work of Bettin and Chandee~\cite{BC18} on trilinear forms involving Kloosterman fractions. For the second type, the coefficient $(c_r)$ has a convolution structure, arising from the Buchstab identity and 
Iwaniec's new construction of the Rosser--Iwaniec sieve. We exploit this special structure of $(c_r)$ as efficiently as possible: it allows us to apply estimates for sums of Kloosterman sums coming from spectral theory of automorphic forms, and in this way we derive a kind of Lindel\"{o}f Hypothesis on average for $R\approx q^{(5-2\vartheta)/8}$ (for further discussions, see Section \ref{sec:comments}).

By employing H\"older's inequality in a different way, we are instead led to bound fourth moment of Dirichlet $L$-functions twisted by Dirichlet polynomials (intimately related to $k=2$ in \eqref{eq:twistedmoment-ingredients}) rather than twisted second moments of character sums described above. In this setting we can apply the recent work of Blomer, Humphries, Khan and Milinovich \cite{BHKM20}, where $q$ is restricted to primes and $R$ can be as large as $q^{1/4}$, and thereby obtain another kind of Lindel\"of Hypothesis on average. This serves as a key ingredient in establishing the second assertion of Theorem~\ref{thm:BT-generalmoduli<1/2}.

We remark that the advantage of applying estimates for sums of Kloosterman sums (to twisted second moment of character sums) is that $q$ is not necessarily a prime, but we are restricted by the possible presence of exceptional spectrum of the non-Euclidean Laplacian. It is also worthwhile to remark that 
if we may factorize the coefficient $(c_r)$ as we wish, this argument can produce the same level without appealing to the work of Blomer, Humphries, Khan and Milinovich \cite{BHKM20}, modulo the Ramanujan--Petersson conjecture.

An alternative approach to study \eqref{eq:twistedmoment-ingredients} is to estimate the character sum individually. In particular, when $q$ is squarefree and has only very small prime factors, the method of arithmetic exponent pairs in \cite{WX21} can guarantee that the sum over $l$ in \eqref{eq:twistedmoment-ingredients} is $o(L)$ for very small $L$ and each $\chi\neq\chi_0$. This, amongst many other things, allows us to show that $C(\varpi)=2$ is admissible in \eqref{eq:BT} for $\varpi\in[1/8,~5/12[$. One will see this constant is best possible 
if assuming the existence of Landau--Siegel zeros.

{\it Additive Approach}: A direct application of Poisson summation will reduce the problem to proving
\begin{align*}
\frac{1}{H}\sum_{h\leqslant H}\mathop{\sum_{m\sim M}\sum_{n\sim N}}_{(mn,q)=1}\alpha_m\beta_n\ue\Big(\frac{h\overline{mn}}{q}\Big)\ll MN(\log x)^{-2024},
\end{align*}
where $H\approx qMN/x.$ It is natural to apply Cauchy--Schwarz to smooth out one variable, say $m,$ and we then arrive at the incomplete Kloosterman sum
\begin{align*} 
\sum_{\substack{m\sim M\\ (m,q)=1}}\ue\Big(\frac{a(h_1n_2-h_2n_1)\overline{n_1n_2m}}{q}\Big).
\end{align*}
The P\'olya--Vinogradov method transforms the incomplete sum to a complete one, and then Weil's bound for individual Kloosterman sums applies. This was exactly done by Iwaniec \cite{Iw82}. Our approach differs from Iwaniec's in two directions:

(1)  If $q$ is squarefree and has only very small prime factors, the factorization of $q$ can be controlled very well, so that the work on arithmetic exponent pairs, developed by Wu, Xi and Sawin \cite{WX21}, can be applied. In such situation, one can save much more than before given the incomplete Kloosterman sum of the same length, and one can also estimate the incomplete sum non-trivially even though $M$ is smaller than $q^{1/2}$.
In particular, this approach works very well when $\varpi$ is around $2/3$. The details can be found in Section \ref{subsec:BT-smoothmoduli2/3}.

(2) If $q$ has no good factorizations, say $q$ is a large prime, we also employ the P\'olya--Vinogradov method, but do not go with Weil's bound directly. Alternatively, we take advantage of averages over $h_1,h_2,n_1,n_2$, as well as the dual sum resulted from the P\'olya--Vinogradov method. By grouping $h_1,h_2,n_1,n_2$ in certain ways, we may create some bilinear and trilinear forms with Kloosterman sums. This produces possibilities to capture sign changes of Kloosterman sums. For the bilinear form, we need to control a sum of product of many Kloosterman sums, which will appeal to \cite{FGKM14} or \cite{Xi17}; while for the trilinear form, we reduce the problem to bounding the bilinear form
\begin{align*}
\sum_m\sum_n K(mn),
\end{align*}
where $K:\bZ/q\bZ\rightarrow\bC$ is given by the Fourier transform of a product of two Kloosterman sums.
Hence the work of Fouvry, Kowalski and Michel \cite{FKM14} on smooth bilinear forms with general trace functions will apply, which is rooted in $\ell$-adic cohomology and spectral theory of Eisenstein series. The details can be found in the first several parts of Section \ref{sec:BT-generalmoduli>1/2}, as well as Section \ref{subsec:tracefunctionfinitefield}.

\subsection{Connections with Landau--Siegel zeros}
It should be mentioned that the reduction of $C(\varpi)$ is of great importance for the location of Landau--Siegel zeros.
Motohashi \cite{Mo79} showed that if for $x>q^c,$
\begin{align*}
\pi(x;q,a)<(2-\delta)\frac{1}{\varphi(q)}\frac{x}{\log(x/q)}
\end{align*}
holds with some constant $\delta>0$, the real zero $\beta$ (if it exists) of $L(s,\chi)$ then satisfies
\begin{align*}
\beta<1-\frac{c'\eta}{\log(q+3)},
\end{align*}
 implying there are no Landau--Siegel zeros. The phenomenon that such improvements over the Brun-Titchmarsh theorem eliminate Landau--Siegel zeros seems first observed by Klimov \cite[Remark 1]{Kl61}; see also \cite{Si83} and \cite{Gr22} for related results.

On the other hand, it is believed that $\pi(x;q,a)\leqslant 2\Li(x)/\varphi(q)$ as long as $q$ is at most a small power of $x$; see \cite{Mo12} for the details and history.
Recently, Maynard \cite{Ma13} demonstrated this expectation is valid for all $q^8\leqslant x.$
His method involves a log-free zero density estimate for Dirichlet $L$-functions and careful analysis of Landau--Siegel zeros.
Up to a very small constant $\varepsilon$, we have shown in Theorem \ref{thm:BT-smoothmoduli<1/2} that $\pi(x;q,a)$ is at most twice the conjectural bound, provided that $q<x^{5/12},$ $q$ is squarefree and all prime factors are sufficiently small.

\subsection*{Notation and convention}  
We write $\ue(t)=\ue^{2\pi\mathrm{i}t}$. We use $\varepsilon$ to denote a small positive number, which might be different at each occurrence. The notation $n\sim N$ means $N<n\leqslant2N.$  We write $A\leqslant B^{-}$ if $A\leqslant B^{1-\delta}$ holds for some constant $\delta>0.$ Denote by $\varphi,$ $\Lambda$ and $\mu$ the Euler, von Mangoldt and M\"obius functions, respectively. For a complex coefficient $\balpha=(\alpha_n),$ we define the $\ell_2$-norm $\|\balpha\|=(\sum_n|\alpha_n|^2)^{\frac{1}{2}}.$


\smallskip

\section{Dirichlet polynomials and character sums}
\label{sec:charactersum}

The following lemma follows from the orthogonality of Dirichlet characters; see \cite[Theorem 6.2]{Mo71}.
\begin{lemma}\label{lm:Dirichletpolynomial-secondmoment}
For any complex sequence $\balpha=(\alpha_n)_{n\sim N},$ we have
\begin{align*}
\sum_{\chi\bmod q}\Big|\sum_{n\sim N} \alpha_n\chi(n)\Big|^2\ll (N+q)\|\balpha\|^2.
\end{align*}
\end{lemma}

We also need Huxley's large value theorem for Dirichlet polynomials; see \cite{Hu74} for the original proof, as well as Jutila \cite{Ju77} for another elegant proof.
\begin{lemma}\label{lm:largevalue-Huxley}
For any complex sequence $\balpha=(\alpha_n)_{n\sim N}$ and $V>0,$ we have
\begin{align*}
\#\Big\{\chi\bmod q:\Big|\sum_{n\sim N} \alpha_n\chi(n)\Big|>V\Big\} \ll \|\balpha\|^2 N V^{-2}+ \|\balpha\|^6q^{1+\varepsilon} N V^{-6}.
\end{align*}
\end{lemma}

In the proof of Theorem \ref{thm:BT-generalmoduli<1/2}, we may reduce the problem to bounding the second moment of Dirichlet $L$-functions weighted by a Dirichlet polynomial:
\begin{align*}
\frac{1}{\varphi(q)}\sum_{\substack{\chi\bmod{q}\\ \chi\neq \chi_{0}}}\Big|\sum_{m} c_m \chi(m)\Big|^2|L(\tfrac{1}{2}+it, \chi)|^2.
\end{align*}
Alternatively, we may also consider the second moment of incomplete character sums weighted by a Dirichlet polynomial:
\begin{align*}
\frac{1}{\varphi(q)}\sum_{\substack{\chi\bmod{q}\\ \chi\neq \chi_{0}}}\Big|\sum_{m} c_m \chi(m)\Big|^2\Big|\sum_{n\leqslant N}\chi(n)n^{it}\Big|^2.
\end{align*}
We hope to prove the Lindel\"of Hypothesis on average with the support of $c_m$ as large as possible. In what follows, we turn to a slightly general situation.

Let $F\in\cC_0^\infty(\bR^2)$ be a Schwartz function with support in $]L,2L]\times]L,2L]$, satisfying
\begin{align*}
\frac{\partial^{\nu_1+\nu_2}}{\partial x_1^{\nu_1}\partial x_2^{\nu_2}}F(x_1,x_2)
&\ll \eta^{\nu_1+\nu_2}L^{-\nu_1-\nu_2}
\end{align*}
for all $\nu_1,\nu_2\geqslant0$ and some $\eta\geqslant1.$ Let $q\ge 1$, and let $\balpha=(\alpha_n)$ and $\bbeta=(\beta_n)$ be complex sequences supported in $(N,2N]$. We put
\begin{align}\label{eq:cR-countingcongruence}
\cR=\cR(\balpha,\bbeta,F;q):=\mathop{\sum\sum\sum\sum}_{\substack{n_1,n_2\sim N,~l_1,l_2\in\bZ\\n_1l_1\equiv n_2l_2\bmod q\\  (n_1n_2l_1l_2,q)=1}}\alpha_{n_1}\beta_{n_2}F(l_1,l_2).
\end{align}

We have an asymptotic formula for the above counting function.

\begin{theorem}\label{thm:countingcongruence-generalcoefficient}
Let $\eta,q,N,L\geqslant1.$ Suppose $\balpha,\bbeta$ are  complex sequences supported in $[N,2N],$ and $F$ is a Schwartz function defined as above.
Then
\begin{align*}
\cR
&=\frac{\varphi(q)}{q^2}\Big(\mathop{\sum\sum}_{\substack{n_1,n_2\sim N\\(n_1n_2,q)=1}}\alpha_{n_1}\beta_{n_2}\Big)\Big(\iint_{\bR^2} F(y,z)\ud y\ud z\Big)\\
&\ \ \ \ \ +O(\eta (NL)^{1+\varepsilon}+\eta^{\frac{7}{2}}N^{\frac{23}{8}}Lq^{-1+\varepsilon}(N^{-\frac{9}{40}}q^{\frac{3}{20}}+1))
\end{align*}
for any $\varepsilon>0,$
with an implied constant depending only on $\varepsilon.$
\end{theorem}

Theorems \ref{thm:countingcongruence-generalcoefficient} and \ref{thm:countingcongruence-specialcoefficient} will be utilized to study a weighted moment of character sums, which can be regarded as a kind of Lindel\"of Hypothesis on average.
The proof of Theorem \ref{thm:countingcongruence-generalcoefficient} relies on the following estimate, due to Bettin and Chandee \cite{BC18}, for trilinear forms involving Kloosterman fractions.

\begin{lemma}\label{lm:trilinearform-BC}
Let $\balpha=(\alpha_m),$ $\bbeta=(\beta_n)$ and $\bdelta=(\delta_k)$ be complex sequences with supports in $[1,M],[1,N]$ and $[1,K],$ respectively. Let $a$ be a non-zero integer. Then
\begin{align*}
\sum_{k\leqslant K}\delta_k\sum_{m\leqslant M}\sum_{n\leqslant N}\alpha_m\beta_n\ue\Big(ak\frac{\overline{m}}{n}\Big)
& \ll\|\balpha\|\|\bbeta\|\|\bdelta\|(KMN)^{\frac{7}{20}+\varepsilon}\Big(1+\frac{|a|K}{M N}\Big)^{\frac{1}{2}} \\
& \qquad\times((M+N)^{\frac{1}{4}}+(KMN)^{\frac{1}{40}}(KM+KN)^{\frac{1}{8}}).
\end{align*}
\end{lemma}

\begin{remark}
The above estimate of Bettin and Chandee \cite{BC18} is in fact a successor of a (weaker) bound due to Duke, Friedlander and Iwaniec \cite{DFI97b} with $K=1$.
\end{remark}

If the coefficients $\balpha,\bbeta$ have suitable convolution structures, we can do a bit better than Theorem \ref{thm:countingcongruence-generalcoefficient}.

\begin{theorem}\label{thm:countingcongruence-specialcoefficient}
Let $\eta,q,N,L\geqslant1.$  Suppose $\balpha,\bbeta$ admit the convolution structure that $\balpha=\balpha'*\balpha'',$ $\bbeta=\bbeta'*\bbeta''$, where for some $M_1,N_1,M_2,N_2\geqslant1$ with $M_1N_1=M_2N_2= 2N$,
the complex sequences $\balpha',$ $\balpha'',$ $\bbeta',$ $\bbeta''$ are supported on $[1,M_1],$ $[1,N_1],$ $[1,M_2]$ and $[1,N_2],$ respectively. Moreover, we assume $\balpha''$ and $\bbeta''$ vanish on the set $\{n\geqslant2:\omega(n)>(\log N)^{1/3}\text{ or }\mu(n)=0\}.$
We further assume that $|\alpha_n|,|\beta_n|\leqslant (\tau(n)\log n)^c$ for some constant $c>0$.
Then
\begin{align}\label{eq:R-asymptotic}
\cR
&=\frac{\varphi(q)}{q^2}\Big(\mathop{\sum\sum}_{\substack{n_1,n_2\sim N\\(n_1n_2,q)=1}}\alpha_{n_1}\beta_{n_2}\Big)\Big(\iint_{\bR^2} F(y,z)\ud y\ud z\Big)+\cE,
\end{align}
where 
\begin{align*}
\cE&\ll \eta^{\frac{9}{2}}N^{\frac{3}{2}}Lq^{\frac{\vartheta-1}{2}}(N_1N_2)^{\frac{1}{2}}(N^2q^{-1}+N_1N_2)^{\frac{1}{4}}\Big(1+\frac{N^2}{qN_1N_2}\Big)^{\frac{1}{4}}(qNL)^{\varepsilon}\\
&\ \ \ \ +\eta (NL)^{1+\varepsilon}+\eta^{\frac{9}{2}}N^2L(qN_1N_2)^{-\frac{1}{2}}(qNL)^\varepsilon
\end{align*}
for any $\varepsilon>0,$
with an implied constant depending only on $(\varepsilon,c)$. Here $\vartheta\leqslant 7/64$ is the exponent towards the Ramanujan--Petersson conjecture for $\GL_2(\bQ).$
\end{theorem}

\begin{remark}
At first glance, the restriction that $\balpha''$ and $\bbeta''$ vanish on the set $\{n\geqslant2:\omega(n)>(\log N)^{1/3}\}$ is quite artificial. In our applications (as well as many other situations) this is harmless since only a minority of integers have been ruled out. We remark that this condition is designed for the later application of a trick due to Fouvry (see Lemma \ref{lm:Fouvrypartition} below).
\end{remark}

The control on the error term in Theorem \ref{thm:countingcongruence-specialcoefficient} depends heavily on well-factorizations of $\balpha,\bbeta.$ The proof will be given in Section \ref{sec:countingcongruence-proof}, and we appeal to Kloostermania in the following version.

\begin{lemma}\label{lm:trilinearKloosterman-inverse}
Let $C,M,N\geqslant1,$ $a,r,s\in\bZ^+$ pairwise coprime with $\mu(rs)\neq0,$ and $g(c,m,n)$ a smooth function
with compact support in $[C,2C]\times[M,2M]\times[N,2N]$ such that
\begin{align*}
\frac{\partial^{\nu_1+\nu_2+\nu_3}}{\partial c^{\nu_1}\partial m^{\nu_2}\partial n^{\nu_3}}g(c,m,n)\ll c^{-\nu_1}m^{-\nu_2}n^{-\nu_3}
\end{align*}
for all $\boldsymbol\nu=(\nu_1,\nu_2,\nu_3)\in\bN^3,$ the implied constant in $\ll$ depending at most on $\boldsymbol\nu.$ 

For any complex sequence $\boldsymbol\beta=(\beta_n),$ we have
\begin{align*}
\mathop{\sum_c\sum_m\sum_n}_{(sc,rm)=1}\beta_n g(c,m,n)\ue\Big(\frac{\pm an\overline{rm}}{sc}\Big)&\ll (arsCMN)^\varepsilon \|\boldsymbol\beta\|\cdot\Delta,
\end{align*}
where
\begin{align*}
\Delta&=a^\vartheta (rsCM)^{\frac{1}{2}}\frac{(1+1/X)^{\vartheta}}{1+X^{\frac{1}{2}}}\Big(1+X+\frac{C}{rM}\Big)^{\frac{1}{2}}\Big(1+X+\frac{N}{rs}\Big)^{\frac{1}{2}}+s^{-1}MN^{\frac{1}{2}}
\end{align*}
with $X=aN/(srCM).$
Here $\vartheta\leqslant7/64$ is the exponent towards the Ramanujan--Petersson conjecture  for $\GL_2(\bQ).$
\end{lemma}

Lemma \ref{lm:trilinearKloosterman-inverse} is essentially due to Deshouillers and Iwaniec \cite{DI82}, and the above version is directly taken from Matom\"aki \cite[Lemma 10]{Ma09}.

Theorems \ref{thm:countingcongruence-generalcoefficient} and \ref{thm:countingcongruence-specialcoefficient} allow us to derive the following 
versions of Lindel\"of Hypothesis (or square-root cancellation) on average.

\begin{lemma}\label{lm:weightedL-secondmoment-generalcoefficient}
Let $N,L,q\geqslant3$ and let $g$ be a smooth function with compact support in $[1,2],$ satisfying $\|g^{(j)}\|_\infty\ll \eta^j$ for some $\eta\geqslant1.$
Let $\blambda=(\lambda_n)$ be a  complex sequence with support in $[1,N].$ Suppose that $N,L\ll q^{O(1)}.$ Then we have
\begin{align*}
\frac{1}{\varphi(q)}&\sum_{\substack{\chi\bmod{q}\\ \chi\neq \chi_{0}}}
\Big|\sum_n\frac{\lambda_n\chi(n)}{\sqrt{n}}\Big|^2\Big|\sum_{l\in\bZ}g\Big(\frac{l}{L}\Big)\frac{\chi(l)}{\sqrt{l}}\Big|^2
\ll \eta q^{\varepsilon}+\eta^{\frac{7}{2}}N^{\frac{15}{8}}q^{-1+\varepsilon}(N^{-\frac{9}{40}}q^{\frac{3}{20}}+1)
\end{align*}
for any $\varepsilon>0,$ where the implied constant depends only on $\varepsilon.$
\end{lemma}

\begin{lemma}\label{lm:weightedL-secondmoment-specialcoefficient}
Let $M,N,L,q\geqslant3$ and let $g$ be a smooth function with compact support in $[1,2],$ satisfying $\|g^{(j)}\|_\infty\ll \eta^j$ for some $\eta\geqslant1.$
Let $\balpha,\bbeta$ be complex sequences with supports in $[1,M],$ $[1,N],$ respectively, and $\|\balpha\|_\infty,\|\bbeta\|_\infty\leqslant1.$ We also assume that $\balpha,\bbeta$ vanish on the set $\{n\geqslant2:\omega(n)>(\log MN)^{1/3} \text{ or } \mu(n)=0\}.$ Suppose that $M,N,L\ll q^{O(1)}.$ Then we have
\begin{align*}
\frac{1}{\varphi(q)}&\sum_{\substack{\chi\bmod{q}\\ \chi\neq \chi_{0}}}
\Big|\sum_m\frac{\alpha_m\chi(m)}{\sqrt{m}}\Big|^2\Big|\sum_n\frac{\beta_n\chi(n)}{\sqrt{n}}\Big|^2\Big|\sum_{l\in\bZ}g\Big(\frac{l}{L}\Big)\frac{\chi(l)}{\sqrt{l}}\Big|^2\\
&\ll \eta^{\frac{9}{2}}(q^{\frac{\vartheta-1}{2}}M^{\frac{1}{2}}N^2+1)(Mq^{-\frac{1}{2}}+1)q^{\varepsilon}
\end{align*}
for any $\varepsilon>0,$ where the implied constant depends only on $\varepsilon.$
\end{lemma}

One may rewrite the product of two Dirichlet polynomials in Lemma \ref{lm:weightedL-secondmoment-specialcoefficient} as one Dirichlet polynomial with convoluted coefficients. This special convolution structure allows us to obtain square-root cancellations on average for the Dirichlet polynomial longer than 
Lemma \ref{lm:weightedL-secondmoment-generalcoefficient}. In this sense, Lemma \ref{lm:weightedL-secondmoment-specialcoefficient} is new. We will give more precise comparisons in Section \ref{subsec:Twistedsecond momentcharacter}.
See also \cite{BPZ20} for the second moment of Dirichlet $L$-functions twisted by Dirichlet polynomials, but with asymptotic formulae.

The following lemma presents a weighted fourth moment of Dirichlet $L$-functions. 

\begin{lemma}\label{lm:weightedL-fourthmoment}
Let $q$ be a large prime and $N\geqslant 1$. Let $\balpha=(\alpha_n)_{n\leqslant N}$ be a complex sequence supported on squarefree numbers and $\|\balpha\|_\infty\leqslant1.$ Then we have
\begin{align*}
\frac{1}{\varphi(q)}\sum_{\substack{\chi\bmod{q}\\ \chi\neq \chi_{0}}}
\Big|\sum_{n\leqslant N}\frac{\alpha_n\chi(n)}{\sqrt{n}}\Big|^2|L(\tfrac{1}{2}+it, \chi)|^4
\ll (1+|t|)^{16} (N^2q^{-\frac{1}{2}}+1)q^\varepsilon
\end{align*}
for any $\varepsilon>0$ and $t\in\bR,$ where the implied constant depends only on $\varepsilon.$
\end{lemma}
Lemma \ref{lm:weightedL-fourthmoment} is proved by Blomer, Humphries, Khan and Milinovich \cite{BHKM20} for $t=0$, and their argument can be easily extended to all $t\in\bR$. It is very useful to relax $q$ to an arbitrary modulus, which in turn can guarantee the choice \eqref{eq:C*(varpi):varpi<1/2,primemoduli} in Theorem \ref{thm:BT-generalmoduli<1/2} for all $q\geqslant3$.

\smallskip

\section{General trace functions and arithmetic exponent pairs}
\label{sec:tracefunction}

Parts of this paper will employ $\ell$-adic cohomology from algebraic geometry to study the necessary exponential and character sums. Many of them would be formulated in the language of trace functions which will be recalled below.

\subsection{Trace functions over finite fields}\label{subsec:tracefunctionfinitefield}
For an odd prime $p$, we consider  middle-extension sheaves on $\bA^1_{\bF_p}$. Recall that a sheaf $\cF$ on $\bA^1_{\bF_p}$ is said to be a middle-extension sheaf if there exists a non-empty open subset $U\subset \bA^1_{\bF_p}$ and an open immersion $j:U\hookrightarrow \bA^1_{\bF_p}$ such that $\cF$ is lisse on $U$ and $\cF\simeq j_*(j^*\cF)$. Fix an isomorphism $\iota : \overline{\bQ}_\ell\rightarrow\bC$ for a $\ell\neq p$. 
Following Katz \cite[Section 7.3.7]{Ka88}, we define trace functions as follows.

\begin{definition}[Trace functions]\label{def:tracefunction}
Let $\cF$ be an $\ell$-adic middle-extension sheaf pure of weight zero, 
which is lisse on an open set $U$. The trace function associated to $\cF$ is defined by
\begin{align*}
K:x\in\bF_p\mapsto\iota(\tr(\Frob_x\mid V_\cF)),
\end{align*}
where $\Frob_x$ denotes the geometric Frobenius at $x\in\bF_p,$ and $V_\cF$ is a finite dimensional $\overline{\bQ}_\ell$-vector space, which corresponds to a continuous finite dimensional Galois representation and unramified at every closed point $x$ of $U.$
\end{definition}

Motivated by many applications in analytic number theory, the conductor of a (middle-extension) sheaf was introduced by Fouvry, Kowalski and Michel \cite{FKM14} to measure the geometric complexity of the trace function.
For an $\ell$-adic middle-extension sheaf $\cF$ on $\bP^1_{\bF_p}$ of rank $\rank(\cF)$, the conductor of $\cF$ is defined to be
\begin{align*}  
\fc(\cF)= \rank(\cF) + \sum_{x\in S(\cF)} (1+\Swan_x(\cF)),
\end{align*}
where $S(\cF)\subset\bP^1(\overline{\bF}_p)$ denotes the $($finite$)$ set of singularities of $\cF$, 
and $\Swan_x(\cF)$ $(\geqslant 0)$ denotes the Swan conductor of $\cF$ at $x$.

We have a lot of practical examples of trace functions arising in modern analytic number theory.
Among many interesting examples, we would like to mention additive characters, multiplicative characters, and Kloosterman sums. 
\begin{itemize}
\item Let $f\in\bF_p(X)$ be a rational function, and $\psi$ a primitive additive character on $\bF_p$, then $\psi(f(x))$ is a trace function of an $\ell$-adic middle-extension sheaf, which is taken to be zero when $x$ is a pole of $f$. More precisely, one can show that there exists an $\ell$-adic middle-extension sheaf modulo $p$, denoted by
$\cL_{\psi(f)},$ such that $x\mapsto\psi(f(x))$ is the trace function of $\cL_{\psi(f)}.$ The conductor can be bounded in terms of the degree of $f$, independent of $p$.

\item  
Let $f\in\bF_p(X)$ be a rational function, and $\chi$ a multiplicative character of order $d >1$. 
If $f$ has no pole or zero of order divisible by $d$, 
then one can show that there exists an $\ell$-adic middle-extension
sheaf, denoted by $\cL_{\chi(f)}$, such that $x\mapsto\chi(f(x))$ is the trace function of $\cL_{\chi(f)}$. 
The conductor can also be bounded in terms of the degree of $f$, independent of $p$.

\item 
Another example is the following {\it normalized} Kloosterman sum 
defined by
\begin{align*}
\kl(\cdot,p) : x\in\bF_p\mapsto \frac{1}{\sqrt{p}}
\sum_{y\in\bF_p^\times}\ue\Big(\frac{y+x/y}{p}\Big).
\end{align*}
According to Deligne, there exists an $\ell$-adic middle-extension sheaf $\cK\ell$, called a Kloosterman sheaf, such that $\kl(\cdot,p)$ gives its trace function over 
$\bF_p^\times.$ The sheaf $\cK\ell$ is geometrically irreducible of rank $2$, with conductor bounded by $5$.
\end{itemize}

A middle-extension sheaf on $\bP_{\bF_p}^1$ that is pointwise pure of weight 0 is said to be a {\it Fourier sheaf} if no geometrically irreducible component is geometrically isomorphic to an Artin--Schreier
sheaf  $\cL_\psi$  attached to some additive character $\psi$ of $\bF_p.$

The following lemma collects some properties of Fourier transforms of Fourier sheaves due to Deligne \cite{De80}, Laumon \cite{La87}, 
Brylinski \cite{Br86}, Katz \cite{Ka88, Ka90} and Fouvry, Kowalski and Michel \cite{FKM15}.

\begin{lemma}\label{lm:Fouriertransform}
Let $\psi$ be a non-trivial additive character of $\bF_p$ and $\cF$ a Fourier sheaf on $\bP_{\bF_p}^1$ with the trace function $K_\cF.$  Then there exists an $\ell$-adic sheaf $\cG=\ft_\psi(\cF)$
called the Fourier transform of $\cF$, which is also an $\ell$-adic Fourier sheaf with the trace function
\begin{align*}
K_{\cG}(y)
= \ft_\psi(K_\cF)(y)
:= -\frac{1}{\sqrt{p}}\sum_{x\in\bF_p}K_\cF(x)\psi(yx).
\end{align*}
Furthermore, we have
\begin{itemize}

\item 
The sheaf $\cG$ is geometrically irreducible, or geometrically isotypic, if and only if $\cF$ is;

\item 
The Fourier transform is involutive, in the sense that we have a canonical arithmetic isomorphism
$\ft_\psi(\cG)\simeq[\times(-1)]^*\cF,$
where $[\times(-1)]^*$ denotes the pull-back by the map $x\mapsto-x;$

\item 
We have $\fc(\ft_\psi(\cF))\leqslant10\fc(\cF)^2.$
\end{itemize}

\end{lemma}

For $a,b,y\in\bF_p$, put
\begin{align*}
V_p(y;a,b)=\frac{1}{\sqrt{p}}\sum_{x\in\bF_p}\kl(ax,p)\kl(bx,p)\ue\Big(\frac{-yx}{p}\Big).
\end{align*}
In the proof of Theorem \ref{thm:BT-generalmoduli>1/2}, we need to show that $y\mapsto V_p(y;a,b)$ is a good trace function in some sense, provided that $a\neq b$.

\begin{lemma}\label{lm:Kloostermantensor} 
Let $\cK\ell$ be the Kloosterman sheaf of rank $2$ as above. For all $a,b\in\bF_p^\times$ with $a\neq b,$ the tensor product
$\cF:=([\times a]^*\cK\ell)\otimes ([\times b]^*\cK\ell)$ is a Fourier sheaf, which is geometrically irreducible, pure of weight zero, of rank $4,$ and has the geometric monodromy group $\SL_2\times \SL_2.$ The conductor of $\cF$ is at most $25.$

Moreover, for each non-trivial additive character $\psi$ of $\bF_p,$ the Fourier transform $\ft_\psi(\cF)$ is also a Fourier sheaf, which is geometrically irreducible and of conductor at most $6250.$
\end{lemma} 

\proof
The first part has been proven by Fouvry and Michel \cite[Section 7]{FM07}. The second part follows from Lemma \ref{lm:Fouriertransform} immediately.
\endproof

Before going into more details on trace functions over finite fields, we would like to mention another definition of Kloosterman sums defined over $\bZ/q\bZ:$ for $q\geqslant1$, put
\begin{align*}
S(m,n;q)=\sideset{}{^*}\sum_{a\bmod q}\ue\Big(\frac{ma+n\overline{a}}{q}\Big),
\end{align*}
where $m,n\in\bZ$ and $a\overline{a}\equiv1\bmod q.$ It is clear that $\kl(x,p)\sqrt{p}=S(x,1;p)=S(1,x;p)$ for each prime $p$, if we identify $\bF_p$ with $\{0,1,\cdots,p-1\}.$
We call $S(*,0;q)$ a Ramanujan sum, which satisfies $|S(m,0;q)|\leqslant (m,q);$ in general, we have the Weil bound 
\begin{align*}
|S(m,n;q)|\leqslant q^{\frac{1}{2}}(m,n,q)^{\frac{1}{2}}\tau(q).
\end{align*}

We now recall a very strong estimate of Fouvry, Kowalski and Michel \cite{FKM14} on smooth bilinear forms with general trace functions.

\begin{lemma}\label{lm:FKM} 
Let $K$ be a trace function of an $\ell$-adic middle-extension sheaf $\cF$ on $\bA_{\bF_p}^1,$ which is geometrically isotypic, pointwise pure of weight $0,$ and no geometrically irreducible component of which is geometrically isomorphic to an Artin--Schreier
sheaf  $\cL_\psi$ attached to some additive character $\psi$ of $\bF_p.$

Let $M, N \geqslant 1$, and let $V_1,V_2$ be two fixed smooth functions with compact supports in $[1,2].$  We then have
\begin{align*}
\mathop{\sum\sum}_{m, n \in\bZ} K(m n) V_1\Big(\frac{m}{M}\Big) V_2\Big(\frac{n}{N}\Big)
\ll M N\Big(1+\frac{p}{M N}\Big)^{\frac{1}{2}} p^{-\frac{1}{8}+\varepsilon}
\end{align*}
for any $\varepsilon>0,$ where the implicit constant depends only on $\varepsilon$ and the conductor of $\cF.$
\end{lemma}

\begin{lemma}\label{lm:Klsum-2nu-thmoment} 
Let $p$ be a large prime and $ \nu\in \bZ^+.$ For each subset $\cN\subseteq[1,p]$ and arbitrary coefficient $\bbeta=(\beta_n)$ with $\|\bbeta\|_\infty\leqslant1,$ we have
\begin{align*}
\sum_{m\bmod p}\Big|\sum_{ n\in \cN}\beta_n \kl(mn,p)\Big|^{2\nu}
\ll |\cN|^{\nu}p+|\cN|^{2\nu}p^{\frac{1}{2}},
\end{align*}
where the implied constant depends only on $\nu.$
\end{lemma}

\proof Denote by $S$ the quantity in question. 
We may write
\begin{align*}
|S|\leqslant \sum_{\bn\in \cN^{2\nu}}|S(\bn)|,
\end{align*}
where for $\bn=(n_1, n_2, \ldots, n_{2 \nu}) \in \cN^{2\nu}$,
\begin{align*}
S(\bn)=\sum_{x \bmod p} \prod_{1 \leqslant j \leqslant 2\nu} \kl(n_j x,p).
\end{align*}
There is a trivial bound that
\begin{align*}
S(\bn)\ll p.
\end{align*}
We also expect there exist square-root cancellations if the coordinates of $\bn$ appear in suitable configurations. More precisely, according to \cite{FGKM14} or \cite{Xi17},
we have
\begin{align*}
S(\bn)\ll \sqrt{p}
\end{align*}
as long as the multiplicity
\begin{align*}
\mu_n(\bn):=|\{1 \leqslant j \leqslant 2\nu:n_j=n\}|
\end{align*}
is odd for some $n\in \cN$. In this way, we employ the above trivial bound if $\mu_n(\bn)$ is even for all $n\in \cN$.
This proves the lemma readily.
\endproof

\subsection{Composite trace functions and arithmetic exponent pairs}\label{subsec:exponentpairs}
By the Chinese remainder theorem, a trace function can  also be defined for squarefree moduli, which allows Wu, Xi and Sawin \cite{WX21} to develop the general form of $q$-analogue of van der Corput method and arithmetic exponent pairs.

In terms of trace functions over finite fields, we may construct {\it composite} trace functions following the manners in \cite{Po14} and \cite{WX21}.
Let $q\geqslant1$ be a squarefree number and $K_q$ a {\it composite} trace function (formally) defined by the product
\begin{align}\label{eq:def-Kq}
K_q(n)=\prod_{p\mid q}K_p(n),\end{align}
where $K_p$ is the trace function of a certain $\ell$-adic middle-extension sheaf $\cF_p$ on $\bA^1_{\bF_p}$ for each $p\mid q$. We adopt the convention that $K_q(n)=1$ for all $n$ if $q=1.$
In practice, the value of $K_p(n)$ may depend on the complementary divisor $q/p.$ 
Originally, $K_p$ is a function living in $\bZ/p\bZ$; one can extend its domain to $\bZ$ by periodicity, and one can do the same thing to $K_q$.

For a given interval $I$, we consider the following average
\begin{align*}
\sum_{n\in I}K_q(n).
\end{align*}
Roughly speaking, the P\'olya--Vinogradov method can be used to bound such sums non-trivially for $|I|>q^{\frac{1}{2}+\varepsilon}$, and it was observed by Heath-Brown \cite{HB78} that one may go further if $q$ has good factorizations, and this led him to prove a Weyl-type subconvexity bound for $L(\tfrac{1}{2}, \chi)$ to smooth moduli. Inspired by Heath-Brown, Wu--Xi--Sawin \cite{WX21} developed the method of arithmetic exponent pairs for averages of composite trace functions that may go far beyond the P\'olya--Vinogradov bound, as long as the factorization of $q$ is good enough. 
We now summarize results of the arithmetic exponent pairs developed in \cite{WX21}. For more background and applications of the $q$-analogue of van der Corput method and arithmetic exponent pairs, the readers may refer to \cite{WX21} directly and references therein.

Following Sawin \cite{WX21}, we call a middle-extension sheaf $\cF_p$ on $\bA_{\bF_p}^1$ {\it universally amiable} if 
\begin{enumerate}[(i)]
\item it is pointwise pure of weight $0$;
\item no geometrically irreducible component of $\cF_p$ is geometrically isomorphic to an Artin--Schreier sheaf $\cL_{\psi}$ attached to some additive characer $\psi$ of $\bF_p$;
\item its local monodromy at $\infty$ has all slopes $\leq 1$. 
\end{enumerate}

In such case, we also say the associated trace function $K_p$ is universally-amiable. A composite trace function $K_q$ is said to be {\it compositely universally amiable} if for each $p\mid q,$ $K_p$ can be decomposed into a sum of universally amiable trace functions, in which case we also say the corresponding sheaf $\cF:=(\cF_p)_{p\mid q}$ is compositely universally amiable.

To mention the main results in \cite{WX21}, we assume that
\begin{itemize}
\item $\eta>0$ is a sufficiently small number;
\item $q$ is a squarefree number with no prime factors exceeding $q^\eta$;
\item $\delta$ is positive integer such that $(\delta,q)=1$;
\item $K_q$ is a compositely universally amiable trace function $\bmod q$;
\item $I$ is an interval with $N:=|I|<q\delta;$
\item there exists some uniform constant $\fc>0,$ such that $\fc(\cF_p)\leqslant \fc$ for each $p\mid q$.
\end{itemize}
We further assume that $\delta$ is also squarefree and no prime factors exceed $\delta^\eta$, as will be satisfied in our applications in this paper.
We are interested in the upper bound of the average
\begin{align*}
\fS(K,W)=\sum_{n\in I}K_q(n)W_\delta(n),
\end{align*}
where $W_\delta:\bZ/\delta\bZ\rightarrow\bC$ is an arbitrary function, which we call deformation factor later.
Let $(\kappa,\lambda,\nu)$ be a tuple such that
\begin{align}\label{eq:exponentpair-estimate}
\mathfrak{S}(K, W)
\ll_{\eta,\varepsilon, \fc} N^{\varepsilon+O(\eta)}\|W_\delta\|_{\infty}(q/N)^{\kappa}N^{\lambda}\delta^{\nu}
\end{align}
for any $\varepsilon>0,$ 
where the implied constant is allowed to depend on $\eta,\varepsilon$ and $\fc.$ 
We define two maps
\begin{align}\label{eq:exponentpair-afterA}
A\cdot(\kappa,\lambda,\nu)
= \Big(\frac{\kappa}{2(\kappa+1)},~\frac{\kappa+\lambda+1}{2(\kappa+1)},~\frac{\kappa}{2(\kappa+1)}\Big)
\end{align}
and
\begin{align}\label{eq:exponentpair-afterB}
B\cdot(\kappa,\lambda,\nu)
=\Big(\lambda-\frac{1}{2}, ~\kappa+\frac{1}{2}, ~\lambda+\nu-\kappa\Big),
\end{align}
resulted from $A$- and $B$-processes in the $q$-analogue of van der Corput method; see \cite{WX21} for details.
Note that \eqref{eq:exponentpair-estimate} holds trivially for $(\kappa,\lambda,\nu)=(0,1,0).$ Under different combinations of $A$ and $B$, we may produce a series of tuples starting from $(\kappa,\lambda,\nu)=(0,1,0).$ Namely, we may introduce
\begin{align}\label{eq:arithmeticexponentpairs-set}
\fP=\{A^{m_1}B^{n_1}A^{m_2}B^{n_2}\cdots A^{m_k}B^{n_k}\cdot (0,1,0):k\geqslant1, m_j\geqslant0,n_j\in\{0,1\}\}.
\end{align}
By induction, one may show that 
\begin{align}\label{eq:exponents-inequalies}
0\leqslant\kappa\leqslant\frac{1}{2}\leqslant\lambda\leqslant1,\ \ 0\leqslant\nu\leqslant1
\end{align}
for each $(\kappa,\lambda,\nu)\in\fP.$ In practice, the exponent $\nu$ is not quite crucial (since it is at most $1$) and we formally call the elements in $\fP$ as arithmetic exponent pairs. Note that this convention is different from that in \cite{WX21}.

We are now ready to formulate the work in \cite{WX21} on arithmetic exponent pairs.
\begin{lemma}\label{lm:arithmeticexponentpairs}
Keep the notation and convention as above. For each compositely universally amiable trace function $K_q,$ the estimate $\eqref{eq:exponentpair-estimate}$ holds for all
\begin{align*}
(\kappa,\lambda,\nu)\in\fP.
\end{align*}

Amongst many choices, one may take $(\kappa,\lambda,\nu)$ to be
\begin{align*}
B\cdot (0,1,0)&=\Big(\frac{1}{2}, \frac{1}{2}, 1\Big),\\
AB\cdot (0,1,0)&=\Big(\frac{1}{6}, \frac{2}{3}, \frac{1}{6}\Big),\\
BA^2B\cdot (0,1,0)&=\Big(\frac{2}{7}, \frac{4}{7}, \frac{11}{14}\Big),\\
A^2BA^2B\cdot (0,1,0)&=\Big(\frac{1}{20},\frac{33}{40},\frac{1}{20}\Big),
\end{align*}
or
\begin{align*}
A^{k-1}B\cdot (0,1,0)&=\Big(\frac{1}{2^{k+1}-2}, 1-\frac{k}{2^{k+1}-2},\frac{1}{2^{k+1}-2}\Big)
\end{align*}
for each positive integer $k\geqslant3$.
\end{lemma}

As a direct consequence of Lemma \ref{lm:arithmeticexponentpairs}, we have the following estimate for incomplete character sums to smooth moduli.

\begin{lemma}\label{lm:incompletecharactersum-smoothmoduli}
Let $\eta>0$ be a sufficiently small number. 
Let $q$ be an $\eta$-smooth squarefree number and $\chi$ a non-trivial Dirichlet character modulo $q.$ 

For each positive integer $k$ and an interval $I,$ we have
\begin{align*}
\sum_{n\in I}\chi(n)\ll q^\kappa |I|^{\lambda-\kappa}q^{O(\eta)}
\end{align*}
with
\begin{align*}
(\kappa,\lambda)=\Big(\frac{1}{2^{k+1}-2}, 1-\frac{k}{2^{k+1}-2}\Big),
\end{align*}
where the implied constant depends only on $(\eta,k).$
Consequently, for each interval $I$ with $|I|\sim q^\theta$ and $\theta\in]0,1],$ there exists some constant $c=c(\theta)>0$ such that
\begin{align*}
\sum_{n\in I}\chi(n)\ll |I|^{1-c+O(\eta)},
\end{align*}
where the implied constant depends only on $(\eta,\theta).$
\end{lemma}

\proof
The first part follows directly from Lemma \ref{lm:arithmeticexponentpairs}. 
We then have 
\begin{align*}
\sum_{n\in I}\chi(n)\ll |I|^{1-c+O(\eta)}
\end{align*}
with
\begin{align*}
c=\frac{1}{2^{k+1}-2}\Big(k+1-\frac{1}{\theta}\Big).
\end{align*}
Taking $k=[3/\theta]$ finishes the proof.
\endproof

\begin{remark}
We would like to mention that the statement in Lemma $\ref{lm:incompletecharactersum-smoothmoduli}$ is not new. One may refer to Graham and Ringrose \cite[Theorem 5]{GR90}
for upper bounds with explicit dependences on the factorization of $q$.
\end{remark}

The following lemma gives an upper bound for incomplete Kloosterman sums to smooth moduli.
\begin{lemma}\label{lm:incompleteKloostermansum-smoothmoduli}
Let $\eta>0$ be a sufficiently small number. 
Let $q$ be an $\eta$-smooth squarefree number. For each integer $h$ and an interval $I,$ we have
\begin{align*}
\sum_{\substack{n\in I\\(n,q)=1}}
\ue\Big(\frac{h\overline{n}}{q}\Big)
\ll\frac{|I|}{q}(h, q)+q^{\kappa}|I|^{\lambda-\kappa+{O(\eta)}}(h,q)^{\nu},
\end{align*}
where $(\kappa,\lambda,\nu)$ can be any arithmetic exponent pair in $\eqref{eq:arithmeticexponentpairs-set}.$
\end{lemma}

\proof
We first assume that $|I|<q$.
Let $d=(h,q), q_1=q/d$ and $h_1=h/d$.  Hence
\begin{align*}
S:=\sum_{\substack{n\in I\\(n,q)=1}}
\ue\Big(\frac{h\overline{n}}{q}\Big)
&=\sum_{\substack{n\in I\\(n,q_1)=1}}
\ue\Big(\frac{h_1\overline{n}}{q_1}\Big)\mathbf{1}_{(n,d)=1}.
\end{align*}
We are now in a good position to apply Lemma \ref{lm:arithmeticexponentpairs} with
\begin{align*}
(q,\delta)\leftarrow(q_1,d),\ \ 
K_{q_1}: x\mapsto 
\ue\Big(\frac{h_1\overline{x}}{q_1}\Big),\ \ 
W_d: x\mapsto \mathbf{1}_{(x,d)=1},
\end{align*}
getting
\begin{align*}
S\ll q^{\kappa}|I|^{\lambda-\kappa+{O(\eta)}}(h,q)^{\nu}.
\end{align*}

For $|I|\geqslant q,$ the average $S$ would be split into $O(|I|/q)$ copies of Ramanujan sum $S(h,0,q)$ and one incomplete sum of length at most $q$. This incomplete sum can be bounded as above, and the lemma is proven by noting that $|S(h,0,q)|\leqslant (h, q)$.
\endproof

\smallskip

\section{Sieve methods}
\label{sec:sievemethods}

\subsection{Linear sieve} We introduce some conventions and fundamental results on the linear Rosser--Iwaniec sieve.
Let $\cA=(a_n)$ be a finite sequence of integers and $\cP$ a set of prime numbers. Define the sifting function
\begin{align*}S(\cA,\cP,z)=\sum_{(n,P(z))=1}a_n
\quad\text{with}\quad
P(z)=\prod_{p<z, \, p\in\cP}p.
\end{align*}
For squarefree $d$ with all its prime factors belonging to $\cP$, 
we consider the subsequence $\cA_d=(a_n)_{n\equiv 0\bmod d}$ and 
the congruence sum
\begin{align*}
A_d=\sum_{n\equiv 0\!\bmod{d}}a_n.
\end{align*}
In applications to sieve methods, it is expected that $\cA$ equidistributes in the special arithmetic progression 
$n\equiv0\bmod d$: 
There are an appropriate approximation $\sX$ to $A_1$ and 
a multiplicative function $g$ supported on squarefree numbers with all its prime factors belonging to $\cP$ verifying
\[0<g(p)<1\quad(p\in\cP)\]
such that
\par
(a) the remainder 
\[r(\cA,d)=A_d-g(d)\sX\]
is small on average over $d\mid P(z)$;

(b) there exists a constant $L>1$ such that
\begin{align*}
\frac{V(z_1)}{V(z_2)}
\leqslant\frac{\log z_2}{\log z_1}\bigg(1+\frac{L}{\log z_1}\bigg)
\quad\text{with}\quad
V(z)=\prod_{p<z, \, p\in\cP}(1-g(p))
\end{align*}
for $2\leqslant z_1< z_2$.

Let $F$ and $f$ be the continuous solutions to the system
\begin{align}\label{eq:sieve-Ff}
\begin{cases}
sF(s) =2\mathrm{e}^\gamma & \text{for $\,0< s\leqslant 2$},
\\
sf(s)=0                                    & \text{for $\,0<s\le 2$},
\\
(sF(s))'=f(s-1)                         & \text{for $\,s>2$},
\\
(sf(s))'=F(s-1)                         & \text{for $\,s>2$},
\end{cases}
\end{align}
where $\gamma\approx 0.57721$ is the Euler constant.

For $k\ge 1$, denote by $\tau_k(n)$ the number of ways of expressing $n$ as
the product of $k$ positive integers. An arithmetic function $\lambda(d)$ is of {\it level}
$D$ and {\it order} $k$, if
\begin{align*}
\lambda(d)=0
\quad(d>D)
\qquad\text{and}\qquad
|\lambda(d)|\leqslant\tau_k(d)
\quad
(d\leqslant D).
\end{align*}
We say that $\lambda$ is {\it well-factorable}, if for every decomposition $D = D_1D_2$ with $D_1, D_2\geqslant1,$ there exist two arithmetic functions 
$\lambda_1, \lambda_2$ such that
\begin{align*}
\lambda=\lambda_1*\lambda_2\end{align*}
with each $\lambda_j$ of level $D_j$ and order $k$.

We now state the following fundamental result of Iwaniec \cite{Iw80}.

\begin{lemma}\label{lm:Iwaniec}
Let $0<\varepsilon<\tfrac{1}{8}$.
Under the above hypothesis, we have
\begin{align*}
S(\cA,\cP,z)
\leqslant \sX V(z)\Big\{F\Big(\frac{\log D}{\log z}\Big)+E^+\Big\}
+\sum_{t\leqslant T}\sum_{d\mid P(z)}\lambda_t^+(d) r(\cA,d)
\end{align*}
and
\begin{align*}
S(\cA,\cP,z)
\geqslant \sX V(z)\Big\{f\Big(\frac{\log D}{\log z}\Big)-E^-\Big\}
-\sum_{t\leqslant T}\sum_{d\mid P(z)}\lambda_t^-(d) r(\cA,d)
\end{align*}
for all $z\geqslant 2,$
where $F,f$ are given by the system \eqref{eq:sieve-Ff}, 
$T$ depends only on $\varepsilon,$ 
$\lambda_t^\pm(d)$ are well-factorable of level $D$ and order $1,$ 
and $E^\pm\ll\varepsilon +\varepsilon^{-8} \mathrm{e}^L(\log D)^{-1/3}.$
\end{lemma}

From the system \eqref{eq:sieve-Ff}, one can show that
\begin{align}\label{eq:sieve-F}
F(s) = \begin{cases}
\displaystyle\frac{2\ue^{\gamma}}{s}, & \text{for }1\leqslant s\leqslant3,
\\\noalign{\vskip 1mm}
\displaystyle\frac{2\ue^{\gamma}}{s}\Big(1+\int_2^{s-1}\frac{\log(t-1)}{t}\ud t\Big), & \text{for }3\leqslant s\leqslant5,
\end{cases}
\end{align}
and
\begin{align}\label{eq:sieve-f}
f(s)=
\begin{cases}
0, & \text{for }0\leqslant s\leqslant2,
\\ \noalign{\vskip 1mm}
\displaystyle\frac{2\ue^{\gamma}}{s}\log(s-1), & \text{for }2\leqslant s\leqslant4.
\end{cases}
\end{align}
Moreover, $F(s)$ decreases (resp. $f(s)$ increases) rapidly to $1$ as $s\rightarrow+\infty.$ In applications, we would like to explore the admissible $D$ as large as possible such that the total remainder
\begin{align*}
\sum_{t\leqslant T}\sum_{d\mid P(z)}\lambda_t^\pm(d) r_d(\cA)
\end{align*}
is still under control. In principle, a larger $D$ would produce a better upper/lower bound for $S(\cA,\cP,z)$ upon a suitable choice for $z.$ The merit of Lemma \ref{lm:Iwaniec} lies in the well-factorization of $\lambda_t^\pm$, so that one can estimate the resultant bilinear forms more effectively to produce a better admissible $D$ (compared with earlier sieve methods before Iwaniec). This bilinear structure is very crucial in many existing applications of Iwaniec's linear sieve.

We also quote the following Buchstab identity, which can be verified directly.
\begin{lemma}\label{lm:Buchstab}
For all $z>w\geqslant2,$ we have
\begin{align*}
S(\cA,\cP,z)=S(\cA,\cP,w)-\sum_{\substack{w\leqslant p<z\\p\in\cP}}S(\cA_p,\cP,p).
\end{align*}
\end{lemma}

\subsection{Sifting primes in $\pi(x;q,a)$}
Let $y=x(\log x)^{-3}$ and introduce a smooth function $\varPhi$ which is supported on $[y,x+y]$, satisfying
\begin{align*}
\begin{cases}
\varPhi(t)=1\ \ &\text{for } t\in[2y,x],\\
\varPhi(t)\geqslant0\ \ &\text{for } t\in\bR,\\
\varPhi^{(j)}(t)\ll_j y^{-j} &\text{for all } j\geqslant0,
\end{cases}
\end{align*}
and $\widehat{\varPhi}(0)=\int_\bR \varPhi(t)\ud t=x.$
Therefore,
\begin{align*}
\pi(x; q, a)=\pi^*(x; q, a)+O(q^{-1}x(\log x)^{-2}),
\end{align*}
where
\begin{align*}
\pi^*(x;q,a)=\sum_{p\equiv a \bmod q}\varPhi(p).
\end{align*}
In view of the definition of $C(\varpi)$ in \eqref{eq:BT}, the above error $O(q^{-1}x(\log x)^{-2})$ is acceptable, and it suffices to derive upper bounds for the smooth sum $\pi^*(x; q, a).$ Denote by $P(t)$ the product of primes up to $t$ and not dividing $q$.
On one hand, we have
\begin{align}\label{eq:Pi*(x;q,a)-upperbound}
\pi^*(x;q,a)\leqslant\sum_{\substack{n\equiv a \bmod q\\(n,P(z))=1}}\varPhi(n)+O(1+q^{-1}x^{\frac{1}{2}})
\end{align}
for all $2\leqslant z<\sqrt{x}.$
By Buchstab identity (Lemma \ref{lm:Buchstab}), we also have
\begin{align}\label{eq:Pi*(x;q,a)-viaBuchstab}
\pi^*(x;q,a)\leqslant
\sum_{\substack{n\equiv a \bmod q\\(n,P(w))=1}}\varPhi(n)-\sum_{w\leqslant p<z}\sum_{\substack{pn\equiv a \bmod q\\(n,P(p))=1}}\varPhi(pn)+O(1+q^{-1}x^{\frac{1}{2}})
\end{align}
for all $2\leqslant w< z<\sqrt{x}.$

For $(d, q)=1$ we consider the congruence sum
\begin{align}\label{eq:A(d)}
A(d)
:= \sum_{\substack{n\equiv a\bmod q\\d\mid n}} \varPhi(n)=
\sum_{n\equiv a\overline{d}\bmod q} \varPhi(dn).
\end{align}
Define
\begin{align}\label{eq:r(d)}
r(d)
&:=A(d)-\frac{x}{dq}.
\end{align}
Let $(\lambda_d)$ be a linear upper bound sieve of level $D$, so that $1*\mu\leqslant1*\lambda$. 
Then 
\begin{align}\label{eq:upperboundsieve}
\pi^*(x;q,a)
&\leqslant \sum_{\substack{d\leqslant D\\d\mid P(z)}} \lambda_d A(d)+O(q^{-1}x(\log x)^{-2})\\
&\leqslant \frac{x}{q}V(z)F\Big(\frac{\log D}{\log z}\Big)
+ \sum_{\substack{d\leqslant D\\d\mid P(z)}} \lambda_d r(d)+O(q^{-1}x(\log x)^{-2}),
\end{align}
where $F$ is defined by \eqref{eq:sieve-Ff} and
\begin{align*}
V(z) 
:= \prod_{p<z, \, p\nmid q} \Big(1-\frac{1}{p}\Big)
= \frac{q}{\varphi(q)} \cdot \frac{\ue^{-\gamma}}{\log z}\Big\{1+O\Big(\frac{1}{\log x}\Big)\Big\}.
\end{align*}
Thanks to Lemma \ref{lm:Iwaniec}, we may choose well-factorable remainder terms in the above application of the linear Rosser--Iwaniec sieve.
We thus conclude the following result.

\begin{proposition}\label{prop:sieveinequality}
Let $\varepsilon>0$ sufficiently small and $D<x^{1-\varepsilon}.$ For any $z$ with  $D^{\frac{1}{3}}\leqslant z<\sqrt{x},$ we have
\begin{align}\label{eq:pi(x;q,a)-sievebound} 
\pi(x; q, a)
\leqslant \{2+O(\varepsilon)\}\frac{x}{\varphi(q)\log D}
+\sum_{t\leqslant T(\varepsilon)}\sum_{\substack{d\leqslant D\\d\mid P(z)}}\lambda_t(d) r(d),\end{align}
where $T(\varepsilon)$ depends only on $\varepsilon,$ $\lambda_t$ is well-factorable of level $D,$ and order $1,$ and $r(d)$ is defined by $\eqref{eq:A(d)}$ and $\eqref{eq:r(d)}.$
\end{proposition}

In place of \eqref{eq:Pi*(x;q,a)-upperbound}, we may obtain the following upper bound for $\pi(x;q,a)$ by virtue of \eqref{eq:Pi*(x;q,a)-viaBuchstab}.

\begin{proposition}\label{prop:Buchstabsieveinequality}
Let $\varepsilon>0$ sufficiently small and $D,D_p<x^{1-\varepsilon}.$ For any $2\leqslant w< z<\sqrt{x},$ we have
\begin{align*}
\pi(x;q,a)
&\leqslant (\ue^{-\gamma}+O(\varepsilon))\frac{x}{\varphi(q)}\Big\{\frac{1}{\log w}F\Big(\frac{\log D}{\log w}\Big)-\sum_{w\leqslant p<z}\frac{1}{p\log p}f\Big(\frac{\log D_p}{\log p}\Big)\Big\}\\
&\ \ \ \ +\sum_{t\leqslant T^+(\varepsilon)}\sum_{\substack{d\leqslant D\\d\mid P(w)}}\lambda_t^+(d) r(d)+\sum_{t\leqslant T^-(\varepsilon)}\sum_{w\leqslant p<z}\sum_{\substack{d\leqslant D_p\\d\mid P(p)}}\lambda_t^-(d) r(pd),
\end{align*}
where $T^\pm(\varepsilon)$ depends only on $\varepsilon,$ and $\lambda_t^+,\lambda_t^-$ are well-factorable of level $D,D_p,$ respectively.
\end{proposition}

\smallskip

\section{Proofs of Theorems \ref{thm:BT-generalmoduli<1/2} and \ref{thm:BT-smoothmoduli<1/2}}
\label{sec:BT-generalmoduli<1/2}

Keep in mind that $q< x^{\frac{1}{2}-}$ throughout this section.
Recall the definitions $\eqref{eq:A(d)}$ and $\eqref{eq:r(d)}.$ Fix $w\leqslant P<z$. Let $\lambda^+$ be a well-factorable function of level $D$, order 1, and $\lambda^-$ a well-factorable function of level $D^*$, order 1 with $D,D^*<x^{1-\varepsilon}$.
In view of Proposition \ref{prop:Buchstabsieveinequality}, we would like to prove that
\begin{align}\label{eq:R+(D)-expectedbound}
R^+(D):=\sum_{d\leq D}\lambda^+(d) r(d)
\ll \frac{x}{q(\log x)^{2024}}
\end{align}
and
\begin{align}\label{eq:R-(P,D*)-expectedbound}
R^-(P,D^*):=\sum_{p\sim P}\sum_{d\leq D^*}\lambda^-(d) r(pd)
\ll \frac{x}{q(\log x)^{2024}}
\end{align} 
with $D,D^*$ as large as possible.

\subsection{Approaching by multiplicative characters}\label{subsection:Approaching by multiplicative characters} 
From orthogonality of Dirichlet characters, it follows that
\begin{align*}
r(d)
&=\frac{1}{\varphi(q)}\sum_{\substack{\chi\bmod{q}\\
\chi\neq \chi_0}}\overline{\chi}(a)\sum_{n\in\bZ} \varPhi(dn)\chi(dn)+O(1).
\end{align*}
By Mellin inversion, we may write
\begin{align*}
r(d)
&=\frac{1}{\varphi(q)}\frac{1}{2\pi i}\int_{(2)}\widetilde{\varPhi}(s)d^{-s}\sum_{\substack{\chi\bmod{q}\\
\chi\neq \chi_0}}\overline{\chi}(a)\chi(d)L(s,\chi)\ud s+O(1),
\end{align*}
where $L(s,\chi)$ is the Dirichlet $L$-function.
In view of the decay
\begin{align*}
\widetilde{\varPhi}(s)
&\ll x^{\Re(s)}(1+x^{-1}y|s|)^{-A}
\end{align*}
for any $A\geqslant0$, we shift the integral to $\Re(s)=1/2$ and sum over $d\leqslant D,$ getting
\begin{align*}
R^+(D)
&=\frac{1}{\varphi(q)}\frac{1}{2\pi}\int_{\bR}\widetilde{\varPhi}(\tfrac{1}{2}+it)\sum_{\substack{\chi\bmod{q}\\
\chi\neq \chi_0}}\overline{\chi}(a)\cD(\tfrac{1}{2}+it,\chi)L(\tfrac{1}{2}+it,\chi)\ud t+O(1),
\end{align*}
where $\cD(s,\chi)$ is a Dirichlet polynomial defined by
\begin{align*}
\cD(s,\chi)=\sum_{d\leqslant D}\lambda^+(d)\chi(d) d^{-s}.
\end{align*}
Note that
\begin{align*}
\int_{\bR}|\widetilde{\varPhi}(\tfrac{1}{2}+it)|(1+|t|)^4\ud t
&\ll x^{\frac{1}{2}}(\log x)^{12},
\end{align*}
it suffices to prove that for $\Re(s)=\frac{1}{2},$
\begin{align}\label{eq:moments-expectedbound}
\sum_{\substack{\chi\bmod{q}\\
\chi\neq \chi_0}}\overline{\chi}(a)\cD(s,\chi)L(s,\chi)
\ll |s|^4x^{\frac{1}{2}}(\log x)^{-2036}.
\end{align}

Due to the well-factorization of $\lambda^+$, we write 
\begin{align}\label{eq:sieveweight-convolution} 
\lambda^+=\balpha*\bbeta,
\end{align}
where $\balpha$ and $\bbeta$ are of level $M$, $N$, respectively, and 
\begin{align*}
M,N\geqslant1,\ \ \ \ MN=D 
\quad\text{and}\quad
\|\balpha\|_\infty,\|\bbeta\|_\infty \leq 1.
\end{align*}
By dyadic device, we henceforth assume that $\balpha$ and $\bbeta$ are supported in $]M,2M]$, $]N,2N]$, respectively. 
Put
\begin{align*}
\cM(s,\chi)
&=\sum_{m\sim M}\alpha_m\chi(m)m^{-s},\\ 
\cN(s,\chi)
&=\sum_{n\sim N}\beta_n\chi(n)n^{-s},\\ 
\cP(s,\chi)
&=\sum_{p\sim P}\chi(p)p^{-s}.
\end{align*}
Hence the proof  \eqref{eq:R+(D)-expectedbound} can be reduced to proving that
\begin{align}\label{eq:moments-expectedbound-L}
\sum_{\substack{\chi\bmod{q}\\
\chi\neq \chi_0}}|\cM(s,\chi)\cN(s,\chi)L(s,\chi)|\ll |s|^4x^{\frac{1}{2}}(\log x)^{-2036}
\end{align}
for $\Re(s)=\frac{1}{2}.$
Recall that $\varPhi$ is supported in $[y,x+y]$. For $m\sim M,n\sim N,$ it is clear that $\varPhi(mnl)$ vanishes unless $l\asymp x/(MN)$. As argued above, an alternative approach to \eqref{eq:R+(D)-expectedbound} is proving
\begin{align}\label{eq:moments-expectedbound-cL}
\sum_{\substack{\chi\bmod{q}\\
\chi\neq \chi_0}}|\cM(s,\chi)\cN(s,\chi)\cL(s,\chi)|\ll |s|^4x^{\frac{1}{2}}(\log x)^{-2036},
\end{align}
where $\Re(s)=\frac{1}{2},$ and 
\begin{align*}
\cL(s,\chi)=\sum_{l\in\bZ}\phi\Big(\frac{l}{L}\Big)\chi(l)l^{-s},\ \ \ L=\frac{x}{MN}.
\end{align*}
Here $\phi$ is a fixed smooth function which dominates the indicator function of $[-1,1].$

In a similar way, the proof 
\eqref{eq:R-(P,D*)-expectedbound} can be reduced to proving that
\begin{align}\label{eq:moments-Buchstabexpectedbound-L}
\sum_{\substack{\chi\bmod{q}\\
\chi\neq \chi_0}}|\cP(s,\chi)\cM(s,\chi)\cN(s,\chi)L(s,\chi)|\ll |s|^4x^{\frac{1}{2}}(\log x)^{-2036}
\end{align}
or
\begin{align}\label{eq:moments-Buchstabexpectedbound-cL}
\sum_{\substack{\chi\bmod{q}\\
\chi\neq \chi_0}}|\cP(s,\chi)\cM(s,\chi)\cN(s,\chi)\cL(s,\chi)|\ll |s|^4x^{\frac{1}{2}}(\log x)^{-2036},
\end{align}
where $\Re(s)=\frac{1}{2}.$

\subsection{Concluding Theorem \ref{thm:BT-generalmoduli<1/2}}\label{sec:concludingBT-generalmoduli<1/2,underconjecture-generalmoduli}

We would like to estimate the error terms in Proposition \ref{prop:Buchstabsieveinequality}. They correspond to the desired inequalities in \eqref{eq:R+(D)-expectedbound} and \eqref{eq:R-(P,D*)-expectedbound}.
The first one will be dealt with using Lemma \ref{lm:weightedL-secondmoment-generalcoefficient}, and in the treatment to the second one we utilize the convolution structures of the coefficients in Dirichlet polynomials, so that Lemma \ref{lm:weightedL-secondmoment-specialcoefficient} applies. Assume that $x^{\frac{9}{20}}\leqslant q< x^{\frac{1}{2}-}$ since the case $q<x^{\frac{9}{20}}$ has been dealt with by Iwaniec \cite{Iw82}.

Denote by $\Xi_1$ the LHS of \eqref{eq:moments-expectedbound-cL}.
By Cauchy--Schwarz inequality and Lemma \ref{lm:Dirichletpolynomial-secondmoment}, we have
\begin{align*}
\Xi_1^2\ll (M+q)\sum_{\substack{\chi\bmod{q}\\
\chi\neq \chi_0}}|\cN(s,\chi)\cL(s,\chi)|^2,
\end{align*}
where $MN=D$ and $M,N\geqslant1.$
From Theorem \ref{thm:countingcongruence-generalcoefficient}, it follows that
\begin{align*}
\Xi_1^2\ll |s|^{\frac{7}{2}}(M+q)(q+N^{\frac{33}{20}}q^{\frac{3}{20}}+N^{\frac{15}{8}})q^{\varepsilon}.
\end{align*}
Note that the choice
\begin{align}\label{eq:MN,Buchstabupperlevel}
M= xq^{-1-\varepsilon},\ \ N=q^{\frac{17}{33}-\varepsilon}
\end{align}
is valid to guarantee the bound \eqref{eq:moments-expectedbound-cL}, from which we obtain a level $D=xq^{-\frac{16}{33}-2\varepsilon}$ in \eqref{eq:R+(D)-expectedbound}.

If appealing to Proposition \ref{prop:sieveinequality} directly, we may choose $C(\varpi)=66/(33-16\varpi)$ in \eqref{eq:BT}, i.e., one can take $C^*(\varpi)=0$ in \eqref{eq:C(varpi):varpi<1/2}. The positive choices of $C^*(\varpi)$ in \eqref{eq:C*(varpi):varpi<1/2} and \eqref{eq:C*(varpi):varpi<1/2,primemoduli} are due to the task of Proposition \ref{prop:Buchstabsieveinequality}, in which the trilinear structure in the error term becomes very crucial.

Denote by $\Xi_2$ the LHS of \eqref{eq:moments-Buchstabexpectedbound-cL}.
It is our next task to apply Lemma \ref{lm:weightedL-secondmoment-specialcoefficient}, for which it is necessary to check the supports of $\balpha=(\alpha_m)$ and $\bbeta=(\beta_n)$.
In fact, we can write the LHS of \eqref{eq:R+(D)-expectedbound} as
\begin{align*}
 \sum_{p\sim P}\sum_{\substack{d\leqslant D^*\\d\mid P(p)}} \lambda_d^- r(pd) =\sum_{p\sim P}\sum_{\substack{d\leqslant D^*\\d\mid P(p)\\ \omega(d)\leqslant (\log x)^{1/3}}} \lambda_d^-  r(pd)+\sum_{p\sim P}\sum_{\substack{d\leqslant D^*\\d\mid P(p)\\ \omega(d)> (\log x)^{1/3}}} \lambda_d^-  r(pd)
=\pi_1^*+\pi_2^*.
\end{align*}
Recall that $r(d)$ is defined by \eqref{eq:r(d)}. Hence 
\begin{align*}
\pi_2^*
&\ll \sum_{p\sim P}\sum_{\substack{d\leqslant D^*\\ \omega(d)> (\log x)^{1/3}}}|\lambda_d^-|\sum_{\substack{n\equiv a\bmod q\\ pd\mid n}} \varPhi(n)+\frac{x}{q}\sum_{p\sim P}\sum_{\substack{d\leqslant D^*\\ \omega(d)> (\log x)^{1/3}}}\frac{|\lambda_d^-|}{pd}\\
&\ll 2^{-(\log x)^{1/3}}\sum_{d\leqslant 2PD^*}2^{\omega(d)}\sum_{\substack{n\equiv a\bmod q\\d\mid n}} \varPhi(n)+\frac{2^{-(\log x)^{1/3}}x}{q}\sum_{d\leqslant 2PD^*}\frac{2^{\omega(d)}}{d}\\
&\ll 2^{-(\log x)^{1/3}}\sum_{\substack{n\leqslant 2x\\ n\equiv a\bmod q}} \tau^2(n)+\frac{2^{-(\log x)^{1/3}}x}{q}\sum_{d\leqslant 2PD^*}\frac{2^{\omega(d)}}{d}.
\end{align*}
Recalling that $P<z<\sqrt{x}$ and $D^*<x^{1-\varepsilon}$, 
the last expression is at most $O(x/(q\log^{2024}x))$ by Shiu \cite[Theorem 2]{Sh80}, provided that $q<x^{0.99}.$ Thus we have
\begin{align*}
 \sum_{p\sim P}\sum_{\substack{d\leqslant D^*\\d\mid P(p)}} \lambda_d^- r(pd)=\sum_{p\sim P}\sum_{\substack{d\leqslant D^*\\d\mid P(p)\\ \omega(d)\leqslant (\log x)^{1/3}}} \lambda_d^-  r(pd)+O\Big(\frac{x}{q(\log x)^{2024}}\Big).
\end{align*}
Hence it suffices to consider those $\lambda_d^-$ supported on $\{d\leqslant D^*:\mu^2(d)=1,\omega(d)\leqslant (\log x)^{1/3}\}.$
With such modifications to the definitions of $\cM(s,\chi)$ and $\cN(s,\chi)$, it follows from Lemma \ref{lm:weightedL-secondmoment-specialcoefficient} that
\begin{align*}
\Xi_2^2
&\ll |s|^{\frac{9}{2}}(M+q)P^{\frac{1}{2}}N^2q^{\frac{\vartheta+1}{2}}(Pq^{-\frac{1}{2}}+1)q^{\varepsilon}+|s| (M+q)q^{1+\varepsilon}+|s|^{\frac{9}{2}}(M+q)Pq^{\frac{1}{2}+\varepsilon}
\end{align*}
for any $\varepsilon>0,$ where $MN=D^*$ and $M,N\geqslant1.$ To guarantee the bound \eqref{eq:moments-Buchstabexpectedbound-cL},
we hope
\begin{align*}
MN^2&\leqslant x^{1-\varepsilon}q^{-\frac{\vartheta+1}{2}}P^{-\frac{1}{2}}(Pq^{-\frac{1}{2}}+1)^{-1},\\
N^2&\leqslant x^{1-\varepsilon}q^{-\frac{\vartheta+3}{2}}P^{-\frac{1}{2}}(Pq^{-\frac{1}{2}}+1)^{-1},\\
M&\leqslant x^{1-\varepsilon}q^{-1},\\
M&\leqslant x^{1-\varepsilon}q^{-\frac{1}{2}}P^{-1}.
\end{align*}
The above inequalities also require an additional condition
$P\leqslant x^{1-\varepsilon}q^{-\frac{3}{2}}$, which we henceforth assume.
In particular, we may take
\begin{align*}
M=x^{1-\varepsilon}q^{-1}\min\{1,P^{-1}q^{\frac{1}{2}}\},\ \ N=x^{-\varepsilon}q^{\frac{1-\vartheta}{4}}P^{-\frac{1}{4}}.
\end{align*}
Hence \eqref{eq:R-(P,D*)-expectedbound} holds with
\begin{align*}
D^*=x^{1-2\varepsilon}q^{-\frac{3+\vartheta}{4}}P^{-\frac{1}{4}}\min\{1,P^{-1}q^{\frac{1}{2}}\},
\end{align*}
so that one may take
\begin{align*}
D_p=x^{1-2\varepsilon}q^{-\frac{3+\vartheta}{4}}p^{-\frac{1}{4}}\min\{1,p^{-1}q^{\frac{1}{2}}\}
\end{align*}
in Proposition \ref{prop:Buchstabsieveinequality}. This gives
\begin{align*}
\pi(x;q,a)
&\leqslant (\ue^{-\gamma}+O(\varepsilon))\frac{x}{\varphi(q)}\Big\{\frac{1}{\log w}F\Big(\frac{\log(xq^{-\frac{16}{33}})}{\log w}\Big)\\
&\ \ \ \ \ -\sum_{w\leqslant p<z}\frac{1}{p\log p}f\Big(\frac{\log(xq^{-\frac{3+\vartheta}{4}}p^{-\frac{1}{4}}\min\{1,p^{-1}q^{\frac{1}{2}}\})}{\log p}\Big)\Big\}
\end{align*}
for $z\leqslant x^{1-\varepsilon}q^{-\frac{3}{2}}.$

We would like to specialize $z,w$ such that the argument in the above $F(\cdot)$ falls into $[3,5]$ and that in $f(\cdot)$ falls into $[2,4].$ To do so, we put
\begin{align*}
z=x^{\kappa_1-\varepsilon},\ \ w=x^{\kappa_2},\ \ \kappa_1=\min\Big\{1-\frac{3}{2}\varpi,~\frac{8-(1+2\vartheta)\varpi}{28}\Big\},\ \ \kappa_2=\frac{1}{5}\Big(1-\frac{16}{33}\varpi\Big).
\end{align*}
By the Prime Number Theorem and partial summation,
\begin{align*}
\pi(x;q,a)
&\leqslant (\ue^{-\gamma}+O(\varepsilon))\frac{x}{\varphi(q)}\Big\{\frac{1}{\log w}F\Big(\frac{\log(xq^{-\frac{16}{33}})}{\log w}\Big)\\
&\ \ \ \ \ -\int_w^z\frac{1}{t(\log t)^2}f\Big(\frac{\log(xq^{-\frac{3+\vartheta}{4}}t^{-\frac{1}{4}}\min\{1,t^{-1}q^{\frac{1}{2}}\})}{\log t}\Big)\ud t\Big\}\\
&= (\ue^{-\gamma}+O(\varepsilon))\frac{c\cdot x}{\varphi(q)\log x},
\end{align*}
where
\begin{align*}
c&=\frac{1}{\kappa_2}F\Big(\frac{1-\frac{16}{33}\varpi}{\kappa_2}\Big)-\int_{\kappa_2}^{\varpi/2}f\Big(\frac{4-(3+\vartheta)\varpi)}{4u}-\frac{1}{4}\Big)\frac{\ud u}{u^2}-\int_{\varpi/2}^{\kappa_1}f\Big(\frac{4-(1+\vartheta)\varpi}{4u}-\frac{5}{4}\Big)\frac{\ud u}{u^2}.
\end{align*}
Recalling the explicit expressions \eqref{eq:sieve-F} and \eqref{eq:sieve-f} for $F$ and $f,$ we write
\begin{align*}
c&=\frac{2\ue^{\gamma}}{1-\frac{16}{33}\varpi}\Big(1+\int_2^4\frac{\log(t-1)}{t}\ud t\Big)-\frac{4\cdot 2\ue^{\gamma}}{4-(3+\vartheta)\varpi}\int_{\frac{8-(7+2\vartheta)\varpi}{4\varpi}}^{\frac{5(4-(3+\vartheta)\varpi)}{4(1-\frac{16}{33}\varpi)}-\frac{1}{4}}\frac{\log(t-1)}{t}\ud t\\
&\ \ \ \ -\frac{4\cdot 2\ue^{\gamma}}{4-(1+\vartheta)\varpi}\int_{\max\{\frac{4-(1+\vartheta)\varpi}{2(2-3\varpi)}-\frac{5}{4},~2\}}^{\frac{8-(7+2\vartheta)\varpi}{4\varpi}}\frac{\log(t-1)}{t}\ud t.
\end{align*}
Hence we may take
\begin{align*}
C(\varpi)&=\frac{66}{33-16\varpi}-C^*(\varpi)
\end{align*}
in \eqref{eq:BT}, where
\begin{align*}
C^*(\varpi)&=-\frac{66}{33-16\varpi}\int_2^4\frac{\log(t-1)}{t}\ud t+\frac{8}{4-(3+\vartheta)\varpi}\int_{\frac{8-(7+2\vartheta)\varpi}{4\varpi}}^{\frac{165(4-(3+\vartheta)\varpi)}{4(33-16\varpi)}-\frac{1}{4}}\frac{\log(t-1)}{t}\ud t\\
&\ \ \ \ +\frac{8}{4-(1+\vartheta)\varpi}\int_{\max\{\frac{4-(1+\vartheta)\varpi}{2(2-3\varpi)}-\frac{5}{4},~2\}}^{\frac{8-(7+2\vartheta)\varpi}{4\varpi}}\frac{\log(t-1)}{t}\ud t.
\end{align*}
This proves the first part of Theorem \ref{thm:BT-generalmoduli<1/2} readily.

We now assume that $q$ is a large prime, and denote by $\Theta$ the LHS of \eqref{eq:moments-Buchstabexpectedbound-L}, i.e.,
\begin{align*}
\Theta=\sum_{\substack{\chi\bmod{q}\\
\chi\neq \chi_0}}|\cM(s,\chi)\cP(s,\chi)\cN(s,\chi)L(s,\chi)|.
\end{align*}
By H\"older's inequality, we have
\begin{align}\label{eq:Theta-Holder}
\Theta^4\ll \Big(\sum_{\substack{\chi\bmod{q}\\
\chi\neq \chi_0}}|\cM(s,\chi)|^2\Big)^2\Big(\sum_{\substack{\chi\bmod{q}\\
\chi\neq \chi_0}}|\cP(s,\chi)|^4\Big)\Big(\sum_{\substack{\chi\bmod{q}\\
\chi\neq \chi_0}}|\cN(s,\chi)L(s,\chi)|^4\Big).
\end{align}
From Lemma \ref{lm:weightedL-fourthmoment} we infer
\begin{align*}
\sum_{\substack{\chi\bmod{q}\\
\chi\neq \chi_0}}|\cN(s,\chi)L(s,\chi)|^4\ll |s|^{16} (N^4q^{-\frac{1}{2}}+1)q^{1+\varepsilon}.
\end{align*}
Together with Lemma \ref{lm:largesievebound-Dirichletpolynomial}, we find
\begin{align*}
\Theta\ll |s|^4 (M+q)^{\frac{1}{2}}(P^2+q)^{\frac{1}{4}}(Nq^{-\frac{1}{8}}+1)q^{\frac{1}{4}+\varepsilon}.
\end{align*}
Hence \eqref{eq:moments-Buchstabexpectedbound-L} is valid as long as
\begin{align*}
M=x^{1-\varepsilon}q^{-1}\min\{1,P^{-1}q^{\frac{1}{2}}\},\ \ N=q^{\frac{1}{8}},
\end{align*}
if $q<x^{\frac{1}{2}-}$, $P\leqslant x^{1-\varepsilon}q^{-\frac{3}{2}}$. Upon the above choice, we may take
\begin{align*}
D_p=x^{1-\varepsilon}q^{-\frac{7}{8}}\min\{1,p^{-1}q^{\frac{1}{2}}\}
\end{align*}
in Proposition \ref{prop:Buchstabsieveinequality}.

As argued above, we arrive at
\begin{align*}
\pi(x;q,a)
&\leqslant (\ue^{-\gamma}+O(\varepsilon))\frac{x}{\varphi(q)}\Big\{\frac{1}{\log w}F\Big(\frac{\log(xq^{-\frac{16}{33}})}{\log w}\Big)\\
&\ \ \ \ \ -\int_w^z\frac{1}{t(\log t)^2}f\Big(\frac{\log(xq^{-\frac{7}{8}}\min\{1,t^{-1}q^{\frac{1}{2}}\})}{\log t}\Big)\ud t\Big\}
\end{align*}
with 
\begin{align*}
z=x^{\min\{1-\frac{3}{2}\varpi,~\frac{1}{3}-\frac{1}{8}\varpi\}-\varepsilon},\ \ \ w=x^{\frac{1}{5}(1-\frac{16}{33}\varpi)}.
\end{align*}
Hence
\begin{align*}
\pi(x;q,a)
&\leqslant (\ue^{-\gamma}+O(\varepsilon))\frac{x}{\varphi(q)\log x}\Big\{\frac{5F(5)}{1-\frac{16}{33}\varpi}-\int_{\frac{1}{5}(1-\frac{16}{33}\varpi)}^{\frac{1}{2}\varpi}f\Big(\frac{1-\frac{7}{8}\varpi}{u}\Big)\frac{\ud u}{u^2}\\
&\ \ \ \ -\int_{\frac{1}{2}\varpi}^{\min\{1-\frac{3}{2}\varpi,~\frac{1}{3}-\frac{1}{8}\varpi\}}f\Big(\frac{1-\frac{3}{8}\varpi}{u}-1\Big)\frac{\ud u}{u^2}\Big\}.
\end{align*}
This establishes the second part of Theorem \ref{thm:BT-generalmoduli<1/2} by substituting the explicit expressions \eqref{eq:sieve-F} and \eqref{eq:sieve-f} for $F$ and $f.$

\subsection{Concluding Theorem \ref{thm:BT-smoothmoduli<1/2}}
In this subsection, we assume $q$ is squarefree and $\eta$-smooth, as required in Theorem \ref{thm:BT-smoothmoduli<1/2}. 
The application of Lemma \ref{lm:largevalue-Huxley} will depend on the sizes of corresponding Dirichlet polynomials, so that we divide all non-trivial characters $\bmod q$  into eight sets. Precisely, we define 
\begin{align*}
\Gamma_{1}^{+}=\{\chi\neq \chi_{0}\bmod{q}:|\cM(s,\chi)|\geqslant q^{\frac{1}{4}},|\cN(s,\chi)|\geqslant q^{\frac{1}{4}},|\cL(s,\chi)|\geqslant q^{\frac{1}{8}}\},\\
\Gamma_{1}^{-}=\{\chi\neq \chi_{0}\bmod{q}:|\cM(s,\chi)|\geqslant q^{\frac{1}{4}},|\cN(s,\chi)|\geqslant q^{\frac{1}{4}},|\cL(s,\chi)|\leqslant q^{\frac{1}{8}}\},\\
\Gamma_{2}^{+}=\{\chi\neq \chi_{0}\bmod{q}:|\cM(s,\chi)|\geqslant q^{\frac{1}{4}},|\cN(s,\chi)|\leqslant q^{\frac{1}{4}},|\cL(s,\chi)|\geqslant q^{\frac{1}{8}}\},\\
\Gamma_{2}^{-}=\{\chi\neq \chi_{0}\bmod{q}:|\cM(s,\chi)|\geqslant q^{\frac{1}{4}},|\cN(s,\chi)|\leqslant q^{\frac{1}{4}},|\cL(s,\chi)|\leqslant q^{\frac{1}{8}}\},\\
\Gamma_{3}^{+}=\{\chi\neq \chi_{0}\bmod{q}:|\cM(s,\chi)|\leqslant q^{\frac{1}{4}},|\cN(s,\chi)|\geqslant q^{\frac{1}{4}},|\cL(s,\chi)|\geqslant q^{\frac{1}{8}}\},\\
\Gamma_{3}^{-}=\{\chi\neq \chi_{0}\bmod{q}:|\cM(s,\chi)|\leqslant q^{\frac{1}{4}},|\cN(s,\chi)|\geqslant q^{\frac{1}{4}},|\cL(s,\chi)|\leqslant q^{\frac{1}{8}}\},\\
\Gamma_{4}^{+}=\{\chi\neq \chi_{0}\bmod{q}:|\cM(s,\chi)|\leqslant q^{\frac{1}{4}},|\cN(s,\chi)|\leqslant q^{\frac{1}{4}},|\cL(s,\chi)|\geqslant q^{\frac{1}{8}}\},\\
\Gamma_{4}^{-}=\{\chi\neq \chi_{0}\bmod{q}:|\cM(s,\chi)|\leqslant q^{\frac{1}{4}},|\cN(s,\chi)|\leqslant q^{\frac{1}{4}},|\cL(s,\chi)|\leqslant q^{\frac{1}{8}}\}.
\end{align*}
For $ 1\leqslant k\leqslant 4$, and $s=\frac{1}{2} +i t$ with $t\in\bR,$ we define 
\begin{align*}
S_{k}^{\pm}
&=   
\sum_{\substack{\chi\in \Gamma_{k}^{\pm}}}|\cM(s,\chi)\cN(s,\chi)L(s,\chi)|,\\
T_{k}^{\pm}
&=   
\sum_{\substack{\chi\in \Gamma_{k}^{\pm}}}|\cM(s,\chi)\cN(s,\chi)\cL(s,\chi)|.
\end{align*} 
Hence it remains to prove that for all $*\in\{+,-\}$ and $1\leqslant k\leqslant 4, $
\begin{align}\label{eq:moments-expectedbound-minST}
\min\{S_{k}^{*},~T_{k}^{*}\}
\ll |s|^3x^{\frac{1}{2}}(\log x)^{-2036} .
\end{align}

To estimate these sums, we first develop some lemmas on mean values of $\cM(s,\chi),$ $\cN(s,\chi),$ $\cL(s,\chi)$ and $L(s,\chi)$ for $s=\frac{1}{2}+it,$  which would be denoted by $\cM,\cN,\cL$ and $L$ for convenience.

The following lemma is a direct consequence of Lemma \ref{lm:Dirichletpolynomial-secondmoment}.

\begin{lemma}\label{lm:largesievebound-Dirichletpolynomial} 
For $1\leqslant k\leqslant 4$, we have
\begin{align*} 
\sum_{\chi \in \Gamma_{k}^{\pm}} |\cM|^2
\ll M+q
\end{align*}
and
\begin{align*}
\sum_{\chi \in \Gamma_{k}^{\pm}} |\cN|^{2}
\ll N+q.
\end{align*}
\end{lemma}

\begin{lemma}\label{lm:largevaluebound-Dirichletpolynomial} 
For $ \Gamma \in \{ \Gamma_{1}^{+},  \Gamma_{1}^{-},  \Gamma_{2}^{+},  \Gamma_{2}^{-} \} ,$ we have
\begin{equation}\label{eq:largevaluebound-Msecond} 
\sum_{\chi \in \Gamma} |\cM|^{2}
\ll Mq^\varepsilon.
\end{equation}
For $ \Gamma \in \{ \Gamma_{1}^{+},  \Gamma_{1}^{-},  \Gamma_{3}^{+},  \Gamma_{3}^{-}  \}, $ we have
\begin{equation}\label{eq:largevaluebound-Nsecond} 
\sum_{\chi \in \Gamma} |\cN|^{2}
\ll N q^\varepsilon.
\end{equation}
For $ \Gamma \in \{ \Gamma_{3}^{+},  \Gamma_{3}^{-},  \Gamma_{4}^{+},  \Gamma_{4}^{-}  \} ,$ we have
\begin{equation}\label{eq:largevaluebound-Msixth} 
\sum_{\chi \in \Gamma} |\cM|^{6}
\ll Mq^{1+\varepsilon}.
\end{equation}
For $ \Gamma \in \{ \Gamma_{2}^{+},  \Gamma_{2}^{-},  \Gamma_{4}^{+},  \Gamma_{4}^{-}  \}, $ we have
\begin{equation}\label{eq:largevaluebound-Nsixth} 
\sum_{\chi \in \Gamma} |\cN|^{6}
\ll Nq^{1+\varepsilon}.
\end{equation}
\end{lemma}

\proof
For $ \Gamma \in \{ \Gamma_{1}^{+},  \Gamma_{1}^{-},  \Gamma_{2}^{+},  \Gamma_{2}^{-}  \} $, 
Lemma \ref{lm:largevalue-Huxley} yields
\begin{align*}
\sum_{\chi \in \Gamma} |\cM|^{2}
&\ll \log x \mathop{\sup}\limits_{q^{\frac{1}{4}}\leq A\ll M^{\frac{1}{2}}}
A^{2}(MA^{-2}+q^{1+\varepsilon}MA^{-6})
\ll Mq^\varepsilon.
\end{align*}
The proofs of other inequalities are similar.
\endproof

Following similar arguments, we also have following estimates.
\begin{lemma}\label{lm:largevaluebound-L} 
For $1\leqslant k\leqslant 4,$ we have
\begin{equation}\label{eq:largevaluebound-Lfourth} 
\begin{aligned}
\sum_{\chi \in \Gamma_{k}^{+}} |\cL|^{4}
\ll L^{2}q^\varepsilon
\end{aligned}
\end{equation}
\begin{equation}\label{eq:largevaluebound-L12th} 
\sum_{\chi \in \Gamma_{k}^{-}} |\cL|^{12}
\ll L^{2}q^{1+\varepsilon}.
\end{equation}
\end{lemma}

We are now ready to prove \eqref{eq:moments-expectedbound-minST}. Without loss of generality, we assume $q>x^{0.001}.$
Fix $\varepsilon>0$ to be a very small constant in terms of $\eta$. We also assume $ L\geq q^{\varepsilon}$, so that
\begin{equation}\label{eq:level-MN}
MN \leq  x^{1-\varepsilon}.
\end{equation}
From Lemma \ref{lm:incompletecharactersum-smoothmoduli}, it follows that
\begin{equation}\label{eq:shortcharactersum} 
\cL(s,\chi)\ll |s| L^{\frac{1}{2}-}.
\end{equation}
We also note that
\begin{align*}
\sum_{\substack{\chi\bmod{q}\\ \chi\neq \chi_{0}}}
|L(s,\chi)|^4\ll q^{1+\varepsilon}.
\end{align*}

{\bf (I) Treatment to $T_{1}^{\pm}$}:
From \eqref{eq:shortcharactersum} it follows that
\begin{align*}
T_{1}^{\pm}\ll |s|L^{\frac{1}{2}-}\sum_{\substack{\chi \in \Gamma_1^{\pm}}} |\cM||\cN|.
\end{align*}
By Cauchy--Schwarz inequality and \eqref{eq:largevaluebound-Msecond}, \eqref{eq:largevaluebound-Nsecond}, we have 
\begin{align*}
T_{1}^{\pm}\ll  |s|L^{\frac{1}{2}-}\Big(\sum_{\substack{\chi \in \Gamma_{1}^{\pm}}} |\cM|^{2}\Big)^{\frac{1}{2}}\Big(\sum_{\substack{\chi \in \Gamma_{1}^{\pm}}} |\cN|^{2}\Big)^{\frac{1}{2}}
\ll |s|x^{\frac{1}{2}-}
\end{align*}
as desired.

{\bf (II) Treatments to $S_2^{\pm}$ and $T_2^{\pm}$}:
As argued above, we have
\begin{align*}
T_2^{\pm}\ll |s|x^{\frac{1}{2}-}
\end{align*}
for  $ N \geq q$.
We henceforth assume that $N<q.$
Using H\"{o}lder's inequality, 
\begin{align*}
S_2^{\pm} 
&\ll \Big(\sum_{\substack{\chi \in \Gamma_{2}^{\pm}}} |\cM|^{2}\Big)^{\frac{1}{2}}
\Big(\sum_{\substack{\chi \in \Gamma_{2}^{\pm}}} |\cN|^{2}\Big)^{\frac{1}{8}} 
\Big(\sum_{\substack{\chi \in \Gamma_{2}^{\pm}}}|\cN|^{6}\Big)^{\frac{1}{8}}
\Big(\sum_{\chi \in \Gamma_{2}^{\pm}} |L|^{4}\Big)^{\frac{1}{4}}.
\end{align*}
Hence Lemma \ref{lm:largesievebound-Dirichletpolynomial} and \eqref{eq:largevaluebound-Nsixth} yield
\begin{align*}
S_{2}^{\pm} 
\ll (\log x)^{2} M^{\frac{1}{2}} q^{\frac{1}{8}}
(Nq)^{\frac{1}{8}}q^{\frac{1}{4}+\varepsilon}
\ll q^{\frac{1}{2}+\varepsilon}M^{\frac{1}{2}}
N^{\frac{1}{8}}.
\end{align*}

Using \eqref{eq:largevaluebound-Lfourth} and \eqref{eq:largevaluebound-L12th}, we obtain
\begin{align*}
T_{2}^{+}
&\ll \Big(\sum_{\substack{\chi \in \Gamma_{2}^{+}}} |\cM|^{2}\Big)^{\frac{1}{2}}\Big(\sum_{\substack{\chi \in \Gamma_{2}^{+}}} |\cN|^{2}\Big)^{\frac{1}{8}} 
\Big(\sum_{\substack{\chi \in \Gamma_{2}^{+}}}|\cN|^{6}\Big)
^{\frac{1}{8}}
\Big(\sum_{\chi \in \Gamma_{2}^{+}} |\cL|^{4}\Big)^{\frac{1}{4}}\\
&\ll (\log x)^{2} M^{\frac{1}{2}} q^{\frac{1}{8}}
(Nq)^{\frac{1}{8}}L^{\frac{1}{2}}
\\
&\ll q^{\frac{1}{4}+\varepsilon}M^{\frac{1}{2}}
N^{\frac{1}{8}}L^{\frac{1}{2}}
\end{align*}
and
\begin{align*}
T_{2}^{-}
&\ll \Big(\sum_{\substack{\chi \in \Gamma_{2}^{-}}} |\cM|^{2}\Big)^{\frac{1}{2}}\Big(\sum_{\substack {\chi \in \Gamma_{2}^{-}}} |\cN|^{2}\Big)^{\frac{3}{8}}    
\Big(\sum_{\substack{\chi \in \Gamma_{2}^{-}}}|\cN|^{6}\Big)
^{\frac{1}{24}}
\Big(\sum_{\chi \in \Gamma_{2}^{-}} |\cL|^{12}\Big)^{\frac{1}{12}}\\
&\ll (\log x)^{2} M^{\frac{1}{2}} q^{\frac{3}{8}}
(Nq)^{\frac{1}{24}}(qL^2)^{\frac{1}{12}}
\\
&\ll q^{\frac{1}{2}+\varepsilon}M^{\frac{1}{2}}
N^{\frac{1}{24}}L^{\frac{1}{6}}.
\end{align*}

Collecting the above bounds, we have
\begin{align*}
\min\{S_{2}^{+},~T_{2}^{+}\}
\ll q^{\frac{1}{2}+\varepsilon}M^{\frac{1}{2}}
N^{\frac{1}{8}}\min\{1,q^{-\frac{1}{4}}L^{\frac{1}{2}}\}
\ll x^{\frac{1}{8}}q^{\frac{7}{16}}
M^{\frac{3}{8}}q^\varepsilon
\end{align*}
and
\begin{align*}
\min\{S_{2}^{-},~T_{2}^{-}\}
\ll q^{\frac{1}{2}+\varepsilon}M^{\frac{1}{2}}
N^{\frac{1}{8}}\min\{1,N^{-\frac{1}{12}}L^{\frac{1}{6}}\}
\ll x^{\frac{1}{12}}q^{\frac{1}{2}}
M^{\frac{5}{12}}q^\varepsilon.
\end{align*}
Hence \eqref{eq:moments-expectedbound-minST} with $*\in\{+,-\}$ and $k=2$ holds as long as
\begin{align}\label{eq:level-M}
M\leq x^{1-\varepsilon}q^{-\frac{6}{5}},\ \ q\leq x^{\frac{5}{11}-\varepsilon}. 
\end{align}
In fact, the restriction $q\leq x^{\frac{5}{11}-\varepsilon}$ is additionally imposed to guarantee that $M$ is allowed to be bigger than $q$.

{\bf (III) Treatments to $S_3^{\pm}$ and $T_3^{\pm}$}:
The treatments are similar to (II) by interchanging $M$ and $ N$.
Hence \eqref{eq:moments-expectedbound-minST} with $*\in\{+,-\}$ and $k=3$ holds as long as
\begin{align}\label{eq:level-N}
N\leq x^{1-\varepsilon}q^{-\frac{6}{5}},\ \ q\leq x^{\frac{5}{11}-\varepsilon}. 
\end{align}

{\bf (IV) Treatments to $S_4^{\pm}$ and $T_4^{\pm}$}:
Firstly, we may follow the above arguments to derive that
\begin{align*}
T_4^{\pm}\ll |s|x^{\frac{1}{2}-}
\end{align*}
for $ M,N \geq q$. The case $N<q<M$ can be treated along with the arguments in (II), and by interchanging $M$ and $N$ one can also deal with the case $M<q<N.$
Hence we also need the restrictions \eqref{eq:level-M} and \eqref{eq:level-N}.

We henceforth assume that $M,N<q$, in which case 
H\"{o}lder's inequality gives
\begin{align*}
S_{4}^{\pm}
&\ll \Big(\sum_{\substack{\chi \in \Gamma_{4}^{\pm}}} |\cM|^{2}\Big)^{\frac{5}{16}}\Big(\sum_{\substack{\chi \in \Gamma_{4}^{\pm}}} |\cN|^{2}\Big)^{\frac{5}{16}} 
\Big(\sum_{\substack{\chi \in \Gamma_{4}^{\pm}}}|\cM|^{6}\Big)^{\frac{1}{16}}\Big(\sum_{\substack{\chi \in \Gamma_{4}^{\pm}}}|\cN|^{6}\Big)
^{\frac{1}{16}}
\Big(\sum_{\chi \in \Gamma_{4}^{+}} |L|^{4}\Big)^{\frac{1}{4}}.
\end{align*}
From Lemma \ref{lm:largesievebound-Dirichletpolynomial}, \eqref{eq:largevaluebound-Msixth} and \eqref{eq:largevaluebound-Nsixth}, it follows that
\begin{align*}
S_{4}^{\pm}
\ll (\log x)^{2} q^{\frac{5}{16}}q^{\frac{5}{16}} 
(Mq)^{\frac{1}{16}}(Nq)^{\frac{1}{16}}q^{\frac{1}{4}+\varepsilon}
\ll q^{1+\varepsilon}(MN)^{\frac{1}{16}}.
\end{align*}
Using \eqref{eq:largevaluebound-Lfourth} and \eqref{eq:largevaluebound-L12th}, we obtain
\begin{align*}
T_{4}^{+}
&\ll  
\Big(\sum_{\substack{\chi \in \Gamma_{4}^{+}}} |\cM|^{2}\Big)^{\frac{5}{16}}\Big(\sum_{\substack{\chi \in \Gamma_{4}^{+}}} |\cN|^{2}\Big)^{\frac{5}{16}} 
\Big(\sum_{\substack{\chi \in \Gamma_{4}^{+}}}|\cM|^{6}\Big)^{\frac{1}{16}}\Big(\sum_{\substack{\chi \in \Gamma_{4}^{+}}}|\cN|^{6}\Big)
^{\frac{1}{16}}
\Big(\sum_{\substack{\chi \in \Gamma_{4}^{+}}}|\cL|^{4}\Big)
^{\frac{1}{4}}\\
&\ll (\log x)^{2}q^{\frac{5}{16}}q^{\frac{5}{16}} 
(Mq)^{\frac{1}{6}}
(Nq)^{\frac{1}{6}} L^{\frac{1}{2}}\\
&\ll q^{\frac{3}{4}+\varepsilon}(MN)^{\frac{1}{16}}L^{\frac{1}{2}}
\end{align*}
and
\begin{align*}
T_{4}^{-}
&\ll  
\Big(\sum_{\substack{\chi \in \Gamma_{4}^{-}}} |\cM|^{2}\Big)^{\frac{7}{16}}
\Big(\sum_{\substack{\chi \in \Gamma_{4}^{-}}} |\cN|^{2}\Big)^{\frac{7}{16}} 
\Big(\sum_{\substack{\chi \in \Gamma_{4}^{-}}}|\cM|^{6}\Big)^{\frac{1}{48}}\Big(\sum_{\substack{\chi \in \Gamma_{4}^{-}}}|\cN|^{6}\Big)
^{\frac{1}{48}}
\Big(\sum_{\substack{\chi \in \Gamma_{4}^{-}}}|\cL|^{12}\Big)
^{\frac{1}{12}}\\
&\ll (\log x)^{2}q^{\frac{7}{16}}q^{\frac{7}{16}} 
(Mq)^{\frac{1}{48}}
(Nq)^{\frac{1}{48}} (qL^2)^{\frac{1}{12}}\\
&\ll q^{1+\varepsilon}(MN)^{\frac{1}{48}}L^{\frac{1}{6}}.
\end{align*}

Collecting the above bounds, we have
\begin{align*}
\min\{S_{4}^{+},~T_{4}^{+}\}
\ll q^{1+\varepsilon}(MN)^{\frac{1}{16}}\min\{1,q^{-\frac{1}{4}}L^{\frac{1}{2}}\}\ll x^{\frac{1}{16}}q^{\frac{31}{32}+\varepsilon}
\end{align*}
and
\begin{align*}
\min\{S_{4}^{-},~T_{4}^{-}\}
\ll q^{1+\varepsilon}(MN)^{\frac{1}{16}}\min\{1,L^{\frac{1}{6}}(MN)^{-\frac{1}{24}}\}\ll x^{\frac{1}{20}}q^{1+\varepsilon}.
\end{align*}
Hence \eqref{eq:moments-expectedbound-minST} with $*\in\{+,-\}$ and $k=4$ holds as long as
$q\leq x^{\frac{9}{20}-\varepsilon}$.

Combining the above treatments in four parts, as well as the restrictions \eqref{eq:level-MN}, \eqref{eq:level-M} and \eqref{eq:level-N}, we choose
\begin{align*}
M=N= x^{1-\varepsilon}q^{-\frac{6}{5}}
\end{align*}
if $q\leq x^{\frac{9}{20}-\varepsilon}$, in which case we may take 
\begin{align*}
D=\min\{x^{1-\varepsilon},~x^{2-\varepsilon}q^{-\frac{12}{5}}\}
\end{align*}
in Proposition \ref{prop:sieveinequality}. This proves Theorem \ref{thm:BT-smoothmoduli<1/2} readily.

\smallskip

\section{Proof of Theorems \ref{thm:BT-generalmoduli>1/2} and \ref{thm:BT-smoothmoduli2/3generalexponentpair}}
\label{sec:BT-generalmoduli>1/2}
Keep in mind that $q\geqslant x^{\frac{1}{2}}$ throughout this section, and we appeal to Proposition \ref{prop:sieveinequality} with the level $D$ as large as possible.
\subsection{Approaching by additive characters}
Recall that
\begin{align*}
A(d)=\sum_{n\equiv a\overline{d}\bmod q} \varPhi(dn)
\end{align*}
and $\widehat{\varPhi}(0)=x.$
It follows from Poisson summation that
\begin{align*}
A(d)
=\frac{1}{dq}\sum_{h\in\bZ} \widehat{\varPhi}\Big(\frac{h}{dq}\Big)\ue\Big(\frac{ah\overline{d} }{q}\Big).
\end{align*}
From the rapid decay of $\widehat{\varPhi}$, we may write
\begin{align*}
r(d)
&=A(d)-\frac{x}{dq}
=\frac{1}{dq} \sum_{0<|h|\leqslant H}\widehat{\varPhi}\Big(\frac{h}{dq}\Big)
\mathrm{e}\Big(\frac{ah\overline{d}}{q}\Big)+O(x^{-1})
\end{align*}
for $H:=MNqX^{-1+\varepsilon}$.

Take $\lambda$ to be well-factorable of degree $2$; that is, 
\begin{align*}
\lambda=\balpha*\bbeta,
\end{align*}
where $\balpha$ and $\bbeta$  are of level $M$ and $N$, respectively, and 
$$
\quad 
MN=D 
\quad\text{and}\quad
\|\balpha\|_\infty,\|\bbeta\|_\infty \leq 1.$$
It suffices to show that
\begin{align*} 
\sum_{m\leqslant M}\sum_{n\leqslant N}\alpha_m\beta_n r(mn)\ll \frac{x}{q(\log x)^{2024}}.
\end{align*}
By dyadic device, it suffices to prove that
\begin{align}\label{eq:R(M,N)-expectedbound} 
\cR(M,N):=\sum_{m\sim M}\sum_{n\sim N} \frac{\alpha_m\beta_n}{mn} \sum_{0<|h|\leqslant H}\widehat{\varPhi}\Big(\frac{h}{mnq}\Big)
\mathrm{e}\Big(\frac{ah\overline{mn}}{q}\Big)\ll \frac{x}{(\log x)^{2026}},
\end{align}
where $MN=D$ and $H=x^{-1}MNq^{1+\varepsilon}$.

\subsection{Transforming exponential sums}
We transform the above multiple exponential sum in $R(M,N)$ into averages of complete Kloosterman sums. To do so,
we first write
\begin{align*} 
\cR(M,N)
&=\sum_{m\sim M}\sum_{n\sim N} \frac{\alpha_m\beta_n}{mn} \sum_{0<|h|\leqslant H}\int_\bR \varPhi(t)\ue\Big(\frac{-ht}{mnq}\Big)
\mathrm{e}\Big(\frac{ah\overline{mn}}{q}\Big)\ud t\\
&=\sum_{m\sim M}\sum_{n\sim N} \frac{\alpha_m\beta_n}{n} \sum_{0<|h|\leqslant H}\int_\bR \varPhi(mt)\ue\Big(\frac{-ht}{nq}\Big)
\mathrm{e}\Big(\frac{ah\overline{mn}}{q}\Big)\ud t\\
&\ll \frac{1}{N}\int_\bR \sum_{m\sim M}\varPhi(mt)\Big|\sum_{0<|h|\leqslant H}\sum_{n\sim N} \beta_{h,n,t}\ue\Big(\frac{ah\overline{mn}}{q}\Big)\Big|\ud t,
\end{align*}
where $|\beta_{h,n,t}|\leqslant 1.$
Therefore,
\begin{align*} 
\cR(M,N)
&\ll \frac{x}{MN}\sum_{m\in\bZ}W\Big(\frac{m}{M}\Big)\Big|\sum_{0<|h|\leqslant H}\sum_{n\sim N} \beta_{h,n,t}\ue\Big(\frac{ah\overline{mn}}{q}\Big)\Big|
\end{align*}
for some $t\in\bR$ and a smooth function $W$ which dominates the indicator function of $[1,2].$

By Cauchy--Schwarz inequality, we obtain
\begin{align*} 
\cR(M,N)^2
&\ll \frac{x^2}{MN^2}\sum_{m\in\bZ}W\Big(\frac{m}{M}\Big)\Big|\sum_{0<|h|\leqslant H}\sum_{n\sim N} \beta_{h,n,t}\ue\Big(\frac{ah\overline{mn}}{q}\Big)\Big|^2.
\end{align*}
Squaring out and switching summations, it follows that
\begin{align*} 
\cR(M,N)^2
&\ll \frac{x^2}{MN^2}\mathop{\sum\sum\sum\sum}_{\substack{0<|h_1|,|h_2|\leqslant H,~n_1,n_2\sim N\\ (n_1n_2,q)=1}}\beta_{h_1,n_1,t}\overline{\beta_{h_2,n_2,t}}\sum_{\substack{m\in\bZ\\ (m,q)=1}}\\
&\ \ \ \ \times W\Big(\frac{m}{M}\Big)\ue\Big(\frac{a(h_1n_2-h_2n_1)\overline{n_1n_2m}}{q}\Big).
\end{align*}
The diagonal terms with $h_1n_2=h_2n_1$ contribute at most
\begin{align*}
&\ll \frac{x^2}{N^2}\mathop{\sum\sum\sum\sum}_{\substack{0<|h_1|,|h_2|\leqslant H,~n_1,n_2\sim N\\ h_1n_2=h_2n_1}}|\beta_{h_1,n_1,t}\beta_{h_2,n_2,t}|\ll x^2HN^{-1}q^\varepsilon.
\end{align*}
To ensure the desired bound \eqref{eq:R(M,N)-expectedbound}, we require that 
\begin{align}\label{eq:R(M,N)-Msize} 
M \leqslant q^{-1} x^{1-\varepsilon}.
\end{align}
We now have
\begin{align}\label{eq:R(M,N)-Sigma}
\cR(M,N)^2
&\ll x^2HN^{-1}q^\varepsilon+\frac{x^2}{MN^2}\varSigma
\end{align}
with
\begin{align}\label{eq:Sigma}
\varSigma
&=\mathop{\sum\sum\sum\sum}_{\substack{0<|h_1|,|h_2|\leqslant H,~n_1,n_2\sim N\\ h_1n_2\neq h_2n_1\\ (n_1n_2,q)=1}}\beta_{h_1,n_1,t}\overline{\beta_{h_2,n_2,t}}\sum_{\substack{m\in\bZ\\ (m,q)=1}}W\Big(\frac{m}{M}\Big)\ue\Big(\frac{a(h_1n_2-h_2n_1)\overline{n_1n_2m}}{q}\Big).
\end{align}

From Poisson summation, it follows that
\begin{align*}
\varSigma
&=\frac{M}{q}\mathop{\sum\sum\sum\sum}_{\substack{0<|h_1|,|h_2|\leqslant H,~n_1,n_2\sim N\\ h_1n_2\neq h_2n_1\\ (n_1n_2,q)=1}}\beta_{h_1,n_1,t}\overline{\beta_{h_2,n_2,t}}\sum_{k\in\bZ}\widehat{W}\Big(\frac{kM}{q}\Big)S(k,a(h_1n_2-h_2n_1)\overline{n_1n_2};q).
\end{align*} 
Thanks to the rapid decay of $\widehat{W},$ we may truncate the $k$-sum at $K:=M^{-1}q^{1+\varepsilon}$ with a negligible error term. 
Note that the zero-th frequency contributes to $\varSigma$ at most
\begin{align*}
&\ll \frac{M}{q}\mathop{\sum\sum\sum\sum}_{\substack{0<|h_1|,|h_2|\leqslant H,~n_1,n_2\leqslant N\\ h_1n_2\neq h_2n_1}}(h_1n_2-h_2n_1,q)\ll q^{-1+\varepsilon}M(HN)^2.
\end{align*}
Hence
\begin{align}\label{eq:Sigma-Sigma*}
\varSigma
&=\varSigma^*+O(q^{-1+\varepsilon}M(HN)^2)
\end{align}
with
\begin{align*}
\varSigma^*
&=\frac{M}{q}\mathop{\sum\sum\sum\sum}_{\substack{0<|h_1|,|h_2|\leqslant H,~n_1,n_2\sim N\\ h_1n_2\neq h_2n_1\\ (n_1n_2,q)=1}}\beta_{h_1,n_1,t}\overline{\beta_{h_2,n_2,t}}\sum_{1\leqslant |k|\leqslant K}\widehat{W}\Big(\frac{kM}{q}\Big)\\
&\ \ \ \ \times S(k,a(h_1n_2-h_2n_1)\overline{n_1n_2};q).
\end{align*} 
By employing Weil's bound for individual Kloosterman sums, one may produce an upper bound for $\varSigma^*$, and thus for $\cR(M,N).$ This was exactly done by Iwaniec \cite{Iw82}. We would like to estimate $\varSigma^*$ by controlling the sign changes of Kloosterman sums effectively. More precisely, we will group certain variables and create bilinear or trilinear forms to capture cancellations among Kloosterman sums.

There exist at least two means to group variables:
\begin{itemize}
\item For $u\bmod q,$ put
\begin{align*}
\rho(u)
&=\mathop{\sum\sum\sum\sum}_{\substack{0<|h_1|,|h_2|\leqslant H,~n_1,n_2\sim N\\ h_1n_2\neq h_2n_1,~(n_1n_2,q)=1\\ (h_1n_2-h_2n_1)\overline{n_1n_2}\equiv u\bmod q}}\beta_{h_1,n_1,t}\overline{\beta_{h_2,n_2,t}}.
\end{align*}

\item For $b\leqslant B$ and $c\leqslant C$ with $B=4HN$ and $C=4N^2,$ put
\begin{align*}
\rho(b,c)
&=\mathop{\sum\sum\sum\sum}_{\substack{0<|h_1|,|h_2|\leqslant H,~n_1,n_2\sim N\\ h_1n_2\neq h_2n_1,~(n_1n_2,q)=1\\ h_1n_2-h_2n_1=b,~n_1n_2=c}}\beta_{h_1,n_1,t}\overline{\beta_{h_2,n_2,t}}.
\end{align*}
\end{itemize}
Hence
\begin{align*}
\varSigma^*
&=\frac{M}{q}\sum_{u\bmod q}\rho(u)\sum_{1\leqslant |k|\leqslant K}\delta_k S(k,au;q)\\
&=\frac{M}{q}\sum_{b\leqslant B}\sum_{c\leqslant C}\rho(b,c)\sum_{1\leqslant |k|\leqslant K}\delta_k S(k,ab\overline{c};q),
\end{align*} 
where $\delta_k=\widehat{W}(kM/q).$
In the remaining of this section, we specialize $q$ to be a large prime with $M>2q^\varepsilon$, so that $K<q/2$ and
\begin{align}
\varSigma^*
&=\frac{M}{\sqrt{q}}\sum_{u\bmod q}\rho(u)\sum_{1\leqslant |k|\leqslant K}\delta_k \kl(aku,q)
\label{eq:Sigma*-bilinear}\\
&=\frac{M}{\sqrt{q}}\sum_{b\leqslant B}\sum_{c\leqslant C}\rho(b,c)\sum_{1\leqslant |k|\leqslant K}\delta_k \kl(akb\overline{c},q)\label{eq:Sigma*-trilinear},
\end{align} 
where $\kl(\cdot,q)=S(\cdot,1;q)/\sqrt{q}.$

\subsection{Bounding $\varSigma^*$ as bilinear forms}
For each $\nu\in\bZ^+,$ we infer from \eqref{eq:Sigma*-bilinear} and H\"older's inequality that
\begin{align*}
\varSigma^*
&\ll\frac{M}{\sqrt{q}}\cB_1^{1-\frac{1}{\nu}}(\cB_2\cB_3)^{\frac{1}{2\nu}},
\end{align*}
where
\begin{align*}
\cB_j=\sum_{u\bmod q}\rho(u)^j,\ \ j=1,2,
\end{align*}
and
\begin{align*}
\cB_3
&=\sum_{u\bmod q}\Big|\sum_{1\leqslant |k|\leqslant K}\delta_k \kl(ku,q)\Big|^{2\nu}.
\end{align*}

Trivially, we have
\begin{align}\label{eq:B1-upperbound}
\cB_1&\ll (HN)^2.
\end{align}
From Lemma \ref{lm:Klsum-2nu-thmoment}, we also find
\begin{align}\label{eq:B3-upperbound}
\cB_3
&\ll K^{\nu}q+K^{2\nu}q^{\frac{1}{2}}.
\end{align}

To bound $\cB_2$, we need the following lemma.
\begin{lemma}\label{lm:B2-upperbound}
With the above notation, we have
\begin{align*}
\cB_2\ll (1+HN^3/q)(1+H/N)(HN)^2q^\varepsilon
\end{align*}
for any $\varepsilon>0.$
\end{lemma}

\proof
We first write
\begin{align*}
\cB_2
&\ll \mathop{\sum\sum\sum\sum\sum\sum\sum\sum}_{\substack{\left(h_1 n_2-h_2 n_1\right)n_1'n_2'\equiv (h_1'n_2'-h_2' n_1')n_1n_2\bmod{q} \\ n_1,n_1',n_2,n_2'\sim N,~1\leqslant |h_1|,|h_1'|,|h_2|,|h_2'|\leqslant H\\(n_1n_1'n_2n_2',q)=1\\
(h_1 n_2-h_2 n_1)(h_1'n_2'-h_2' n_1')\neq0}}1.
\end{align*}
Note that the above congruence restriction $(h_1 n_2-h_2 n_1)n_1'n_2'\equiv (h_1'n_2'-h_2' n_1')n_1n_2
\bmod{q}$ implies
\begin{align*}
(h_1 n_2-h_2 n_1)n_1'n_2'= (h_1'n_2'-h_2' n_1')n_1n_2+qt,\ \ t\in\bZ\cap[-T,T]
\end{align*}
with $T:=4HN^3/q.$
Denote by $\cB_{2,0}$ and $\cB_{2,1}$ the contributions from $t=0$ and $t\neq0$, respectively.

We now consider $\cB_{2,1}$, which vanishes unless $T\geqslant1$. Picking out the g.c.d. of $n_1,n_2$ and that of $n_1',n_2'$, each of which divides $t\neq0$, we obtain
\begin{align*}
\cB_{2,1}
&\ll \log^2N\sup_{D,D'\leqslant 2N}\cB_{2,1}(D,D')
\end{align*}
with
\begin{align*}
\cB_{2,1}(D,D')
&=\mathop{\sum\sum\sum\sum\sum\sum\sum\sum\sum\sum\sum}_{\substack{d\sim D,~d'\sim D',~n_1,n_2\sim N/d,~n_1',n_2'\sim N/d',~1\leqslant |h_1|,|h_1'|,|h_2|,|h_2'|\leqslant H,~1\leqslant |t|\leqslant T/DD'\\(n_1,n_2)=(n_1',n_2')=1\\
(h_1 n_2-h_2 n_1)(h_1'n_2'-h_2' n_1')\neq0}}1,
\end{align*}
where the summations are also restricted to
\begin{align*}
(h_1 n_2-h_2 n_1)d'n_1'n_2'= (h_1'n_2'-h_2' n_1')dn_1n_2+qt.
\end{align*}

Given an integer $\ell$ with $0<|\ell|\leqslant2HN/D,$  consider the equations
\begin{align}\label{eq:h1h2n1n2''''t-equation1}
h_1 n_2-h_2 n_1=\ell,
\end{align}
and
\begin{align}\label{eq:h1h2n1n2''''t-equation2}
\ell d'n_1'n_2'=(h_1'n_2'-h_2' n_1')dn_1n_2+qt.
\end{align}
If fixing $n_1,n_2$ with $(n_1,n_2)=1$, the number of pairs $(h_1,h_2)$ with $1\leqslant |h_1|,|h_2|\leqslant H$ satisfying \eqref{eq:h1h2n1n2''''t-equation1} is $\ll 1+H/n_1\ll 1+DH/N.$ 
Regarding the equation \eqref{eq:h1h2n1n2''''t-equation2}, we turn to consider the congruence equation
\begin{align}\label{eq:h1h2n1n2''''t-equation3}
(h_1'n_2'-h_2' n_1')dn_1n_2\equiv -qt\bmod{d'n_1'n_2'}.
\end{align}
Given $h_1',h_2',n_1',n_2',t$, the number of tuples $(d,n_1,n_2)$ satisfying \eqref{eq:h1h2n1n2''''t-equation3}
with $d\sim D, n_1,n_2\sim N/d$ is 
\begin{align*}
\ll \Big(\frac{N^2}{D}\frac{(h_1'n_2'-h_2' n_1',d'n_1'n_2')}{d'n_1'n_2'}+1\Big)N^{\varepsilon}
\ll \Big(\frac{D'}{D}(h_1'n_2'-h_2' n_1',d'n_1'n_2')+1\Big)N^{\varepsilon},\end{align*}
where the factor $N^\varepsilon$ comes from the bound $\tau_3(m)\ll m^\varepsilon.$

Collecting all above arguments, we find
\begin{align*}
\cB_{2,1}(D,D')
&\ll q^\varepsilon (1+DH/N)\frac{T}{DD'}\mathop{\sum\sum\sum\sum}_{\substack{d'\sim D',~n_1',n_2'\sim N/d',~1\leqslant |h_1'|,|h_2'|\leqslant H\\
h_1'n_2'-h_2' n_1'\neq0, ~ (n_1',n_2')=1}}\Big(\frac{D'}{D}(h_1'n_2'-h_2' n_1',d'n_1'n_2')+1\Big)\\
&\ll q^\varepsilon (1+DH/N)\frac{T}{DD'}(HN)^2,
\end{align*}
so that
\begin{align*}
\cB_{2,1}
&\ll (1+H/N)(HN)^2Tq^\varepsilon.
\end{align*}

Following similar arguments, we obtain
\begin{align*}
\cB_{2,0}
&\ll (1+H/N)(HN)^2q^\varepsilon.
\end{align*}
Therefore, we conclude
\begin{align*}
\cB_2
&\ll (1+H/N)(1+T)(HN)^2q^\varepsilon
\end{align*}
as desired.
\endproof

We now come back to the estimate for $\varSigma^*.$
From \eqref{eq:B1-upperbound}, \eqref{eq:B3-upperbound} and Lemma \ref{lm:B2-upperbound}, it follows that
\begin{align*}
\varSigma^*
&\ll\frac{M}{\sqrt{q}}(HN)^{2-\frac{1}{\nu}}(1+H/N)^{\frac{1}{2\nu}}(1+HN^3/q )^{\frac{1}{2\nu}}(K^{2\nu}q^{\frac{1}{2}}+K^{\nu}q)^{\frac{1}{2\nu}}q^\varepsilon.
\end{align*}
Recall that $H=x^{-1}MNq^{1+\varepsilon}$, $K=M^{-1}q^{1+\varepsilon}$ and $M \leqslant q^{-1} x^{1-\varepsilon}$. It is useful to make a restriction $N\ll q^{\frac{1}{4}}$, in which case we find $H\ll N$ and $HN^3\ll q^{1+\varepsilon}.$
Hence
\begin{align*}
\varSigma^*
&\ll (M^{2-\frac{1}{\nu}}q^{\frac{5}{2}-\frac{3}{4\nu}}+ M^{\frac{5}{2}-\frac{1}{\nu}}q^{2-\frac{1}{2\nu}})x^{-2+\frac{1}{\nu}}N^{4-\frac{2}{\nu}}q^\varepsilon.
\end{align*}

Taking $\nu=3$, we arrive at
\begin{align*}
\varSigma^*
&\ll (M^{\frac{5}{3}}q^{\frac{9}{4}}+ M^{\frac{13}{6}}q^{\frac{11}{6}})x^{-\frac{5}{3}}N^{\frac{10}{3}}q^\varepsilon.
\end{align*}
 Together with \eqref{eq:R(M,N)-Sigma} and \eqref{eq:Sigma-Sigma*}, this yields
\begin{align*}
\cR(M,N)^2
&\ll (xqM+x^{\frac{1}{3}}q^{\frac{9}{4}}M^{\frac{2}{3}}N^{\frac{4}{3}}+x^{\frac{1}{3}}q^{\frac{11}{6}}M^{\frac{7}{6}}N^{\frac{4}{3}}+qM^2N^2)q^{\varepsilon}.
\end{align*}
Hence \eqref{eq:R(M,N)-expectedbound} holds for
\begin{align}\label{eq:bilinearforms-choose-MN1}
M\leq x^{1-\varepsilon}q^{-1},\ \ N\leq \min\{x^{\frac{3}{4}}q^{-\frac{19}{16}}, ~ x^{\frac{3}{8}}q^{-\frac{1}{2}},~q^{\frac{1}{4}}\}x^{-\varepsilon}
\end{align}
since  $q\geqslant x^{\frac{1}{2}}$.

Taking $\nu=4$, we arrive at
\begin{align*}
\varSigma^*
&\ll (M^{\frac{7}{4}}q^{\frac{37}{16}}+ M^{\frac{9}{4}}q^{\frac{15}{8}})x^{-\frac{7}{4}}N^{\frac{7}{2}}q^\varepsilon.
\end{align*}
Together with \eqref{eq:R(M,N)-Sigma} and \eqref{eq:Sigma-Sigma*}, this yields
\begin{align*}
\cR(M,N)^2
&\ll (xqM+x^{\frac{1}{4}}q^{\frac{37}{16}}M^{\frac{3}{4}}N^{\frac{3}{2}}+x^{\frac{1}{4}}q^{\frac{15}{8}}M^{\frac{5}{4}}N^{\frac{3}{2}}+qM^2N^2)q^{\varepsilon}.
\end{align*}
Hence \eqref{eq:R(M,N)-expectedbound} holds for
\begin{align}\label{eq:bilinearforms-choose-MN2}
M\leq x^{1-\varepsilon}q^{-1},\ \ N\leq \min\{x^{\frac{2}{3}}q^{-\frac{25}{24}}, ~ x^{\frac{1}{3}}q^{-\frac{5}{12}},~q^{\frac{1}{4}}\}x^{-\varepsilon}
\end{align}
since  $q\geqslant x^{\frac{1}{2}}$.

\subsection{Bounding $\varSigma^*$ as trilinear forms}
We infer from \eqref{eq:Sigma*-trilinear} and Cauchy--Schwarz inequality that
\begin{align}\label{eq:Sigma*-T1T2}
\varSigma^*
&\ll \frac{M}{\sqrt{q}}(\cT_1\cT_2)^{\frac{1}{2}},
\end{align}
where
\begin{align*}
\cT_1&=\sum_{b\leqslant B}\sum_{c\leqslant C}\rho(b,c)^2,\\
\cT_2&=\mathop{\sum\sum}_{\substack{b,c\in\bZ\\(c,q)=1}}W\Big(\frac{b}{B}\Big)W\Big(\frac{c}{C}\Big)\Big|\sum_{1\leqslant |k|\leqslant K}\delta_k \kl(akb\overline{c},q)\Big|^2.
\end{align*}

Squaring out and switching summations, 
\begin{align*}
\cT_2&\ll \mathop{\sum\sum}_{1\leqslant |k_1|,|k_2|\leqslant K}|\cT_2(k_1,k_2)|,
\end{align*}
where
\begin{align*}
\cT_2(k_1,k_2)&=\mathop{\sum\sum}_{\substack{b,c\in\bZ\\(c,q)=1}}W\Big(\frac{b}{B}\Big)W\Big(\frac{c}{C}\Big)\kl(ak_1b\overline{c},q)\kl(ak_2b\overline{c},q).
\end{align*}
A trivial bound reads $\cT_2(k_1,k_2)\ll BC$, which will be applied for diagonal terms. For off-diagonal terms, i.e., $k_1\neq k_2$, we would like to bound $\cT_2(k_1,k_2)$ by appealing to the work of Fouvry, Kowalski and Michel \cite{FKM14} on non-correlations of divisor functions and Frobenius trace functions; see Lemma \ref{lm:FKM}.

From Poisson summation, it follows that
\begin{align*}
\sum_{b\in\bZ}W\Big(\frac{b}{B}\Big)\kl(ak_1b\overline{c},q)\kl(ak_2b\overline{c},q)
&=\frac{B}{\sqrt{q}}\sum_{l\in\bZ}\widehat{W}\Big(\frac{Bl}{q}\Big)K(lc,q;k_1,k_2),
\end{align*}
where
\begin{align*}
K(z,q;k_1,k_2)=\frac{1}{\sqrt{q}}\sum_{u\bmod q}\kl(k_1u,q)\kl(k_2u,q)\ue\Big(\frac{-zu}{q}\Big).
\end{align*}
This gives
\begin{align*}
\cT_2(k_1,k_2)&=\frac{B}{\sqrt{q}}\mathop{\sum\sum}_{\substack{l,c\in\bZ\\(c,q)=1}}\widehat{W}\Big(\frac{Bl}{q}\Big)W\Big(\frac{c}{C}\Big)K(lc,q;k_1,k_2).
\end{align*}

Observe that
\begin{align*}
K(0,q;k_1,k_2)
&=\frac{1}{\sqrt{q}}\sum_{u\bmod q}\kl(k_1u,q)\kl(k_2u,q)=q^{-\frac{1}{2}}S(k_1-k_2;0;q).
\end{align*}
We employ the trivial bound for $l=0$, and appeal to Lemma \ref{lm:FKM} for $l\neq0$ with the aid of Lemma \ref{lm:Kloostermantensor}.
Hence the off-diagonal terms contribute at most
\begin{align*}
\cT_2(k_1,k_2)
&\ll \frac{BC}{q}(k_1-k_2, q)+\frac{B}{\sqrt{q}}\frac{qC}{B}\Big(\frac{B}{C}+1\Big)^{\frac{1}{2}}q^{-\frac{1}{8}+\varepsilon}\\
&\ll \frac{HN^3}{q}(k_1-k_2, q)+\Big(\frac{H}{N}+1\Big)^{\frac{1}{2}}N^2q^{\frac{3}{8}+\varepsilon}.
\end{align*}
Summing over $k_1,k_2$, we find
\begin{align}
\cT_2
&\ll HN^3K+HN^3K^2q^{-1+\varepsilon}+\Big(\frac{H}{N}+1\Big)^{\frac{1}{2}}(NK)^2q^{\frac{3}{8}+\varepsilon}\nonumber\\
&\ll \frac{N^4}{x}q^{2+\varepsilon}+\frac{N^2}{M^2}q^{\frac{19}{8}+\varepsilon}\label{eq:T2-upperbound}
\end{align}
since $M<x/q.$

Note that
\begin{align*}
\cT_1
&\ll\mathop{\sum\sum\sum\sum\sum\sum\sum\sum}_{\substack{h_1,h_2,h_1',h_2'\leqslant H,~n_1,n_2,n_1',n_2'\sim N\\ h_1n_2- h_2n_1=h_1'n_2'- h_2'n_1',~n_1n_2=n_1'n_2'}}1
\ll (HN)^2\Big(\frac{H}{N}+1\Big)q^\varepsilon
\ll \Big(\frac{qMN^2}{x}\Big)^2q^\varepsilon,
\end{align*}
from which and \eqref{eq:T2-upperbound}, \eqref{eq:Sigma*-T1T2} we derive that
\begin{align*}
\varSigma^*
&\ll \Big(\frac{q^{\frac{3}{2}}M^2N^4}{x^{\frac{3}{2}}}+\frac{q^{\frac{27}{16}}MN^3}{x}\Big)q^\varepsilon.
\end{align*}
Together with \eqref{eq:R(M,N)-Sigma} and \eqref{eq:Sigma-Sigma*}, this yields
\begin{align*}
\cR(M,N)^2
&\ll (xqM+x^{\frac{1}{2}}q^{\frac{3}{2}}MN^2+xq^{\frac{27}{16}}N+qM^2N^2)q^{\varepsilon}.
\end{align*}
Hence \eqref{eq:R(M,N)-expectedbound} holds for
\begin{align}\label{eq:trilinearforms-choose-MN}
M\leqslant x^{1-\varepsilon}q^{-1},\ \ N\leqslant\min\{xq^{-\frac{27}{16}},~(x/q)^{\frac{1}{4}}\}x^{-\varepsilon}
\end{align}
since  $q\geqslant x^{\frac{1}{2}}$.

In view of \eqref{eq:bilinearforms-choose-MN1}, \eqref{eq:bilinearforms-choose-MN2} and \eqref{eq:trilinearforms-choose-MN}, we find \eqref{eq:R(M,N)-expectedbound} is valid if taking $M=x^{1-\varepsilon}q^{-1}$ and
\begin{align*}
N=
\begin{cases}
\min\{xq^{-\frac{27}{16}},~(x/q)^{\frac{1}{4}}\}x^{-\varepsilon},\ \ &\varpi\in[1/2,~32/61[,\\
\min\{x^{\frac{2}{3}}q^{-\frac{25}{24}}, ~x^{\frac{1}{3}}q^{-\frac{5}{12}}\}x^{-\varepsilon},&\varpi\in[32/61,~7/13[,\\
\min\{x^{\frac{3}{4}}q^{-\frac{19}{16}}, ~x^{\frac{3}{8}}q^{-\frac{1}{2}}\}x^{-\varepsilon}, &\varpi\in[7/13,~4/7[.
\end{cases}
\end{align*}
Taking $D=MN$ in Proposition \ref{prop:sieveinequality}, we conclude the proof of Theorem \ref{thm:BT-generalmoduli>1/2}.

\begin{remark}
The estimates for diagonal terms become better when we take larger $\nu$ in H\"older's inequality. In particular, if we take $\nu =2,$ the diagonal terms will lead to a constraint $MN\leq x^{\frac{3}{2}}q^{-\frac{7}{4}},$ which coincides with Iwaniec’s result. On the other hand, estimates for off-diagonal terms become worse when we take larger $\nu.$
\end{remark}

\begin{remark}
Each of the above choices of $N$ fails to improve the work of Iwaniec \cite{Iw82} when $\varpi=1/2.$ This is due to the fact that we have explored the bilinear or trilinear structures in $\varSigma^*,$ which requires that $K>q^{4\varepsilon}.$ This is impossible by recalling that $K=M^{-1}q^{1+\varepsilon}$ and $M=x^{1-\varepsilon}q^{-1}.$ We will beat this barrier in the second paper of this series {\rm(\cite{XZ25}),} i.e., we are able to improve the work of Iwaniec  \cite{Iw82} when $\varpi=1/2.$
\end{remark}

\subsection{Proof of Theorem $\ref{thm:BT-smoothmoduli2/3generalexponentpair}$}
\label{subsec:BT-smoothmoduli2/3}

Recall that
\begin{align*}
\cR(M,N)^2
&\ll x^2HN^{-1}q^\varepsilon+\frac{x^2}{MN^2}\varSigma
\end{align*}
with
\begin{align*}
\varSigma
&=\mathop{\sum\sum\sum\sum}_{\substack{0<|h_1|,|h_2|\leqslant H,~n_1,n_2\sim N\\ h_1n_2\neq h_2n_1\\ (n_1n_2,q)=1}}\beta_{h_1,n_1,t}\overline{\beta_{h_2,n_2,t}}\sum_{\substack{m\in\bZ\\ (m,q)=1}}W\Big(\frac{m}{M}\Big)\ue\Big(\frac{a(h_1n_2-h_2n_1)\overline{n_1n_2m}}{q}\Big).
\end{align*}
as in \eqref{eq:R(M,N)-Sigma} and \eqref{eq:Sigma}.
Applying Lemma \ref{lm:incompleteKloostermansum-smoothmoduli} with
$h=a(h_1 n_2-h_2 n_1) \overline{n_1 n_2}$ and partial summation to the inner sum, 
we get
\begin{align*}
\varSigma
&\ll \mathop{\sum\sum\sum\sum}_{\substack{0<|h_1|,|h_2|\leqslant H,~n_1,n_2\sim N\\ h_1n_2\neq h_2n_1}}\Big(\frac{M}{q}(h_1n_2-h_2n_1,q)+q^{\kappa}M^{\lambda-\kappa+{O(\eta)}}(h_1n_2-h_2n_1,q)^{\nu}\Big).
\end{align*}
Note that
\begin{align*}
\mathop{\sum\sum\sum\sum}_{\substack{0<|h_1|,|h_2|\leqslant H,~n_1,n_2\sim N\\ h_1n_2\neq h_2n_1}}(h_1n_2-h_2n_1,q)&\ll (HN)^2q^\varepsilon,
\end{align*}
we thus obtain
\begin{align*}
\varSigma
&\ll (x^{-2}M^3 N^4 q+x^{-2}M^{2+\lambda-\kappa}N^4q^{\kappa+2})q^{O(\eta)},
\end{align*}
which yields
\begin{align*}
\cR(M,N)\leqslant
(x^{\frac{1}{2}}q^{\frac{1}{2}}M^{\frac{1}{2}}+q^{\frac{1}{2}}MN+ q^{1+\frac{\kappa}{2}}M^{\frac{1-\kappa+\lambda}{2}}N) q^{O(\eta)}.
\end{align*}

Hence \eqref{eq:R(M,N)-expectedbound} holds for
\begin{align*}
M=q^{-1}x^{1-\varepsilon},\ \ N=\min\{x^{\frac{1+\kappa-\lambda}{2}} q^{-\frac{1+2\kappa-\lambda}{2}},~q^{\frac{1}{2}}\}x^{-\varepsilon}
\end{align*}
with $\varepsilon\gg\eta.$
In particular, for $\varpi\geqslant \frac{1+\kappa-\lambda}{2+2\kappa-\lambda}$, the choice $N=x^{\frac{1+\kappa-\lambda}{2}-\varepsilon} q^{-\frac{1+2\kappa-\lambda}{2}}$ is admissible, so that we may take
\begin{align*}
D=x^{\frac{3+\kappa-\lambda}{2}-2\varepsilon} q^{-\frac{3+2\kappa-\lambda}{2}}
\end{align*}
in Proposition \ref{prop:sieveinequality}. Note that van Lint and Richert \cite{LR65} obtained a level $q^{-1}x^{1-\varepsilon}$, and the above choice of $D$ is better as long as
$\varpi \leq \frac{1+\kappa-\lambda}{1+2\kappa-\lambda}$.
The proof of Theorem \ref{thm:BT-smoothmoduli2/3generalexponentpair} is now finished.

\smallskip

\section{Proof of Theorems \ref{thm:countingcongruence-generalcoefficient} and \ref{thm:countingcongruence-specialcoefficient}}
\label{sec:countingcongruence-proof}

The proofs of Theorems \ref{thm:countingcongruence-generalcoefficient} and \ref{thm:countingcongruence-specialcoefficient} will follow the same lines, except that Lemma \ref{lm:trilinearform-BC} applies in the proof of Theorem \ref{thm:countingcongruence-generalcoefficient}, and we utilize convolution structures of the coefficients when proving Theorem \ref{thm:countingcongruence-specialcoefficient}, so that Kloostermania (Lemma \ref{lm:trilinearKloosterman-inverse}) applies.

\subsection{Preparations}
Denote by $\omega(n)$ the number of distinct prime factors of $n$. We first recall a very crucial observation due to Fouvry \cite[Lemme 6]{Fo84a}.

\begin{lemma}\label{lm:Fouvrypartition}
For $\xi\in\bZ^+$ and $M,N \geqslant 2,$ denote
\begin{align*}
\sF_\xi(M,N)&=\{(m,n)\in[1,M]\times[1,N]:(m,n)=1\text { and } \omega(m),\omega(n)\leqslant \xi\}.
\end{align*}
Then there exists a partition of $\sF_\xi(M,N)$ into at most $(2\log(3MN))^{\xi^2}$ subsets $\sF^*_\xi(M,N)$ having the following property:
\begin{align*}
(m_1,n_1),(m_2,n_2) \in \sF^*_\xi(M,N)\Longrightarrow(m_1,n_2)=1.
\end{align*}
\end{lemma}

Let $H,M,N,R,S\geqslant1$ and $a\in\bZ^*.$ For arbitrary complex coefficients $\balpha=(\alpha_{m,r}),$ $\bbeta=(\beta_{h,n,s}),$ we introduce a quintilinear form with Kloosterman fractions
\begin{align}\label{eq:quintilinearform-fraction}
\cQ_a(\balpha,\bbeta)
&=\sum_{h\leqslant H}\mathop{\sum_{m\leqslant M}\sum_{n\leqslant N}\sum_{r\leqslant R}\sum_{s\leqslant S}}_{(mn,rs)=(ns,a)=1}\alpha_{m,r}\beta_{h,n,s}\ue\Big(\frac{ah\overline{mn}}{rs}\Big).
\end{align}
The proof of Theorem \ref{thm:countingcongruence-specialcoefficient} will require the following estimate.
\begin{lemma}\label{lm:quintilinearform-Kloostermanfraction}
Let $H,M,N,R,S,\xi\geqslant1$ and $a\in\bZ^*.$ We assume that $\beta_{h,n,s}=0$ unless $(n,s)\in\sF_\xi(N,S)$ and $\mu(ns)\neq0.$ Then we have
\begin{align*}
\cQ_a(\balpha,\bbeta)
&\ll \|\balpha\|\|\bbeta\|_{\infty}(HNS)^{\frac{1}{2}}((HNS)^{\frac{1}{4}}+(NS)^{\frac{1}{2}})(HMNRS)^{\varepsilon}(2\log (3NS))^{\xi^2}\Bigg[NS(MR)^{\frac{1}{2}}\\
&\ \ \ \ \times\frac{(|a|+\frac{MNRS}{H})^{\vartheta}}{1+(\frac{|a|H}{MNRS})^{\frac{1}{2}}
}\Big(1+\frac{|a|H}{MNRS}+\frac{R}{MN^2}\Big)^{\frac{1}{2}}\Big(1+\frac{|a|H}{MNRS}+\frac{H}{NS}\Big)^{\frac{1}{2}}+S^{-2}M (HNS)^{\frac{1}{2}}\Bigg]^{\frac{1}{2}}\\
&\ \ \ \ +\|\balpha\|\|\bbeta\|_{\infty}(HMNRS)^{\frac{1}{2}+\varepsilon}(2\log (3NS))^{\xi^2},
\end{align*}
where the implied constant depends only on $\varepsilon,\xi,$ and $\vartheta\leqslant 7/64$ is the exponent towards the Ramanujan--Petersson conjecture for $\GL_2(\bQ).$
\end{lemma}

\proof Without loss of generality we assume $\|\bbeta\|_{\infty}\leqslant1.$ 
According to Lemma \ref{lm:Fouvrypartition}, we split the sums over $n,s$ into at most $O((\log NS)^{2\xi})$ sub-sums of the shape 
\begin{align*}
\cQ^*=\cQ_a^*\{\sF^*_\xi(N,S)\}
&=\sum_{h\leqslant H}\mathop{\sum_{m\leqslant M}\sum_{n\leqslant N}\sum_{r\leqslant R}\sum_{s\leqslant S}}_{\substack{(mn,rs)=(ns,a)=1\\ (n,s)\in \sF^*_\xi(N,S)}}\alpha_{m,r}\beta_{h,n,s}\ue\Big(\frac{ah\overline{mn}}{rs}\Big),
\end{align*}
so that
\begin{align}\label{eq:cQ-partition}
\cQ_a(\balpha,\bbeta)
&=\sum_{\sF^*_\xi(N,S)}\cQ_a^*\{\sF^*_\xi(N,S)\}.
\end{align}

Fix $W$ as a smooth function dominating the indicator function of $[1,2].$ By Cauchy--Schwarz inequality, we have
\begin{align*}
|\cQ^*|^2
&\leqslant \|\balpha\|^2\mathop{\sum\sum}_{\substack{r,m\in\bZ\\(r,m)=1}}W\Big(\frac{r}{R}\Big)W\Big(\frac{m}{M}\Big)\Big|\mathop{\sum_{h\leqslant H}\sum_{n\leqslant N}\sum_{s\leqslant S}}_{\substack{(n,rs)=(ns,a)=1\\ (n,s)\in \sF^*_\xi(N,S)}}\beta_{h,n,s}\ue\Big(\frac{ah\overline{mn}}{rs}\Big)\Big|^2.
\end{align*}
Squaring out and switching summations,
\begin{align*}
|\cQ^*|^2
&\leqslant \|\balpha\|^2
\mathop{\sum\sum\sum\sum\sum\sum}_{\substack{h_1,h_2\leqslant H,~n_1,n_2\leqslant N,~s_1,s_2\leqslant S\\ (n_1n_2,s_1s_2)=1,~ (n_1,s_1),(n_2,s_2)\in \sF^*_\xi(N,S)}}\beta_{h_1,n_1,s_1}\overline{\beta_{h_2,n_2,s_2}}\\
&\ \ \ \ \times\mathop{\sum\sum}_{\substack{r,m\in\bZ\\(ar,n_1n_2)=(am,s_1s_2)=(r,m)=1}}W\Big(\frac{r}{R}\Big)W\Big(\frac{m}{M}\Big)\ue\Big(\frac{ah_1\overline{mn_1}}{rs_1}-\frac{ah_2\overline{mn_2}}{rs_2}\Big).
\end{align*}
The diagonal terms, i.e., $h_1n_2s_2=h_2n_1s_1,$ give a contribution at most
\begin{align*}
&\ll \|\balpha\|^2 MR
\mathop{\sum\sum\sum\sum\sum\sum}_{\substack{h_1,h_2\leqslant H,~n_1,n_2\leqslant N,~s_1,s_2\leqslant S\\ h_1n_2s_2=h_2n_1s_1}}1\\
&\ll \|\balpha\|^2 (HMNRS)^{1+\varepsilon}.
\end{align*}
Denote by $\cQ_1$ the contribution from off-diagonal terms, so that
\begin{align}\label{eq:cQ*^2-cQ1}
|\cQ^*|^2
&\leqslant \cQ_1+O(\|\balpha\|^2 (HMNRS)^{1+\varepsilon}).
\end{align}

For tuples $(n_1,s_1),(n_2,s_2)\in \sF_\xi^*(N,S),$ we find $n_1,n_2,s_1,s_2$ are pairwise coprime according to Lemma \ref{lm:Fouvrypartition}. 
Hence the above fraction in the exponential is
\begin{align*}
\frac{a\overline{mn_1n_2}(h_1n_2s_2-h_2n_1s_1)}{rs_1s_2}\bmod1.
\end{align*}
For given $n_1,n_2,s_1,s_2$ with $(n_1n_2,s_1s_2)=1$ and $(n_1,s_1),(n_2,s_2)\in \sF_\xi^*(N,S)$, we introduce a coefficient $\bgamma_{n_1,n_2,s_1,s_2}=(\gamma_{k;n_1,n_2,s_1,s_2})$, indexed by $k$, as
\begin{align*}
\gamma_{k;n_1,n_2,s_1,s_2}
&=
\mathop{\sum\sum}_{\substack{h_1,h_2\leqslant H\\
h_1n_2s_2-h_2n_1s_1=k}}\beta_{h_1,n_1,s_1}\overline{\beta_{h_2,n_2,s_2}}.
\end{align*}
This convention allows us to write
\begin{align*}
\cQ_1
&= \|\balpha\|^2\mathop{\sum\sum\sum\sum}_{\substack{n_1,n_2\leqslant N,~s_1,s_2\leqslant S\\ (n_1n_2,s_1s_2)=1,~(n_1,s_1),(n_2,s_2)\in \sF_\xi^*(N,S)}}\sum_{1\leqslant |k|\leqslant HNS}\gamma_{k;n_1,n_2,s_1,s_2}\\
&\ \ \ \ \times
\mathop{\sum\sum}_{\substack{r,m\in\bZ\\ (ar,n_1n_2)=(am,s_1s_2)=(r,m)=1}}W\Big(\frac{r}{R}\Big)W\Big(\frac{m}{M}\Big)\ue\Big(\frac{ak\overline{mn_1n_2}}{rs_1s_2}\Big).
\end{align*}
We are now in a position to apply Lemma \ref{lm:trilinearKloosterman-inverse} with
\begin{align*}
(a,r,s)\leftarrow (a,n_1n_2,s_1s_2),\ \ 
(c,m,n)\leftarrow (r,m,k),\ \ 
(C,M,N)\leftarrow (R,M,HNS),
\end{align*}
getting
\begin{align}\label{eq:cQ1-Delta}
\cQ_1
&\ll \|\balpha\|^2\mathop{\sum\sum\sum\sum}_{\substack{n_1,n_2\leqslant N,~s_1,s_2\leqslant S\\ (n_1n_2,s_1s_2)=1}}\|\bgamma_{n_1,n_2,s_1,s_2}\|\cdot \Delta,
\end{align}
where
\begin{align*}
\Delta
&=|a|^\vartheta \sqrt{MRS^2N^2} \frac{(1+\frac{MRS^2N^2}{|a|HNS})^{\vartheta}}{1+(\frac{|a|HNS}{MRS^2N^2})^{\frac{1}{2}}}\Big(1+\frac{|a|HNS}{MRS^2N^2}+\frac{R}{MN^2}\Big)^{\frac{1}{2}}\\
&\ \ \ \times\Big(1+\frac{|a|HNS}{MRS^2N^2}+\frac{HNS}{N^2S^2}\Big)^{\frac{1}{2}}
+S^{-2}M (HNS)^{\frac{1}{2}}\\
&=\sqrt{MN^2RS^2} \frac{(|a|+\frac{MNRS}{H})^{\vartheta}}{1+(\frac{|a|H}{MNRS})^{\frac{1}{2}}}\Big(1+\frac{|a|H}{MNRS}+\frac{R}{MN^2}\Big)^{\frac{1}{2}}\Big(1+\frac{|a|H}{MNRS}+\frac{H}{NS}\Big)^{\frac{1}{2}}\\
&\ \ \ \ +S^{-2}M (HNS)^{\frac{1}{2}}.
\end{align*}

\begin{lemma}\label{lm:L1norm}
\begin{align*}
\mathop{\sum\sum\sum\sum}_{\substack{n_1,n_2\leqslant N,~s_1,s_2\leqslant S\\ (n_1n_2,s_1s_2)=1}}\|\bgamma_{n_1,n_2,s_1,s_2}\|
\ll (H(NS)^2+(HNS)^{\frac{3}{2}})(HNS)^\varepsilon.
\end{align*}
\end{lemma}

\proof We first observe that
\begin{align*}
\mathop{\sum\sum\sum\sum}_{\substack{n_1,n_2\leqslant N,~s_1,s_2\leqslant S\\ (n_1n_2,s_1s_2)=1}}\|\bgamma_{n_1,n_2,s_1,s_2}\|^2
&\ll \mathop{\sum\sum\sum\sum\sum\sum\sum\sum}_{\substack{h_1,h_2,h_1',h_2'\leqslant H,~n_1,n_2\leqslant N,~s_1,s_2\leqslant S
\\ (h_1-h_1')n_2s_2=(h_2-h_2')n_1s_1}}1\\
&\ll ((HNS)^2+H^3NS)(HNS)^\varepsilon.
\end{align*}
The lemma then follows from the Cauchy--Schwarz inequality.
\endproof

From \eqref{eq:cQ1-Delta} and Lemma \ref{lm:L1norm}, it follows that
\begin{align*}
\cQ_1
&\ll \|\balpha\|^2HNS((HNS)^{\frac{1}{2}}+NS)\Bigg[NS(MR)^{\frac{1}{2}}\frac{(|a|+\frac{MNRS}{H})^{\vartheta}}{1+(\frac{|a|H}{MNRS})^{\frac{1}{2}}}\Big(1+\frac{|a|H}{MNRS}+\frac{R}{MN^2}\Big)^{\frac{1}{2}}\\
&\ \ \ \ \times\Big(1+\frac{|a|H}{MNRS}+\frac{H}{NS}\Big)^{\frac{1}{2}}+S^{-2}M (HNS)^{\frac{1}{2}}\Bigg](HMNRS)^\varepsilon,
\end{align*}
from which and \eqref{eq:cQ*^2-cQ1}, \eqref{eq:cQ-partition} we conclude Lemma \ref{lm:quintilinearform-Kloostermanfraction}.
\endproof

\subsection{Initial reductions}
We now begin to prove Theorem \ref{thm:countingcongruence-specialcoefficient},
and assume that $ \|\balpha\|_\infty,\|\bbeta\|_\infty\leqslant1$.
We first write
\begin{align*}
\cR
=\sum_{\substack{d\leqslant 2N\\ (d,q)=1}}\mathop{\sum\sum\sum\sum}_{\substack{n_1,n_2\sim N/d,~l_1,l_2\in\bZ\\n_1l_1\equiv n_2l_2\bmod q\\  (n_1n_2l_1l_2,q)=(n_1,n_2)=1}}\alpha_{dn_1}\beta_{dn_2}F(l_1,l_2).
\end{align*}
Note that the congruence restriction $n_1l_1\equiv n_2l_2\bmod q$
is equivalent to 
\begin{align*}
n_1l_1=n_2l_2+qr,\ \ 0\leqslant |r|\leqslant R/d
\end{align*}
with $R:=2NL/q.$ Therefore, 
\begin{align*}
\cR
&=\sum_{\substack{d\leqslant 2N\\ (d,q)=1}}\sum_{0\leqslant |r|\leqslant R/d}\mathop{\sum\sum}_{\substack{n_1,n_2\sim N/d\\  (n_1n_2,q)=(n_1,n_2)=1}}\alpha_{dn_1}\beta_{dn_2}\sum_{\substack{l_1\equiv qr\overline{n_1}\bmod{n_2}\\ (l_1,q)=1}}F\Big(l_1,\frac{n_1l_1-qr}{n_2}\Big).
\end{align*}
From Poisson summation, the last $l_1$-sum can be transformed as
\begin{align*}
\sum_{l_1}
&=\sum_{\delta\mid q}\mu(\delta)\sum_{l_1\equiv qr\overline{\delta n_1}\bmod{n_2}}F\Big(\delta l_1,\frac{\delta n_1l_1-qr}{n_2}\Big)\\
&=\frac{1}{n_2}\sum_{\delta\mid q}\frac{\mu(\delta)}{\delta}\sum_{h\in\bZ}\widecheck{F}\Big(\frac{h}{n_2\delta}; n_1,n_2,r\Big)\ue\Big(\frac{qrh\overline{\delta n_1}}{n_2}\Big),
\end{align*}
where
\begin{align*}
\widecheck{F}(\lambda; n_1,n_2,r)=\int_\bR F\Big(y,\frac{n_1y-qr}{n_2}\Big)\ue(-\lambda y)\ud y.
\end{align*}
Note that the integration by parts yields
\begin{align}\label{eq:F-Fourier}
\sup_{n_1,n_2,r}|\widecheck{F}(\lambda; n_1,n_2,r)|&\ll L(1+\eta^{-1}|\lambda|L)^{-A}
\end{align}
for all $A\geqslant0$.
Hence
\begin{align*}
\cR
&=\sum_{\substack{d\leqslant 2N\\ (d,q)=1}}\sum_{0\leqslant |r|\leqslant R/d}\mathop{\sum\sum}_{\substack{n_1,n_2\sim N/d\\(n_1n_2,q)=(n_1,n_2)=1}}\frac{\alpha_{dn_1}\beta_{dn_2}}{n_2}\sum_{\delta\mid q}\frac{\mu(\delta)}{\delta}\\
&\ \ \ \ \times \sum_{h\in\bZ}\widecheck{F}\Big(\frac{h}{n_2\delta}; n_1,n_2,r\Big)\ue\Big(\frac{qrh\overline{\delta n_1}}{n_2}\Big).
\end{align*}
The term with $h=0$ gives exactly
\begin{align}\label{eq:tidleR-h=0}
&=\frac{\varphi(q)}{q}\sum_{\substack{d\leqslant 2N\\ (d,q)=1}}\mathop{\sum\sum}_{\substack{n_1,n_2\sim N/d\\(n_1n_2,q)=(n_1,n_2)=1}}\frac{\alpha_{dn_1}\beta_{dn_2}}{n_2}\sum_{0\leqslant |r|\leqslant R/d}\widecheck{F}(0; n_1,n_2,r).
\end{align}
Note that the $r$-sum is 
\begin{align*}
&=n_2\sum_{r\in\bZ}\int_\bR F\Big(n_2y,n_1y-\frac{qr}{n_2}\Big)\ud y=n_2\iint_{\bR^2} F\Big(n_2y,n_1y-\frac{qz}{n_2}\Big)\ud y\ud [z].
\end{align*}
The above Stieltjes integral is equal to 
\begin{align*}
&=n_2\iint_{\bR^2} F\Big(n_2y,n_1y-\frac{qz}{n_2}\Big)\ud y\ud z-q\iint_{\bR^2}\{z\} F'\Big(n_2y,n_1y-\frac{qz}{n_2}\Big)\ud y\ud z,
\end{align*}
where $F'(u,v)=\frac{\partial}{\partial v}F(u,v)$. The second term contributes to \eqref{eq:tidleR-h=0} at most $O(\eta (NL)^{1+\varepsilon}),$ and the first term is exactly
\begin{align*}
\frac{n_2}{q}\iint_{\bR^2} F(y,z)\ud y\ud z.
\end{align*}
Hence
\begin{align}\label{eq:R=R*+...}
\cR=\cR^*+\frac{\varphi(q)}{q^2}\mathop{\sum\sum}_{\substack{n_1,n_2\sim N\\(n_1n_2,q)=1}}\alpha_{n_1}\beta_{n_2}\iint_{\bR^2} F(y,z)\ud y\ud z+O(\eta (NL)^{1+\varepsilon}),
\end{align}
where
\begin{align*}
\cR^*
&=\sum_{\delta\mid q}\frac{\mu(\delta)}{\delta}\sum_{\substack{d\leqslant 2N\\ (d,q)=1}}\sum_{0\leqslant |r|\leqslant R/d}\mathop{\sum\sum}_{\substack{n_1,n_2\sim N/d\\(n_1n_2,q)=(n_1,n_2)=1}}\frac{\alpha_{dn_1}\beta_{dn_2}}{n_2}\\
&\ \ \ \ \times\sum_{h\neq0}\widecheck{F}\Big(\frac{h}{n_2\delta}; n_1,n_2,r\Big)\ue\Big(\frac{qrh\overline{\delta n_1}}{n_2}\Big).
\end{align*}

Denote by $\cR^*_1$ the contribution to $\cR^*$ from the term with $r=0$.
Hence
\begin{align*}
\cR_1^*
&=\sum_{\delta\mid q}\frac{\mu(\delta)}{\delta}\sum_{\substack{d\leqslant 2N\\ (d,q)=1}}\mathop{\sum\sum}_{\substack{n_1,n_2\sim N/d\\(n_1n_2,q)=(n_1,n_2)=1}}\frac{\alpha_{dn_1}\beta_{dn_2}}{n_2}\sum_{h\neq0}\widecheck{F}\Big(\frac{h}{n_2\delta}; n_1,n_2,0\Big).
\end{align*}
Back to the situation before Poisson summation, we find
\begin{align*}
\cR_1^*
&=\sum_{\substack{d\leqslant 2N\\ (d,q)=1}}\mathop{\sum\sum}_{\substack{n_1,n_2\sim N/d\\  (n_1n_2,q)=(n_1,n_2)=1}}\alpha_{dn_1}\beta_{dn_2}\sum_{\substack{l_1\equiv 0\bmod{n_2}\\ (l_1,q)=1}}F\Big(l_1,\frac{n_1l_1}{n_2}\Big)\\
&\ \ \ \ -\sum_{\delta\mid q}\frac{\mu(\delta)}{\delta}\sum_{\substack{d\leqslant 2N\\ (d,q)=1}}\mathop{\sum\sum}_{\substack{n_1,n_2\sim N/d\\(n_1n_2,q)=(n_1,n_2)=1}}\frac{\alpha_{dn_1}\beta_{dn_2}}{n_2}\widecheck{F}(0; n_1,n_2,0)\\
&=\sum_{\substack{d\leqslant 2N\\ (d,q)=1}}\mathop{\sum\sum}_{\substack{n_1,n_2\sim N/d\\  (n_1n_2,q)=(n_1,n_2)=1}}\alpha_{dn_1}\beta_{dn_2}\sum_{(l_1,q)=1}F(n_2l_1,n_1l_1)\\
&\ \ \ \ -\sum_{\delta\mid q}\frac{\mu(\delta)}{\delta}\sum_{\substack{d\leqslant 2N\\ (d,q)=1}}\mathop{\sum\sum}_{\substack{n_1,n_2\sim N/d\\(n_1n_2,q)=(n_1,n_2)=1}}\frac{\alpha_{dn_1}\beta_{dn_2}}{n_2}\widecheck{F}(0; n_1,n_2,0)\\
&\ll (NL)^{1+\varepsilon},
\end{align*}
so that
\begin{align}\label{eq:R*=R*1+R*2}
\cR^*
&=\cR^*_2+O(\eta (NL)^{1+\varepsilon}),
\end{align}
where
\begin{equation}\label{eq:R*2}
\begin{split}
\cR^*_2
&=\sum_{\delta\mid q}\frac{\mu(\delta)}{\delta}\sum_{\substack{d\leqslant 2N\\ (d,q)=1}}\sum_{1\leqslant |r|\leqslant R/d}\mathop{\sum\sum}_{\substack{n_1,n_2\sim N/d\\(n_1n_2,q)=(n_1,n_2)=1}}\frac{\alpha_{dn_1}\beta_{dn_2}}{n_2}\sum_{h\neq0}\widecheck{F}\Big(\frac{h}{n_2\delta}; n_1,n_2,r\Big)\ue\Big(\frac{qrh\overline{\delta n_1}}{n_2}\Big).
\end{split}
\end{equation}

\subsection{Estimates for exponential sums in $\cR^*_2$}

To apply Lemma \ref{lm:quintilinearform-Kloostermanfraction}, we split the $r$-sum in $\cR^*_2$ into dyadic segments and utilize the convolution structure of $\balpha,\bbeta$, i.e., for some $M_1,N_1,M_2,N_2\geqslant1$ with $M_1N_1=M_2N_2\leqslant 2N$ and $\balpha=\balpha'*\balpha'',$ $\bbeta=\bbeta'*\bbeta''$, we have
\begin{align*}
\cR^*_2
&\ll(qNL)^\varepsilon \sum_{\delta\mid q}\frac{1}{\delta}\sum_{\substack{d_1d_2=d_1'd_2'=d\leqslant 2N\\ (d,q)=1}}\sup|\cR^*_2(A_1,B_1,A_2,B_2,R';d_1,d_2,d_1',d_2',\delta)|+1,
\end{align*}
where the supremum is over 
\begin{align*}
1\leqslant A_1\leqslant \frac{M_1}{d_1},~1\leqslant B_1\leqslant \frac{N_1}{d_2}, ~1\leqslant A_2\leqslant \frac{M_2}{d_1'},~1\leqslant B_2\leqslant \frac{N_2}{d_2'}, ~1\leqslant R'\leqslant \frac{R}{d},
\end{align*}
and
\begin{align*}
\cR^*_2(\cdots)
&=\sum_{|r|\sim R'}\mathop{\sum\sum\sum\sum}_{\substack{m_1\sim A_1,~n_1\sim B_1,~m_2\sim A_2,~n_2\sim B_2\\(m_1m_2n_1n_2,q)=(m_1n_1,m_2n_2)=(m_1,d_2)=(m_2,d_2')=1}}\frac{\alpha'_{d_1m_1}\alpha''_{d_2n_1}\beta'_{d_1'm_2}\beta''_{d_2'n_2}}{m_2n_2}\\
&\ \ \ \ \times\sum_{1\leqslant |h|\leqslant H}\widecheck{F}\Big(\frac{h}{m_2n_2\delta}; m_1n_1,m_2n_2,r\Big)\ue\Big(\frac{(q/\delta) hr\overline{m_1n_1}}{m_2n_2}\Big).
\end{align*}
Here we take $H=\eta \delta A_2B_2L^{-1+\varepsilon}.$

We would like explore the cancellations among the exponentials $\ue(\frac{(q/\delta)hr\overline{m_1n_1}}{m_2n_2})$ using Lemma \ref{lm:quintilinearform-Kloostermanfraction}.
It is now necessary to separate variables. To do so, we introduce the Fourier transform of $F$:
\begin{align*}
\widehat{F}(y_1,y_2)
&=\iint_{\bR^2}F(\lambda_1,\lambda_2)\ue(-\lambda_1y_1-\lambda_2y_2)\ud \lambda_1\ud \lambda_2.
\end{align*}
The Fourier inversion gives
\begin{align*}
F(\lambda_1,\lambda_2)
&=\iint_{\bR^2}\widehat{F}(y_1,y_2)\ue(\lambda_1y_1+\lambda_2y_2)\ud y_1\ud y_2.
\end{align*}
Note that for all $\nu_1,\nu_2\in \bN,$
\begin{align}\label{eq:F-Fouriertransform-decay}
\frac{\partial^{\nu_1+\nu_2}}{\partial y_1^{\nu_1}\partial y_2^{\nu_2}}\widehat{F}(y_1,y_2)\ll L^{\nu_1+\nu_2+2}(1+\eta^{-1}L|y_1|)^{-c_1}(1+\eta^{-1}L|y_2|)^{-c_2}
\end{align}
with any $c_1,c_2\geqslant0$, and with an implied constant depending only on $(\nu_1,\nu_2,c_1,c_2)$.

We now start to separate variables in $\widecheck{F}(\frac{h}{m_2n_2\delta}; m_1n_1,m_2n_2,r).$
Again by Fourier inversion,
\begin{align*}
&\ \ \ \ \widecheck{F}\Big(\frac{h}{m_2n_2\delta}; m_1n_1,m_2n_2,r\Big)\\
&=m_2n_2\delta\int_{\frac{rL}{8A_2B_2R'\delta}}^{\frac{2rL}{A_2B_2R'\delta}} F\Big(m_2n_2y\delta,m_1n_1y\delta-\frac{qr}{m_2n_2}\Big)\ue(-h y)\ud y\\
&=m_2n_2\delta\int_{\frac{rL}{8A_2B_2R'\delta}}^{\frac{2rL}{A_2B_2R'\delta}}\ue(-h y)\ud y\iint_{\bR^2} \widehat{F}(y_1,y_2)\ue\Big(m_2n_2y_1y\delta+m_1n_1y_2y\delta-\frac{qry_2}{m_2n_2}\Big)\ud y_1\ud y_2\\
&=\frac{m_2n_2\delta}{r^2}\int_{\frac{rL}{8A_2B_2R'\delta}}^{\frac{2rL}{A_2B_2R'\delta}}\ue(-h y)\ud y\iint_{\bR^2} \widehat{F}\Big(\frac{m_1y_1}{n_2r},\frac{y_2n_2}{m_1r}\Big)\ue\Big(\frac{(m_1m_2y_1+n_1n_2y_2)y\delta}{r}-\frac{qy_2}{m_1m_2}\Big)\ud y_1\ud y_2\\
&=\frac{m_2n_2\delta}{r}\int_{\frac{L}{8A_2B_2R'\delta}}^{\frac{2L}{A_2B_2R'\delta}}\ue(-hry)\ud y\iint_{\bR^2} \widehat{F}\Big(\frac{m_1y_1}{n_2r},\frac{y_2n_2}{m_1r}\Big)\ue\Big((m_1m_2y_1+n_1n_2y_2)y\delta-\frac{qy_2}{m_1m_2}\Big)\ud y_1\ud y_2.
\end{align*}

We need further to separate variables in $\widehat{F}(*,*)$ by Mellin inversion. To this end, we truncate the integral over $y_1$ over five segments 
$$I_1=]-\infty,-V_1],~I_2=[-V_1,-U],~I_3=[-U,U],~I_4=[U,V_1],~I_5=[V_1,+\infty[$$
with
\begin{align*}
U=\frac{1}{N^3L}, \ \ V_1=\frac{\eta R'B_2}{A_1L}L^\varepsilon.
\end{align*}
In view of the decay \eqref{eq:F-Fouriertransform-decay} of $\widehat{F},$ we find the contributions from $I_1,I_5$ are negligible, and we may also bound the integral over $I_3$ trivially, contributing to $\cR^*_2(\cdots)$ at most
\begin{align*}
\ll\frac{\eta^2\delta N^\varepsilon}{N^3}\frac{A_1^2B_1}{B_2}\ll \eta^2\delta d^{-1} N^\varepsilon,
\end{align*}
which is also negligible. It thus suffices to treat the integral over $y_1\in I_2\cup I_4$. In a similar manner, we may also reduce the integral over $y_2$ to $[-V_2,-U]\cup[U,V_2]$, where
\begin{align*}
V_2=\frac{\eta R'A_1}{B_2L}L^\varepsilon.
\end{align*}
By symmetry, we only transform the integrals over $y_1\in[U,V_1]$ and $y_2\in[U,V_2]$. To do so, we attach an auxiliary function $\varDelta(\cdot,\cdot)$ that is the indicator function of $[\frac{A_1U}{B_2R'},\frac{A_1V_1}{B_2R'}]\times [\frac{A_1U}{B_2R'},\frac{A_1V_2}{B_2R'}]$.
To ease the presentation, we also denote the new contributions by $\widecheck{F}(\frac{h}{m_2n_2\delta}; m_1n_1,m_2n_2,r)$:
\begin{align*}
\widecheck{F}\Big(\frac{h}{m_2n_2\delta}; m_1n_1,m_2n_2,r\Big)
&=\frac{m_2n_2\delta}{r}\int_{\frac{L}{8A_2B_2R'\delta}}^{\frac{2L}{A_2B_2R'\delta}}\ue(-hry)\ud y \iint_{[U,V_1]\times[U,V_2]} (\varDelta\widehat{F})\Big(\frac{m_1y_1}{n_2r},\frac{y_2n_2}{m_1r}\Big)\\
&\ \ \ \ \ \times\ue\Big((m_1m_2y_1+n_1n_2y_2)y\delta-\frac{qy_2}{m_1m_2}\Big)\ud y_1\ud y_2.
\end{align*}

Note that, for $w_1,w_2>0,$ the Mellin inversion gives
\begin{align*}
(\varDelta\widehat{F})(w_1,w_2)
&=\frac{1}{4\pi^2}\iint_{\bR^2}\widetilde{F}(it_1,it_2)w_1^{-it_1}w_2^{-it_2}\ud t_1\ud t_2,
\end{align*}
where 
\begin{align*}
\widetilde{F}(s_1,s_2) 
&=\iint_{(\bR^+)^2}(\varDelta\widehat{F})(w_1,w_2)w_1^{s_1-1}w_2^{s_2-1}\ud w_1\ud w_2.
\end{align*}
Hence
\begin{align*}
&\ \ \ \ \widecheck{F}\Big(\frac{h}{m_2n_2\delta}; m_1n_1,m_2n_2,r\Big)\\
&=\frac{m_2n_2\delta}{4\pi^2 r}\int_{\frac{L}{8A_2B_2R'\delta}}^{\frac{2L}{A_2B_2R'\delta}}\ue(-hry)\ud y
\int_U^{V_1} \int_U^{V_2} \ue\Big((m_1m_2y_1+n_1n_2y_2)y\delta-\frac{qy_2}{m_1m_2}\Big)\ud y_1\ud y_2\\
&\ \ \ \ \times \iint_{\bR^2}\widetilde{F}(it_1,it_2)\Big(\frac{n_2r}{m_1y_1}\Big)^{it_1}\Big(\frac{m_1r}{y_2n_2}\Big)^{it_2}\ud t_1\ud t_2.
\end{align*}

Thanks to the separation of variables as above, we get
\begin{align*}
\cR^*_2(\cdots)
&=\frac{\delta}{4\pi^2}\iint_{\bR^2}\widetilde{F}(it_1,it_2)\ud t_1\ud t_2
\int_{\frac{L}{8A_2B_2R'\delta}}^{\frac{2L}{A_2B_2R'\delta}}\ud y
\int_U^{V_1}  y_1^{-it_1}\ud y_1
\int_U^{V_2} y_2^{-it_2}\ud y_2\\
&\ \ \ \ \times \sum_{1\leqslant |h|\leqslant H}\sum_{|r|\sim R'}\ue(-hry)r^{i(t_1+t_2)-1}\\
&\ \ \ \ \times \mathop{\sum\sum\sum\sum}_{\substack{m_1\sim A_1,~n_1\sim B_1,~m_2\sim A_2,~n_2\sim B_2\\(m_1m_2n_1n_2,q)=(m_1n_1,m_2n_2)=(m_1,d_2)=(m_2,d_2')=1}}\alpha'_{d_1m_1}\alpha''_{d_2n_1}\beta'_{d_1'm_2}\beta''_{d_2'n_2}\\
&\ \ \ \ \times \Big(\frac{n_2}{m_1}\Big)^{i(t_1-t_2)}\ue\Big((m_1m_2y_1+n_1n_2y_2)y\delta-\frac{qy_2}{m_1m_2}\Big)
\ue\Big(\frac{(q/\delta) hr\overline{m_1n_1}}{m_2n_2}\Big).
\end{align*}

Note that the integration by parts and \eqref{eq:F-Fouriertransform-decay} yield
\begin{align*}
\widetilde{F}(it_1,it_2)\ll (L\log L)^2(1+\eta^{-1}|t_1|)^{-c_1}(1+\eta^{-1}|t_2|)^{-c_2}
\end{align*}
with any $c_1,c_2\geqslant0$.
This implies
\begin{align*}
\iint_{\bR^2} |\widetilde{F}(it_1,it_2)|\ud t_1\ud t_2\ll (\eta L\log L)^2.
\end{align*}
We are now in a good position to apply Lemma \ref{lm:quintilinearform-Kloostermanfraction} with
\begin{align*}
(a,h,m,n,r,s)&\leftarrow (q/\delta,k,m_1,n_1,m_2,n_2),\\
(H,M,N,R,S,\xi)&\leftarrow (R'H,A_1,B_1,A_2,B_2,(\log N)^{\frac{1}{3}}),\\
\alpha_{m,r}&\leftarrow \alpha'_{d_1m_1}\beta'_{d_1'm_2}m_1^{i(t_2-t_1)}\ue\Big(m_1m_2y_1y\delta-\frac{qy_2}{m_1m_2}\Big),\\
\beta_{h,n,s}&\leftarrow \alpha''_{d_2n_1}\beta''_{d_2'n_2}n_2^{i(t_1-t_2)}\ue(n_1n_2y_2y\delta)\mathop{\sum\sum}_{\substack{hr=k\\ 1\leqslant |h|\leqslant H,~|r|\sim R'}}\ue(-hry)r^{i(t_1+t_2)-1},
\end{align*}
so that
\begin{align*}
\cR^*_2
&\ll(qNL)^\varepsilon \eta^2L\Big(\frac{\eta N^3}{q}\Big)^{\frac{1}{2}}\Big[\Big(\frac{\eta N^2N_1N_2}{q}\Big)^{\frac{1}{4}}+(N_1N_2)^{\frac{1}{2}}\Big]\Big[\eta q^{\vartheta}(N_1N_2)^{\frac{1}{2}}\Big(1+\frac{N^2}{qN_1N_2}\Big)^{\frac{1}{2}}\Big]^{\frac{1}{2}}\\
&\ \ \ \  +\eta^4 NL(\eta N^2/qN_1N_2)^{\frac{1}{2}+\varepsilon}\\
&\ll \eta^{\frac{9}{2}}N^{\frac{3}{2}}Lq^{\frac{\vartheta-1}{2}}(N_1N_2)^{\frac{1}{2}}(N^2q^{-1}+N_1N_2)^{\frac{1}{4}}\Big(1+\frac{N^2}{qN_1N_2}\Big)^{\frac{1}{4}}(qNL)^{\varepsilon}+\eta^{\frac{9}{2}}N^2L(qN_1N_2)^{-\frac{1}{2}}(qNL)^\varepsilon,
\end{align*}
from which and 
\eqref{eq:R*=R*1+R*2}, \eqref{eq:R=R*+...}, 
we then conclude the proof of Theorem \ref{thm:countingcongruence-specialcoefficient}.

Moreover, in the settings of Theorem \ref{thm:countingcongruence-generalcoefficient}, the coefficients $\balpha,\bbeta$ are arbitrary, and we appeal to Lemma \ref{lm:trilinearform-BC} instead to bound the sums with Kloosterman fractions in \eqref{eq:R*2}.
To this end, one may argue as above to group the variables, and the details are omitted here.

\smallskip

\section{Comments and remarks}
\label{sec:comments}

\subsection{Twisted fourth moment of Dirichlet $L$-functions}
As a (weaker) substitute for the Lindel\"of Hypothesis, it is a reasonable approach to consider moments of Dirichlet $L$-functions.  It has been known for a long time that
\begin{align}\label{eq:Lfourthmoment}
\frac{1}{\varphi(q)}\sum_{\substack{\chi\bmod{q}\\ \chi\neq \chi_{0}}}|L(\tfrac{1}{2}+it, \chi)|^4
\ll q^\varepsilon (1+|t|)^{O(1)}
\end{align}
for any $\varepsilon>0$ and $t\in\bR$. In addition, we may also expect, for any complex sequence  $\balpha=(\alpha_m)$  and any $\varepsilon>0$, that
\begin{align}\label{eq:Hypothesis-Lfourthmoment}
\frac{1}{\varphi(q)}\sum_{\substack{\chi\bmod{q}\\ \chi\neq \chi_{0}}}\Big|\sum_{m\leqslant q^\delta} \alpha_m \chi(m)\Big|^2|L(\tfrac{1}{2}+it, \chi)|^4
\ll (q(1+|t|))^{\varepsilon} q^\delta\|\balpha\|_\infty 
\end{align}
with $\delta>0$ as large as possible.

By extending Motohashi's identity to link fourth moment of Dirichlet $L$-functions and cubic moment of $\GL_2$ $L$-functions, Blomer, Humphries, Khan and Milinovich \cite{BHKM20} prove that \eqref{eq:Hypothesis-Lfourthmoment}
is valid as long as $q$ is prime, $t=0$, $\delta\in[0,1/4]$ and $\balpha$ is supported on squarefree numbers. Their bound can, together with a suitable amplification, lead to a new proof of Burgess's bound $L(\frac{1}{2},\chi)\ll q^{3/16+\varepsilon}$ for prime $q$.
In applications to Brun--Titchmarsh theorem, it is highly desirable to extend \eqref{eq:Hypothesis-Lfourthmoment} for a larger admissible $\delta.$
Also, it is desirable to extend \cite{BHKM20} to all moduli $q$, which then readily implies that the choice \eqref{eq:C*(varpi):varpi<1/2,primemoduli} is admissible in  \eqref{eq:C(varpi):varpi<1/2} for all such $q$.

Recall that we have utilized \eqref{eq:Lfourthmoment} in the proof of Theorem \ref{thm:BT-smoothmoduli<1/2}. If one were allowed to use \eqref{eq:Hypothesis-Lfourthmoment} with smooth and squarefree $q$, it would be possible to obtain the desirable constant $2$ in Theorem \ref{thm:BT-smoothmoduli<1/2} in a larger range of $\varpi$.

\subsection{Brun--Titchmarsh theorem on average}

Hooley \cite{Ho72} was the first who was able to study upper bounds of $\pi(x;q,a)$ 
non-trivially through extra averaging over $q$. This was largely motivated by the study of greatest prime factors of shifted primes, as approximations to the twin prime conjecture. 
For an individual $q$, one may also consider the average over the residue class $a\bmod q$. See also \cite{Mo74,Ho75,Mo79,Iw82,Fo84b,Fo85,BH96,FR22} for these developments.

\subsection{Extending Theorem $\ref{thm:BT-generalmoduli>1/2}$ to general moduli}

In his ETH thesis, L\"offel \cite{Lo16} generalized Lemma \ref{lm:FKM}, as well as other works of Fouvry, Kowalski and Michel \cite{FKM14,FKM15}, to squarefree moduli. This also allows us to relax $q$ in Theorem $\ref{thm:BT-generalmoduli>1/2}$ to be squarefree. It is also possible to generalize Theorem $\ref{thm:BT-generalmoduli>1/2}$ to arbitrary moduli along the same lines. As one may see from the proof of Lemma \ref{lm:FKM}, it is necessary to evaluate averages of products of four Kloosterman sums to general moduli.

In another direction, Lin and Michel \cite{LM24} considered a variant of Lemma \ref{lm:FKM} with modulus $q$, which is exactly a product of two distinct primes. When the two prime factors are of suitable sizes, they are able to capture stronger non-correlations between algebraic trace functions and Hecke eigenvalues (holomorphic or Maass). By extending their arguments to the situation of Eisenstein spectrum, it is also possible to study the Brun--Titchmarsh theorem with factorable moduli, which can be compared with the direct application of the method of arithmetic exponent pairs.

\subsection{Twisted second moment of character sums}\label{subsec:Twistedsecond momentcharacter}
Recall that Theorems \ref{thm:countingcongruence-generalcoefficient}, \ref{thm:countingcongruence-specialcoefficient} and Lemmas \ref{lm:weightedL-secondmoment-generalcoefficient}, \ref{lm:weightedL-secondmoment-specialcoefficient}
are intimately related to the twisted second moment of character sums
\begin{align*}
\cZ:=\frac{1}{\varphi(q)}\sum_{\substack{\chi\bmod q\\ \chi\neq\chi_0}}\Big|\sum_{m\leqslant M} c_m\chi(m)\Big|^2\Big|\sum_{l\leqslant L} \chi(l)\Big|^2
\end{align*}
with $|c_m|\ll m^\varepsilon.$
We expect that 
$\cZ\ll MLq^\varepsilon$ holds with $M$ as large as possible, building on the square-root cancellation philosophy.
\begin{itemize}
\item  Friedlander and Iwaniec \cite{FI92} succeeded for $M\leqslant q^{1/2+1/22-}$ and $(c_m)$ is a convolution of two coefficients of level $q^{3/11-}.$
\item Duke, Friedlander and Iwaniec \cite{DFI97a} succeeded for $M\leqslant q^{1/2+1/190-}$ and $c_m$ is general.
\item Lemma \ref{lm:weightedL-secondmoment-generalcoefficient}  succeeds for $M\leqslant q^{1/2+1/66-}$ and $c_m$ is general.
\item Lemma \ref{lm:weightedL-secondmoment-specialcoefficient} succeeds for $M\leqslant q^{(5-2\vartheta)/8-}$ and $c_m$ is well-factorable.

\end{itemize}

Recall that Theorem \ref{thm:countingcongruence-generalcoefficient} and Lemma \ref{lm:weightedL-secondmoment-generalcoefficient} are proven with the aid of Lemma 
\ref{lm:trilinearform-BC} on trilinear forms with Kloosterman fractions. However,
Theorem \ref{thm:countingcongruence-specialcoefficient} and Lemma \ref{lm:weightedL-secondmoment-specialcoefficient}
rely on the control of the quintilinear form in \eqref{eq:quintilinearform-fraction}:
\begin{align*}
\cQ_a(\balpha,\bbeta)
&=\sum_{h\leqslant H}\mathop{\sum_{m\leqslant M}\sum_{n\leqslant N}\sum_{r\leqslant R}\sum_{s\leqslant S}}_{(mn,rs)=(ns,a)=1}\alpha_{m,r}\beta_{h,n,s}\ue\Big(\frac{ah\overline{mn}}{rs}\Big).
\end{align*}
Due to the application of Kloostermania, the size of $a$ plays an important role in the settings of Theorem \ref{thm:countingcongruence-specialcoefficient}. Note that Lemma \ref{lm:trilinearform-BC} follows from a careful application of amplification, so that only elementary means together with Weil's bound for complete algebraic exponential sums are needed. A natural problem is to ask for an elementary bound for $\cQ_a(\balpha,\bbeta)$ which is uniform in $a\neq0$, and $\balpha,\bbeta$ are of convolution structures if necessary.

\subsection{Hooley's Hypothesis R$^*$ and its impacts}
In his studies on the Brun--Titchmarsh theorem and greatest prime factors of cubic polynomials, Hooley \cite{Ho72,Ho78} proposed the following conjecture:

\begin{conjecture}[Hypothesis R$^*$]\label{conj:HypothesisR*}
Let $q$ be a positive integer and $h\in\bZ.$ Let $I$ be an interval with $1<|I|\leqslant q.$ Then
\begin{align*}
\sum_{\substack{n\in I\\ (n,q)=1}}\ue\Big(\frac{h\overline{n}}{q}\Big)
&\ll |I|^{\frac{1}{2}}(h,q)^{\frac{1}{2}}q^\varepsilon
\end{align*}
for any $\varepsilon>0.$
\end{conjecture}
Under Hypothesis R$^*$, Iwaniec \cite[Theorems 8 \& 9]{Iw82} showed the choices
\begin{align}\label{eq:BT-HypothesisR*}
C(\varpi)=
\begin{cases}
6/(5-6\varpi), \quad &\varpi\in]4/9,~7/12[,\\
5/(3-3\varpi), &\varpi\in]7/12,~1[,
\end{cases}
\end{align}
are admissible in \eqref{eq:BT}.
Hypothesis R$^*$ exhibits square-root cancellations within incomplete Kloosterman sums. It seems far beyond the current techniques to prove Conjecture \ref{conj:HypothesisR*} for general $q$, even in the case $|I|\approx q^{0.99}.$
On the other hand, Conjecture \ref{conj:HypothesisR*} corresponds to the exponent pair hypothesis in the theory of arithmetic exponent pairs, for which we hope $(\kappa,\lambda)=(0,\frac{1}{2})$ is admissible for \eqref{eq:exponentpair-estimate}.
This illustrates that all of our previous arguments starting from additive characters should never beat \eqref{eq:BT-HypothesisR*} unconditionally, even when imposing factorization conditions on $q$.

However, as $\varpi\rightarrow 9/20-$ for instance, Theorem \ref{thm:BT-smoothmoduli<1/2}
shows that $C(\varpi)=5/(5-6\varpi)$ is admissible, as long as $q$ is squarefree and smooth. This is obviously sharper than the conjectural bound \eqref{eq:BT-HypothesisR*} which takes $C(\varpi)=6/(5-6\varpi).$
As one may see from previous sections, Theorem \ref{thm:BT-smoothmoduli<1/2} is proved  using the orthogonality of multiplicative characters, and hence Lemma \ref{lm:incompletecharactersum-smoothmoduli} on character sums can be applied to complete the argument. This phenomenon indicates that for $\varpi<1/2-$ it is clever to use multiplicative characters, and for $\varpi>1/2+$ additive ones might be more powerful instead. The exact transition range is very mysterious to the authors.

\bibliographystyle{plain}

\begin{thebibliography}{abcdefghi}


\bibitem[Ba96]{Ba96}
R. C. Baker,
The Brun--Titchmarsh theorem,
\emph{J. Number Theory} \textbf{56} (1996), 343--365.



\bibitem[BH96]{BH96}
R. C. Baker \& G. Harman,
The Brun--Titchmarsh theorem on average,
Analytic Number Theory, Vol. 1 (Allerton Park, IL, 1995),
Progr. Math. 138, Birkh\"auser, Boston, 1996,  39--103.




\bibitem[BC18]{BC18}
S. Bettin \& V. Chandee,
Trilinear forms with Kloosterman fractions,
\emph{Adv. Math.} \textbf{328} (2018), 1234--1262.



\bibitem[BHKM20]{BHKM20}
V. Blomer, P. Humphries, R. Khan \& M. B. Milinovich, 
Motohashi's fourth moment identity for non-archimedean test functions and applications, 
\emph{Compos. Math.} \textbf{156} (2020), 1004--1038.  

\bibitem[BG14]{BG14}
J. Bourgain \& M. Z. Garaev,
Kloosterman sums in residue rings,
\emph{Acta Arith.} \textbf{164} (2014), 43--64.

\bibitem[Br86]{Br86} 
J.-L. Brylinski, 
Transformations canoniques, dualit\'{e} projective, th\'{e}{e}orie de Lefschetz, 
transformations de Fourier et sommes trigonom\'{e}triques, 
\emph{Ast\'{e}risque} {\bf 140-141} (1986), 3--134.



\bibitem[BPZ20]{BPZ20} 
H. M. Bui, K. Pratt \& A. Zaharescu,
Breaking the $1/2$-barrier for the twisted second
moment of Dirichlet $L$-functions,
\emph{Adv. Math.} \textbf{370} (2020) 107--175.


\bibitem[Bu63]{Bu63} 
D. A. Burgess, 
On character sums and $L$-series. II, 
\emph{Proc. London Math. Soc.} \textbf{13} (1963), 524--536.


\bibitem[Bu86]{Bu86} 
D. A. Burgess, 
The character sum estimate with $r = 3,$
\emph{J. London Math. Soc.} \textbf{33} (1986), 219--226.


\bibitem[Ch75]{Ch75}
J. R. Chen,
On the distribution of almost primes in an interval,
\emph{Sci. Sinica} \textbf{18} (1975), 611--627.


\bibitem[De80]{De80}
P. Deligne, 
La conjecture de Weil II, 
\emph{Publ. Math. IH\'{E}S} \textbf{52} (1980), 137--252.



\bibitem[DI82]{DI82} 
J.-M. Deshouillers \& H. Iwaniec, 
Kloosterman sums and Fourier coefficients of cusp forms, 
\emph{Invent. Math.} \textbf{70} (1982/83), 171--188.

\bibitem[DFI97a]{DFI97a} 
W. Duke, J. B. Friedlander \& H. Iwaniec,
Representations by the determinant and mean values of L-functions.
Sieve Methods, Exponential Sums, and Their Applications in Number Theory, Cardiff, 1995, London Math. Soc. Lecture Note Ser., vol. 237, 1997, pp. 109--115.


\bibitem[DFI97b]{DFI97b} 
W. Duke, J. B. Friedlander \& H. Iwaniec,
Bilinear forms with Kloosterman fractions, 
\emph{Invent. Math.} \textbf{128} (1997), 23--43.



\bibitem[Fo84a]{Fo84a}
\'{E}. Fouvry, 
Autour du th\'eor\`eme de Bombieri--Vinogradov,
\emph{Acta Math.} \textbf{152} (1984), 219--244.

\bibitem[Fo84b]{Fo84b}
\'{E}. Fouvry, 
Sur le th\'eor\`eme de Brun--Titchmarsh,
\emph{Acta Arith.} \textbf{43} (1984), 417--424.



\bibitem[Fo85]{Fo85}
\'{E}. Fouvry, 
Th\'eor\`eme de Brun--Titchmarsh: application au th\'eor\`eme de Fermat,
\emph{Invent. Math.} \textbf{79} (1985), 383--407. 

\bibitem[FGKM14]{FGKM14} 
\'E. Fouvry, S. Ganguly, E. Kowalski \& Ph. Michel, 
Gaussian distribution for the divisor function and Hecke eigenvalues in arithmetic progressions, 
\emph{Comment. Math. Helv.} \textbf{89} (2014), 979--1014.


\bibitem[FKM14]{FKM14}
\'{E}. Fouvry, E. Kowalski \& Ph. Michel, 
Algebraic trace functions over the primes, 
\emph{Duke Math. J.} \textbf{163} (2014), 1683--1736.



\bibitem[FKM15]{FKM15} 
\'{E}. Fouvry, E. Kowalski \& Ph. Michel, 
Algebraic twists of modular forms and Hecke orbits, 
\emph{GAFA} \textbf{25} (2015), 580--657.



\bibitem[FM07]{FM07}
\'E. Fouvry \& Ph. Michel,
Sur le changement de signe des sommes de Kloosterman,
\emph{Annals of Math.} \textbf{165} (2007), 675--715.

\bibitem[FR22]{FR22}
\'{E}. Fouvry \& M. Radziwiłł, 
Level of distribution of unbalanced convolutions,
\emph{Ann. Sci. \'Ec. Norm. Sup\'er.} \textbf{55} (2022), 537--568.




\bibitem[FI92]{FI92} 
J. B. Friedlander \& H. Iwaniec,
A mean-value theorem for character sums, 
\emph{Michigan Math. J.} \textbf{39} (1992), 153--159.

\bibitem[FI97]{FI97} 
J. B. Friedlander \& H. Iwaniec, 
The Brun--Titchmarsh theorem. Analytic Number Theory (Kyoto, 1996), 
London Math. Soc. Lecture Note Ser., 247, Cambridge Univ. Press, Cambridge, 1997, 85--93.


\bibitem[FI10]{FI10} 
J. B. Friedlander \& H. Iwaniec, 
Opera de Cribro, 
\emph{Amer. Math. Soc. Colloq. Publ.}, Vol. {\bf 57}, AMS, Providence, RI, 2010.



\bibitem[Go75]{Go75} 
D. Goldfeld, 
A further improvement of the Brun--Titchmarsh theorem,
\emph{J. London Math. Soc.} \textbf{11} (1975), 434--444.




\bibitem[GK91]{GK91} 
S. W. Graham \& G. Kolesnik, 
Van der Corput Method of Exponential Sums, 
London Mathematical Society Lecture Note Ser., 126, Cambridge, Cambridge University Press, 1991.


\bibitem[GR90]{GR90} 
S. W. Graham \& C. J. Ringrose, 
Lower bounds for least quadratic nonresidues, 
Analytic Number Theory (Allerton Park, IL, 1989), Progr. Math. 85, Birkh\"auser, Boston, 1990, 269--309.



\bibitem[Gr22]{Gr22} 
A. Granville,
Sieving intervals and Siegel zeros,
\emph{Acta Arith.} \textbf{205} (2022), 1--19.

\bibitem[HB78]{HB78} 
D. R. Heath-Brown, 
Hybrid bounds for Dirichlet $L$-functions, 
\emph{Invent. Math.} \textbf{47} (1978), 149--170.



\bibitem[Ho72]{Ho72} 
C. Hooley, 
On the Brun--Titchmarsh theorem, 
\emph{J. Reine Angew. Math.} \textbf{255} (1972), 60--79.



\bibitem[Ho75]{Ho75} 
C. Hooley, 
On the Brun--Titchmarsh theorem. II,  
\emph{Proc. London Math. Soc.} \textbf{30} (1975), 114--128. 


\bibitem[Ho78]{Ho78} 
C. Hooley, 
On the greatest prime factor of a cubic polynomial, 
\emph{J. Reine Angew. Math.} \textbf{303}/\textbf{304} (1978), 21--50.

\bibitem[Hu74]{Hu74}
M. N. Huxley,
Large values of Dirichlet polynomials, III,
\emph{Acta Arith.} \textbf{26} (1974/75), 435--444.

\bibitem[Iw80]{Iw80} 
H. Iwaniec, 
A new form of the error term in the linear sieve, \emph{Acta Arith.} \textbf{37} (1980), 307--320.

\bibitem[Iw82]{Iw82} 
H. Iwaniec, 
On the Brun--Titchmarsh theorem, 
\emph{J. Math. Soc. Japan} \textbf{34} (1982), 95--123.

\bibitem[IK04]{IK04} 
H. Iwaniec \& E. Kowalski, 
Analytic Number Theory, 
\emph{Amer. Math. Soc. Colloq. Publ.}, Vol. 53, AMS, Providence, RI, 2004.


\bibitem[Ju77]{Ju77}
M. Jutila,
Zero-density estimates for $L$-functions, 
\emph{Acta Arith.} \textbf{32} (1977), 55--62.


\bibitem[Ka95]{Ka95} 
A. A. Karatsuba, 
Fractional parts of functions of a special form, 
\emph{Izvestiya: Mathematics} \textbf{59} (1995), 721--740.


\bibitem[Ka88]{Ka88} 
N. M. Katz, Gauss Sums, Kloosterman Sums, And Monodromy Groups, 
Ann. of Math. Stud., Vol. 116, Princeton University Press, Princeton, NJ, 1988.


\bibitem[Ka90]{Ka90} 
N. M. Katz, 
Exponential Sums and Differential Equations, 
Ann. of Math. Stud., Vol. 124, Princeton University Press, Princeton, NJ, 1990.


\bibitem[KS03]{KS03}
H. Kim \& P. Sarnak,
Appendix: Refined estimates towards the Ramanujan and Selberg Conjectures,
\emph{J. Amer. Math. Soc.} \textbf{16} (2003), 175--181.


\bibitem[Kl61]{Kl61} 
N. I. Klimov,
Almost prime numbers,
\emph{Uspehi Mat. Nauk} \textbf{16} (1961), 181--188.



\bibitem[La87]{La87} 
G. Laumon, 
Transformation de Fourier, constantes d'\'{e}quations fonctionnelles et conjecture de Weil, 
\emph{Publ. Math. IH\'ES} \textbf{65} (1987), 131--210.

\bibitem[LM24]{LM24} 
Y. Lin \& Ph. Michel,
On algebraic twists with composite moduli, 
\emph{Ramanujan J.} \textbf{63} (2024), 803--837.

\bibitem[LR65]{LR65}
J. H. van Lint \& H.-E. Richert, 
On primes in arithmetic progressions, 
\emph{Acta Arith.} \textbf{11} (1965), 209--216.

\bibitem[Lo16]{Lo16}
B. L\"offel,
Algebraic twists of modular forms and sums over primes to squarefree moduli,
ETHZ thesis, 2016, \url{https://www.research-collection.ethz.ch/handle/20.500.11850/155985}.

\bibitem[Ma09]{Ma09}
K. Matom\"aki,
The distribution of $\alpha p$ modulo one,
\emph{Math. Proc. Camb. Phil. Soc.} \textbf{247} (2009), 267--283.



\bibitem[Ma13]{Ma13}
J. Maynard,
On the Brun--Titchmarsh theorem,
\emph{Acta Arith.} \textbf{157} (2013), 249--296.



\bibitem[Mo71]{Mo71}
H. L. Montgomery,
Topics in Multiplicative Number Theory,
Lecture Notes in Math. 227, Springer, Berlin, 1971. 


\bibitem[MV73]{MV73}
H. L. Montgomery \& R. C. Vaughan,
The large sieve,
\emph{Mathematika} \textbf{20} (1973), 119--134.

\bibitem[Mo73]{Mo73} 
Y. Motohashi, 
On some improvements of the Brun--Titchmarsh theorem. II (in Japanese), 
\emph{RIMS Kokyuroku} \textbf{193} (1973), 97--109. 


\bibitem[Mo74]{Mo74} 
Y. Motohashi, 
On some improvements of the Brun--Titchmarsh theorem, 
\emph{J. Math. Soc. Japan} \textbf{26} (1974), 306--323. 

\bibitem[Mo79]{Mo79}
Y. Motohashi,
A note on Siegel's zeros,
\emph{Proc. Japan Acad. Ser. A Math. Sci.} \textbf{55} (1979), 190--192.

\bibitem[Mo99]{Mo99}
Y. Motohashi,
On the remainder term in the Selberg sieve, Number Theory in Progress, Vol. 2 (Zakopane-Ko\'scielisko, 1997), 1053--1064,
Walter de Gruyter \& Co., Berlin, 1999.


\bibitem[Mo12]{Mo12} 
Y. Motohashi, 
On some improvements of the Brun--Titchmarsh theorem. IV, Preprint (2012), arXiv:1201.3134 [math.NT].


\bibitem[Po14]{Po14} 
D. H. L. Polymath, 
New equidistribution estimates of Zhang type, 
\emph{Algebra Number Theory} \textbf{8} (2014), 2067--2199.


\bibitem[Se91]{Se91} 
A. Selberg, 
Collected papers. Vol. II.
Springer-Verlag, Berlin, 1991. 


\bibitem[Sh80]{Sh80}
P. Shiu,
A Brun--Titchmarsh theorem for multiplicative functions,
\emph{J. Reine Angew. Math.} \textbf{313} (1980), 161--170.



\bibitem[Si83]{Si83} 
H. Siebert, 
Sieve methods and Siegel's zeros,
Studies in Pure Mathematics, To the Memory Paul Tur\'an, Birkh\"auser, Basel, 1983, 659--668. 


\bibitem[Ti30]{Ti30}
E. C. Titchmarsh, 
A divisor problem,
{\em Rend. Circ. Math. Palermo} \textbf{54} (1930), 414--429. Correction: ibid. \textbf{57} (1933), 478--479.


\bibitem[WX21]{WX21}
J. Wu \& P. Xi, 
Arithmetic exponent pairs for algebraic trace functions and applications, with an appendix by W. Sawin,
\emph{Algebra Number Theory} \textbf{15} (2021), 2123--2172.

\bibitem[Xi17]{Xi17}
P. Xi,
Gaussian distributions of Kloosterman sums: vertical and horizontal,
\emph{Ramanujan J.} \textbf{43} (2017), 493--511.

\bibitem[XZ25]{XZ25}
P. Xi \& J. Zheng,
On the Brun--Titchmarsh theorem. II, \emph{IMRN}, Issue 21, November 2025, rnaf33.



\end{thebibliography}

\end{document}